\documentclass[leqno]{amsart}
\usepackage{egastyle}
\usepackage{mymacros}
\usepackage[all,tips]{xy}
\usepackage{url}
\usepackage{bm}
\usepackage{verbatim}
\usepackage{epsfig}
\usepackage{graphicx}

\CompileMatrices

\input{paper1-export.aux}

\newwrite\exportaux
\gdef\exportauxname{paper2-export.aux}
\immediate\openout\exportaux\exportauxname

\makeatletter
\let\oldlabel\label
\def\label#1{\@bsphack\immediate\write\exportaux%
  {\string\newlabel{#1}{{II.\@currentlabel}{\thepage}}}\@esphack%
  \oldlabel{#1}}
\makeatother

\def\cC{{\mathcal C}}

\def\P{\bP}

\newcommand{\PP}{\mathbf{P}}

\newcommand{\RR}{\mathbf{R}}

\newcommand{\T}{\mathbf{T}}

\DeclareMathOperator\val{val}

\renewcommand\red{\operatorname{red}}

\newcommand\an{{\operatorname{an}}}

\renewcommand\tilde\widetilde

\title[Lifting harmonic morphisms II]{Lifting harmonic morphisms II:
  tropical curves and metrized complexes}

\author{Omid Amini} 
\email{oamini@math.ens.fr}
\address{CNRS-DMA, \'Ecole Normale Sup\'erieure, 45 Rue d'Ulm, Paris}

\author{Matthew Baker} 
\email{mbaker@math.gatech.edu}
\address{School of Mathematics, Georgia Institute of Technology, Atlanta GA 30332-0160, USA}

\author{Erwan Brugall\'e} 
\email{brugalle@math.jussieu.fr}
\address{Universit\'e Pierre et Marie Curie, Paris 6, 4 Place Jussieu, 75 005 Paris, France}

\author{Joseph Rabinoff}
\email{rabinoff@math.gatech.edu}
\address{School of Mathematics, Georgia Institute of Technology, Atlanta GA 30332-0160, USA}

\begin{document}

\begin{abstract}  
  In this paper we prove several lifting theorems for
  morphisms of tropical curves.  We interpret the obstruction to lifting a
  finite harmonic morphism of
  augmented metric graphs to a morphism of algebraic curves as the
  non-vanishing of certain Hurwitz
  numbers, and we give various conditions under which this obstruction
  does vanish.  In particular we show that any finite harmonic morphism of
  (non-augmented) metric graphs lifts.  
  We also give various applications of these results.  For example, we
  show that linear equivalence of divisors on a tropical curve $C$
  coincides with the equivalence relation generated by declaring that the
  fibers of every finite harmonic morphism from $C$ to the tropical
  projective line are equivalent.  We study liftability of metrized
  complexes equipped with a finite group action, and use this to
  classify all augmented metric graphs arising as the tropicalization of
  a hyperelliptic curve.  We prove that there exists a $d$-gonal tropical
  curve that does not lift to a $d$-gonal algebraic curve.

  This article is the second in a series of two.
\end{abstract}

\thanks{We are grateful to Andrew Obus for a number of useful comments based on a careful reading of the first arXiv version of this manuscript. We thank Ye Luo for allowing us to include Example~\ref{ex:luo}.
M.B. was partially supported by NSF grant DMS-1201473.
E.B. was partially supported by the ANR-09-BLAN-0039-01.}

\maketitle

{\small Throughout this paper, unless explicitly stated otherwise,
$K$ denotes a complete algebraically closed non-Archimedean field with nontrivial valuation 
$\val:K\to\R\cup\{\infty\}$.  Its valuation ring is denoted $R$, its
maximal ideal is $\fm_R$, and the residue field is $k = R/\fm_R$.  We
denote the value group of $K$ by $\Lambda = \val(K^\times)\subset\R$.}

\section{Introduction}

This article is the second in a series of two.  The first, entitled
\emph{Lifting harmonic morphisms I: metrized complexes and Berkovich skeleta}, 
will be cited as~\cite{abbr:lifting1}; references of the form
``Theorem~I.1.1'' will refer to Theorem~1.1 in~\cite{abbr:lifting1}. 

\paragraph
The basic motivation behind the investigations in this paper is to
understand the relationship between tropical and algebraic curves. 
A fundamental problem along these lines is to understand 
which morphisms between tropical curves arise as tropicalizations\footnote{In the present paper tropicalization is 
defined via Berkovich's theory of analytic spaces (see also \cite{payne:analytification}, \cite{bpr:trop_curves},~\cite{CD12}). Another framework for tropicalization has been proposed by Kontsevich-Soibelman~\cite{kontsevich_soibelman:SYZ} and Mikhalkin~(see for example~\cite{mikhalkin:tropical_geometry}), 
where the link between tropical geometry and complex algebraic geometry is provided by
 real one-parameter families of complex varieties. For some conjectural
 relations between the two approaches see~\cite{kontsevich_soibelman:SYZ, kontsevich_soibelman:affinestructures}.} of morphisms of 
algebraic curves.
More precisely, we are interested in the following  question:

\begin{itemize}
\item[(Q)]
Given a curve $X$ with tropicalization $C$, can we classify the branched covers of $X$ in terms of (a suitable notion of) branched covers of $C$?
\end{itemize}

In addition to this lifting problem for morphisms of tropical curves, we also study 
questions such as ``Which tropical curves arise as tropicalizations of hyperelliptic curves?''.
This naturally leads us to study group actions on tropical curves and how notions such as gonality change under tropicalization.  

\smallskip

In this paper we will consider three different kinds of ``tropical'' objects which one can associate to 
a smooth, proper, connected algebraic curve $X/K$, each depending on the choice of a \emph{triangulation} of $X$.
Roughly speaking, a triangulation $(X,V\cup D)$ of $X$ (with respect to a finite set of punctures $D\subset X(K)$)
is a finite set $V$ of points in the Berkovich analytification $X^\an$ of $X$ whose removal partitions $X^\an$ into open balls
and finitely many open annuli (with the punctures belonging to distinct open balls).
Triangulations of $(X,D)$ are naturally in
one-to-one correspondence with semistable models $\fX$ of
$(X,D)$; see Section~\ref{sec:simultaneous.ss.reduction}.  
The skeleton of a triangulated curve is the dual graph of the special fiber $\fX_k$ of the corresponding semistable
model, with infinite rays for the punctures, equipped with its canonical metric.

\smallskip

To any  triangulated curve, one may associate the three following
``tropical'' objects, at each  step adding some additional structure:
\begin{enumerate}
\item a \emph{metric graph} $\Gamma$: this is the skeleton of the
  triangulated curve $(X,V\cup D)$;
\item an \emph{augmented metric graph} $(\Gamma,g)$, i.e., a metric
  graph $\Gamma$ enhanced with a genus
  function $g:\Gamma\to \Z_{\ge 0}$ which is non-zero only at finitely
  many points: this is the above metric graph together with the function $g$ satisfying $g(p)=0$ for $p \not\in V$
  and $g(p)={\rm genus}(C_p)$ for $p \in V$, where $C_p$ is the (normalization of the) irreducible component of $\fX_k$
  corresponding to $p$;
\item a \emph{metrized complex of curves} $\cC$, i.e., an augmented
  metric graph $\Gamma$ equipped with a vertex set $V$ and a punctured algebraic curve over $k$
  of genus $g(p)$ for each point $p\in V$, with the punctures in bijection with the tangent directions to $p$ in $\Gamma$: 
  this is the above metric graph, together with the curves $C_p$ for $p\in V$ and punctures given by the singular points of
  $\fX_k$.
\end{enumerate}

\smallskip

An (augmented) metric graph or metrized complex
of curves arising from a triangulated curve by the above procedure is
said to be \emph{liftable}.
If $(X,V\cup D)$ and $(X,V'\cup D')$ are triangulations of the same curve $X$, with $D' \subset D$ and $V' \subset V$,
then the corresponding metric graphs are related by a so-called \emph{tropical modification}.
Tropical modifications generate an equivalence relation on the set of
(augmented) metric graphs, and an equivalence class for this relation is called an \emph{(augmented)
  tropical curve}. The (augmented) {\em tropicalization} of a $K$-curve $X$ is by definition the
(augmented) tropical curve $C$ corresponding to any triangulation of $X$.
Tropical curves and augmented tropical curves can be thought of as ``purely combinatorial''
objects, whereas metrized complexes are a mixture of combinatorial objects (which one thinks of as living over the
value group $\Lambda$ of $K$) and algebro-geometric objects over the residue field $k$ of $K$. 

\smallskip

There is a natural notion of {\em finite harmonic morphism} between metric graphs which induces a natural 
notion of {\em tropical morphism} between tropical curves.
There is a corresponding notion of tropical morphism for {\em augmented} tropical curves, where in addition to the 
harmonicity condition one imposes a ``Riemann--Hurwitz condition'' that the ramification divisor is effective.
There is also a natural notion of finite harmonic morphism for metrized complexes of curves.
Each kind of object (metric graphs, tropical curves, augmented tropical curves, metrized complexes) forms
a category with respect to the corresponding notion of morphism.
 The  construction of an (augmented) tropical curve $C$ 
(resp.\ metrized complex  $\cC$) out of a  $K$-curve $X$ 
(resp.\ triangulated $K$-curve $(X,V\cup D)$)
is functorial, in the sense that a finite morphism of curves induces in a
natural way a  tropical morphism $C'\to C$
(resp.\ a finite harmonic morphism $\cC'\to \cC$).

\paragraph
Our original question (Q) now breaks up into the following two separate questions.   
\begin{itemize}
\item[(Q1)] Which finite harmonic morphisms $\cC' \to \cC$ of metrized complexes can be lifted to finite morphisms of triangulated curves (with a pre-specified
lift $X$ of $\cC$)?
\item[(Q2)] Which tropical morphisms between
  augmented tropical curves can be lifted to finite harmonic morphisms of 
metrized complexes?
\end{itemize}

One can also forget the augmentation function $g: \Gamma \to  \Z_{\ge 0}$ and ask the following variant of (Q2):
\begin{itemize}
\item[(Q$2^\prime$)] Which tropical morphisms between tropical curves can be lifted to finite harmonic morphisms of 
metrized complexes?
\end{itemize}

A consequence of the results of~\cite{abbr:lifting1} is that the answer to
question (Q1) is essentially ``all'', so the situation here is rather
satisfactory; there is no obstruction to lifting a finite harmonic
morphism $\cC' \to \cC$ to a branched cover of $X$, at least assuming
everywhere tame ramification when $k$ has characteristic $p>0$.
In particular, if $\chr(k)=0$ then there are no tameness issues, and we have the
following result:

\begin{thm*}
  Assume ${\rm char}(k)=0$ and let $\phi: \Sigma'\to\Sigma$ be a finite
  harmonic morphism of $\Lambda$-metrized complexes of $k$-curves.  Then
  there exists a finite morphism of triangulated punctured curves lifting
  $\phi$.
\end{thm*}

\smallskip

This follows immediately from Proposition~\ref{prop:lifting} and
Theorem~\ref{thm:graph.to.curve}.
We stress that the genus and degree are automatically preserved by
such lifts.

\smallskip

\begin{figure}[h]
\scalebox{.21}{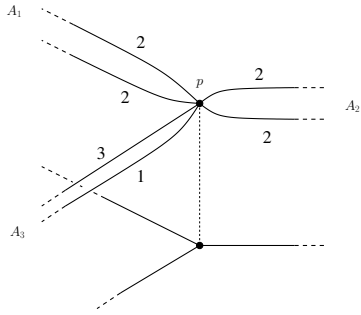}
\caption{A tropical morphism of degree four which cannot be promoted to a degree $4$ morphism of
metrized complexes of curves.  The labels on the edges are the ``expansion factors'' of the harmonic morphism.
See Definition \ref{def:morph metric}.}\label{fig:starmap}
\end{figure}

Essentially by definition, (Q2) reduces to an existence problem for
ramified coverings $\phi_{p'}:C'_{p'}\to C_p$ of a given degree with some
prescribed ramification profiles. Hence the answer to (Q2) is intimately
linked with the question of non-vanishing of Hurwitz numbers.  
See Proposition~\ref{prop:lift augmented to complex}. 
In particular one can easily construct tropical morphisms between augmented
tropical curves which cannot be promoted to a finite harmonic morphism of
metrized complexes (and hence cannot be lifted to a finite morphism of
smooth proper curves over $K$).  The simplest example of such a tropical
morphism is depicted in Figure~\ref{fig:starmap}, and corresponds to the
classical fact that although it would not violate the Riemann--Hurwitz
formula, there is no degree $4$ map of smooth proper connected curves over
${\mathbf C}$ having ramification profile $\{ (2,2),(2,2),(3,1) \}$; this
is a consequence of the (easy part of the) Riemann Existence Theorem (see
Example~\ref{rem:Hurwitznumbers} below for more details).

\smallskip

Understanding when Hurwitz numbers vanish remains mysterious in general,
so at present there is no satisfying ``combinatorial'' answer to question
(Q2), in which we require that the genus of the objects in question be
preserved by our lifts.  However, if we drop the latter condition, i.e.,
if we consider instead question (Q$2^\prime$), we will see that the answer
to (Q$2^\prime$) is also ``all'' (see Theorem~\ref{thm:lifting harm}):

\begin{thm*}
  Any finite harmonic morphism $\overline \phi:\Gamma'\to \Gamma$ of
  $\Lambda$-metric graphs is liftable if ${\rm char}(k)=0$.
\end{thm*}

\paragraph[Applications]
We prove a number of additional results which supplement and 
provide applications of the above results.  Some of these results are as
follows. 

\subparagraph[Tame group actions]
Let $\cC$ be a metrized complex and let $H$ be a finite subgroup of
$\mathrm{Aut}(\cC)$.    
We say the action of  $H$ on $\cC$ is {\it tame} if
for any vertex $p$ of $\Gamma$,  
the stabilizer group $H_p$ acts freely on a dense open subset of $C_p$,
and for any point $x$ of $C_p$, the stabilizer
subgroup $H_x$ of $H$ is cyclic of the form $\mathbf{Z} / d \mathbf{Z}$ 
for some integer $d$,  with $(d,p) =1 $ if $\chr(k) = p > 0$ (see Remark~\ref{rem:tameactions} for further explanation of this condition).
It follows from  Theorem \ref{thm:lifting2} (in its strong form, i.e.\
using the calculation of the automorphism group of a lift)
that we can lift $ \cC$
together with a tame action of $H$ 
 if and
only if the quotient $\cC/H$ exists
in the category of metrized complexes. We characterize when such a quotient exists
in Theorem \ref{thm:liftinggroupaction},  of which the following result is a special case:

\begin{thm*}
Suppose that the action of $H$ is tame and has
no isolated fixed points on the underlying metric graph of $\cC$.
 Then there exists a smooth, proper, and
geometrically connected algebraic $K$-curve $X$ lifting $\cC$ which is equipped with an
action of $H$ commuting with the tropicalization map.
\end{thm*}

\smallskip
In presence of isolated fixed points, there are additional hypothesis
on the action of $H$ to be liftable  to a $K$-curve.
As a concrete example, we prove the following
characterization of all augmented tropical curves arising
as the tropicalization of a hyperelliptic $K$-curve (see Corollary~\ref{cor:characterization hyperelliptic}):

\begin{thm*}
Let $\Gamma$ be an augmented metric graph of genus $g \geq 2$ having no infinite vertices or degree one vertices
of genus $0$.  Then there is a smooth proper hyperelliptic curve $X$ over $K$ of genus $g$ having $\Gamma$ as its minimal skeleton if and only if (a) there exists an involution $s$ on $\Gamma$ such that $s$ fixes all the points $p \in  \Gamma$ with $g(p) > 0$ and the quotient 
$ \Gamma/s$ is a metric tree, and (b) for every $p \in \Gamma$ the number of bridge edges adjacent to $p$
is at most $2g(p)+2$.
\end{thm*}

\subparagraph[Gonality of tropical curves] 
The \emph{tropical projective line} is the augmented tropical curve
$\T\P^1$ represented by any tree with genus function identically zero.
See Example \ref{def:rational tropical}.
An augmented tropical curve $C$
is called \emph{$d$-gonal} if there exists a tropical morphism of degree
$d$ from $C$ to $\T\P^1$.  By Corollary~\ref{cor:morphism.to.harmonic},
the gonality of an augmented tropical curve is always a lower bound for
the gonality of any lift to a smooth proper curve over $K$.  (See
Remark~\ref{rem:gonality.definitions} for a discussion of the various
notions of gonality of tropical curves existing in the literature.)
We prove in Section \ref{sec:examples} that none of the lower
bounds provided by tropical ranks and gonality
are sharp.  For example:

\begin{thm*}
\begin{enumerate}
\item There exists an augmented tropical curve $C$ of  gonality $4$ 
such that the
gonality of any lifting of $C$ is at least $5$. 
\item There exists an effective divisor $D$ on a tropical curve $C$ such that $D$ has
tropical rank equal to $1$, but any effective lifting of $D$ has rank $0$.
\end{enumerate}
\end{thm*}

\smallskip

The construction in~(1) uses the vanishing of the degree $4$ Hurwitz number
$H_{0,0}^4((2,2),(2,2),(3,1))$.
In fact we prove in Theorem \ref{thm:nonlift.gonality} a much stronger
statement: we 
exhibit an augmented (non-metric) graph $G$ such that \emph{none} of the augmented
tropical curves with $G$ as underlying augmented graph can be lifted to
a $4$-gonal $K$-curve.  This means that there is a finite graph with {\em stable gonality} $4$
(in the sense of \cite{cornelissen:graph_li-yau}) which is not the (augmented) dual graph of any $4$-gonal curve $X/K$.

The proof of (2) is based on our lifting results and an explicit example, due
to  Luo (see Example~\ref{ex:luo}), of  a degree $3$ and rank $1$ base-point free divisor $D$
on a tropical curve $C$ which does not appear as the fiber of any degree $3$ tropical morphism from $C$ to $\T \PP^1$.

\subparagraph[Linear equivalence of divisors] 
When the target curve has
genus zero, we investigate in~\parref{par:target genus zero} a variant of
question (Q$2^\prime$) in which 
the genus of the source curve may be prescribed, at the cost of losing
control over the degree of the morphism.  As an application, we show
in Theorem~\ref{thm:equ effective/linear} that
linear equivalence of divisors on a tropical curve $C$ coincides with the
equivalence relation generated by declaring that the fibers of every
tropical morphism from $C$ to the tropical projective line $\T\PP^1$ are
equivalent: 

\begin{thm*}
  Let $\Gamma$ be a metric graph.
  Linear equivalence of divisors on $\Gamma$ is the additive equivalence
  relation  generated by (the retraction to $\Gamma$ of) fibers of finite
  harmonic morphisms from a tropical modification of $\Gamma$ to a metric
  graph of genus zero. 
\end{thm*}

\paragraph[Organization of the paper]
The paper is organized as follows. 
Precise definitions of tropical modifications and tropical curves are
given in Section~\ref{sec:prelim}, along with various kinds of morphisms
between these objects.  In that section we also use 
results from~\cite{abbr:lifting1} to define tropicalizations of morphisms
of curves, and provide a number of
examples.   Lifting results for (augmented) metric graphs
and tropical curves are proved in 
Section~\ref{sec:lifting morphisms metric graphs}.  Section
\ref{sec:applications} contains applications of 
our lifting results.  For example, lifting results for
metrized complexes equipped with a finite group action are discussed
in~\parref{par:quotient}.  In~\parref{par:quotient} we also give a
complete classification of all hyperelliptic augmented tropical curves
which can be realized as the minimal skeleton of a hyperelliptic curve.
Finally, in Section \ref{sec:examples} we study tropical rank and gonality
and related lifting questions.

\paragraph[Related work] 
The definition of effective harmonic morphisms of augmented metric graphs
that we use is the same as in \cite{Br13}.
The closely related, but slightly different, notion of 
an ``indexed harmonic morphism'' between weighted graphs was considered in \cite{caporaso:gonality}.
The indexed pseudo-harmonic (resp.\ harmonic) morphisms
in \cite{caporaso:gonality} are closely related to
harmonic (resp.\ effective harmonic) morphisms in our sense when the
vertex sets are fixed (see Definition~\ref{def:morph metric}), and
non-degenerate morphisms in the sense of  \cite{caporaso:gonality} 
correspond to finite morphisms in our sense.
One notable difference is that in \cite{caporaso:gonality}, 
only the combinatorial type of the metric graphs are fixed;
the choice of positive indices in an indexed pseudo-harmonic morphism
determines the length of the edges in the source graph once the edge lengths in the target
are fixed.

\smallskip

Tropical modifications and the 
``up-to-tropical-modification'' point of view were introduced by 
Mikhalkin~\cite{mikhalkin:tropical_geometry}. 

\smallskip 

In~\parref{par:gonality} we propose a definition for the stable gonality of a graph which coincides with the one used
by Cornelissen et.\ al.\ in their recent preprint \cite{cornelissen:graph_li-yau}.
A slightly different notion of gonality for graphs was introduced by
Caporaso in \cite{caporaso:gonality}. 
We also define the gonality of an augmented tropical curve, 
which strikes us as a more natural and perhaps more useful notion
than the stable gonality of a graph (where the lengths of the edges in the source and target metric graphs are not pre-specified).
We emphasize the importance of considering the dual graph of the special fiber of a semistable model of a smooth proper $K$-curve as an (augmented) {\em metric graph} and
not just as a (vertex-weighted) graph.
Keeping track of the natural edge lengths allows us to avoid pathological examples like Example 2.18
in~\cite{caporaso:gonality} of a 3-gonal graph which is not divisorally 3-gonal.

\smallskip

The question of lifting effective harmonic morphisms of metric graphs also occurs naturally (in a related but different Archimedean framework) when one
considers degenerating families of complex algebraic dynamical systems; see for example \cite[Theorems 1.2 and 7.1]{DeMarco-McMullen}
where DeMarco and McMullen prove a lifting theorem for polynomial-like endomorphisms of (locally finite) simplicial trees which has applications to studying dynamical compactifications of the moduli space of degree $d$ polynomial maps.
Our Theorem~\ref{thm:polynomial like} was inspired by the results of DeMarco--McMullen.

\section{Algebraic and tropical curves} \label{sec:prelim}

In this section we introduce tropical curves and morphisms between them.
We use the results of~\cite{abbr:lifting1} to define functorial
``intrinsic tropicalizations'' of algebraic curves.  We will freely use
the definitions and notations in 
Section~\ref{sec:definitions}.  We reproduce some of them here for the
convenience of the reader. 

\paragraph[Metric graphs]
A \emph{$\Lambda$-metric graph} is a metric graph whose edge lengths
belong to $\Lambda$.
The length of an embedded segment $e$ in a metric graph $\Gamma$ is
denoted $\ell(e)$.  The set of tangent directions at a point $p$ of
$\Gamma$ is denoted $T_p(\Gamma)$.
To a harmonic morphism $\phi:\Gamma'\to\Gamma$ of
metric graphs are associated its degree $\deg(\phi)$, 
its degree at a point $d_{p'}(\phi)$, the degree of $\phi$ along
an edge (also called the expansion factor) $d_{e'}(\phi)$, the directional
derivative of $\phi$ along a tangent direction at a vertex $d_{v'}(\phi)$, 
and the induced map on tangent spaces
$d\phi(p')$ when $d_{p'}(\phi)\neq 0$.

The group of divisors on a metric graph $\Gamma$ is denoted
$\Div(\Gamma)$.  A harmonic morphism of metric graphs
$\phi:\Gamma'\to\Gamma$ gives rise to push-forward and pull-back
homomorphisms $\phi_*:\Div(\Gamma')\to\Div(\Gamma)$
and $\phi^*:\Div(\Gamma)\to\Div(\Gamma')$ defined by
\[ \phi^*(p) = \sum_{p'\mapsto p} d_{p'}(\phi)\,(p') \quad
\mbox{and}\quad  \phi_*(p') = (\phi(p')) \]
and extending linearly.  It is clear that for $D\in\Div(\Gamma)$ we have
$\deg(\phi^*(D)) = \deg(\phi)\cdot\deg(D)$ and
$\deg(\phi_*(D)) = \deg(D)$.

\paragraph[Augmented metric graphs]
An augmented metric graph $\Gamma$ has a genus function
$g:\Gamma\to\Z_{\geq0}$.  We say that $\Gamma$ is 
\emph{totally degenerate} provided that $g$ is identically zero.  The
genus of $\Gamma$ is 
\[ g(\Gamma) = h_1(\Gamma) + \sum_{p\in \Gamma} g(p), \]
where $h_1(\Gamma)$ is the first Betti number of $\Gamma$.  
If $g(\Gamma) = 0$ then we say that $\Gamma$ is \emph{rational}.
The canonical divisor of an augmented metric graph $\Gamma$ is 
\[ K_\Gamma = \sum_{p\in \Gamma} (\val(p) + 2g(p) - 2)\, (p). \]
The degree of $K_\Gamma$ is $\deg(K_\Gamma) = 2g(\Gamma)-2$.

Let  $\phi:\Gamma'\to\Gamma$ be a harmonic morphism of
augmented metric graphs.
The \emph{ramification divisor} of $\phi$ is the divisor 
$R = \sum R_{p'} (p')$, 
where  for $p'\in \Gamma'$,
\[ R_{p'} = d_{p'}(\phi)\cdot\big(2-2g(\phi(p'))\big) - \big(2-2g(p')\big)
  - \sum_{v' \in T_{p'}(\Gamma')} \big(d_{v'}(\phi)-1\big). 
\]
We have the Riemann--Hurwitz formula
\[ K_{\Gamma'} = \phi^*(K_\Gamma) + R. \]
We say that $\phi$ is \emph{generically \'etale} if $R$ is supported on the set
of infinite vertices of $\Gamma$ and is \emph{\'etale} if $R = 0$.

\paragraph[Effective harmonic morphisms]
We will use the following Riemann--Hurwitz condition 
in formulating lifting 
problems for harmonic morphisms of augmented metric graphs.
Given a vertex $p'\in V(\Gamma')$ with $d_{p'}(\phi)\ne 0$, we define the
\emph{ramification degree} of $\phi$ at $p'$ to be 
\[ r_{p'}=R_{p'} - \#\{v'\in T_{p'}(\Gamma') ~:~ d_{v'}(\phi)=0\}. \]

Clearly $r_{p'}\le R_{p'}$, with $r_{p'}=R_{p'}$ if and only if 
$d_{v'}(\phi)>0$ for any $v'\in T_{p'}(\Gamma')$, i.e.
the distinction between ramification divisors and ramification degrees
only makes sense for non-finite harmonic morphisms. 
Our motivation  not to restrict ourselves to finite harmonic
morphisms is that  non-finite harmonic morphisms show up naturally in
many practical situations.

\begin{defn}
A harmonic morphism of augmented $\Lambda$-metric
graphs  $\phi:\Gamma'\to\Gamma$ is said to be \emph{effective} if 
$r_{p'}\ge 0$ for every finite vertex $p'$ of $\Gamma'$ with
$d_{p'}(\phi)\ne 0$.
\end{defn}

The significance
of the number $r_{p'}$ is given in Remark \ref{rem:RH}. In particular, 
only effective  harmonic morphisms of augmented metric graphs have a
chance to be liftable to a harmonic morphism of
  metrized complexes of curves, and possibly to a morphism of
  triangulated punctured $K$-curves. See Remark \ref{rem:lift effecitve}.

\smallskip 

Note that a generically \'etale morphism of augmented metric
graphs is effective.

\begin{eg} \label{eg:effective.examples}
Consider the harmonic morphisms of graphs $\phi:\Gamma'\to\Gamma$
represented in Figure~\ref{fig:ex morph augmented}.
We use the following conventions in our pictures: black dots represent
vertices of $\Gamma'$ 
and $\Gamma$; we label an edge with its degree if and only if the degree is different from $0$ and $1$;
we do not specify the lengths of edges of $\Gamma'$ and $\Gamma$.

The morphisms in Figure~\ref{fig:ex morph augmented}(a,b,c) are
effective provided that all the target graphs are totally degenerate.
Suppose that all $1$-valent vertices are infinite vertices in Figure 
\ref{fig:ex morph augmented}~(d,e), and that $g(p)=0$ in Figure
 \ref{fig:ex morph augmented}(e) and
$g(p)=1$ in Figure \ref{fig:ex morph augmented}(e).
Then $r_{p'}=2g(p')-1$ and $r_{p'_i}=2g(p'_i)-2$, so the morphism depicted
in  \ref{fig:ex morph augmented}(d) is effective if and only if 
$g(p')\ge 1$, and the morphism depicted
in  \ref{fig:ex morph augmented}(e) is effective if and only if both
vertices $p'_1$ and $p'_2$ have genus at least $1$.

The morphism in Figure~\ref{fig:starmap} is effective when
both graphs are totally degenerate.

\noindent
\begin{figure}[h]
\begin{tabular}{ccccc}
\scalebox{.31}{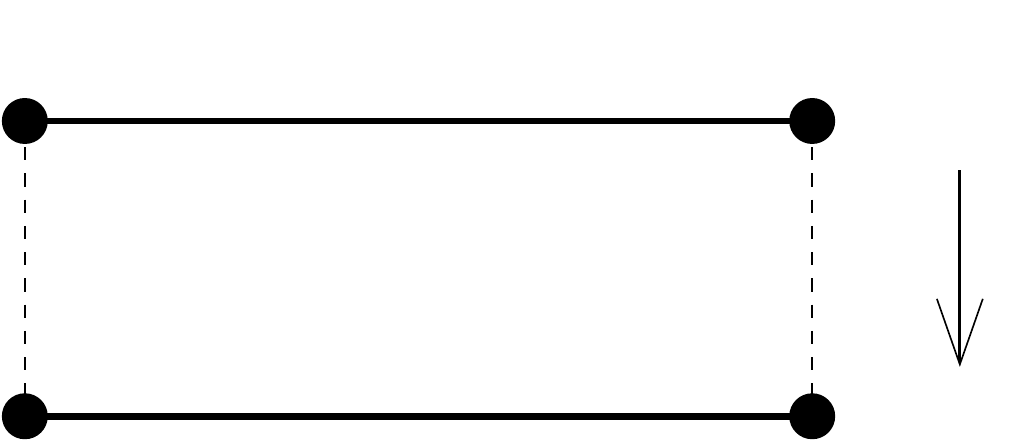}&  \hspace{3ex}&
\scalebox{.31}{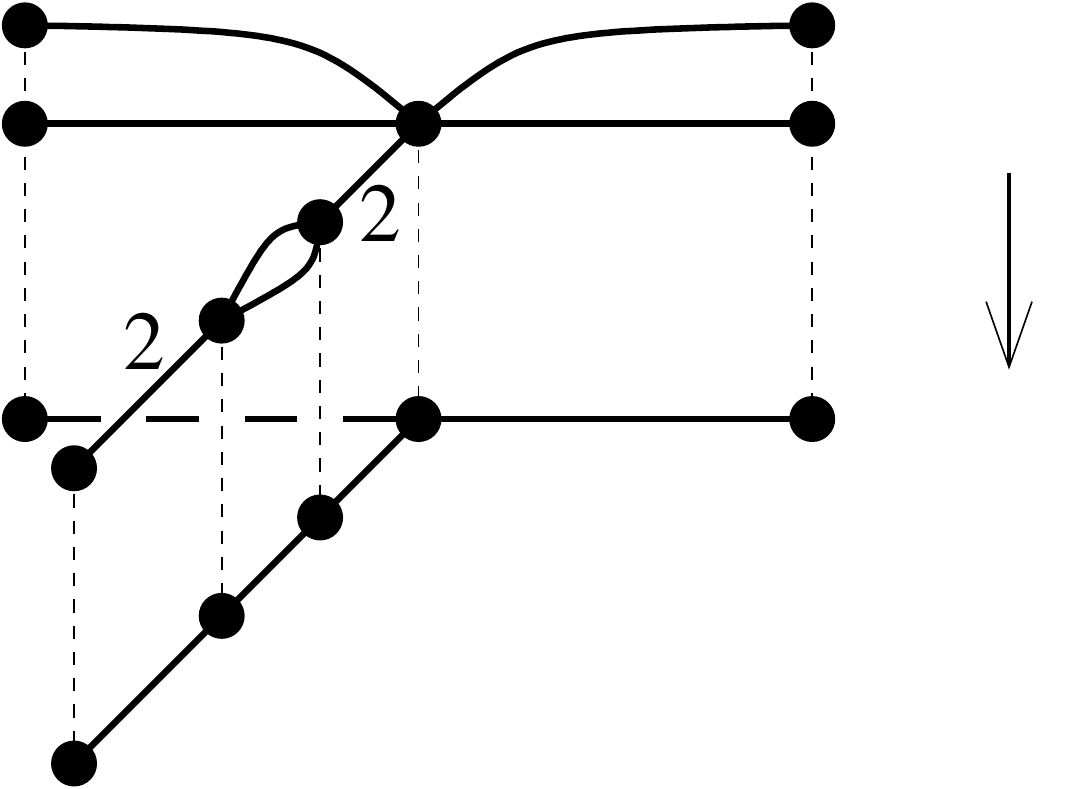}&  \hspace{3ex}&
\scalebox{.31}{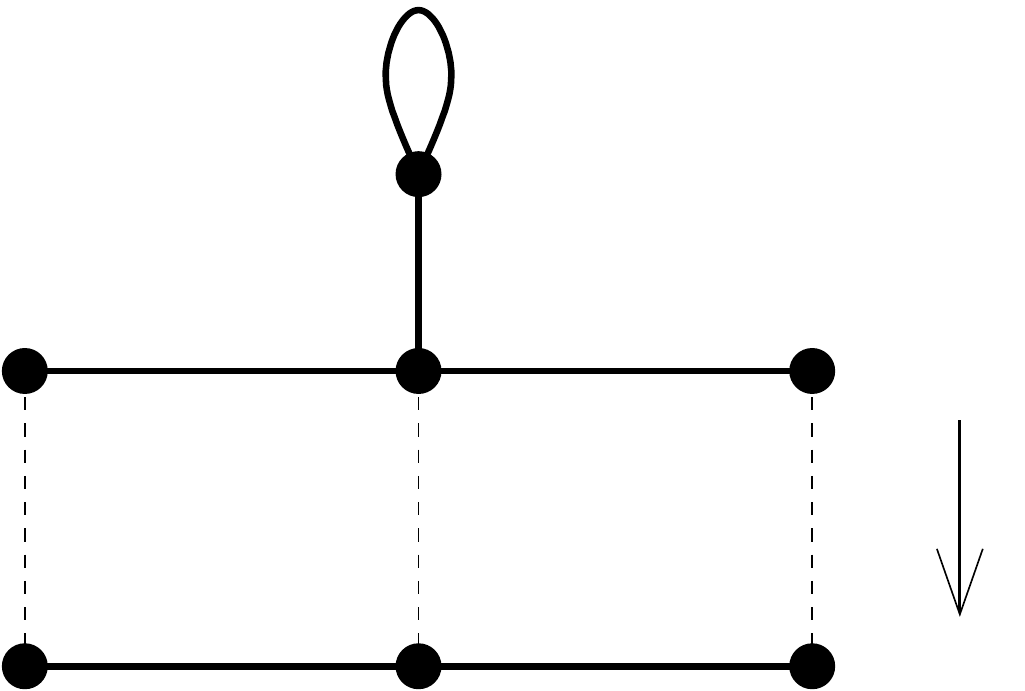}
\\ \\a)   && b) && c) 
\end{tabular}

\smallskip

\noindent
\begin{tabular}{ccc}
\scalebox{.31}{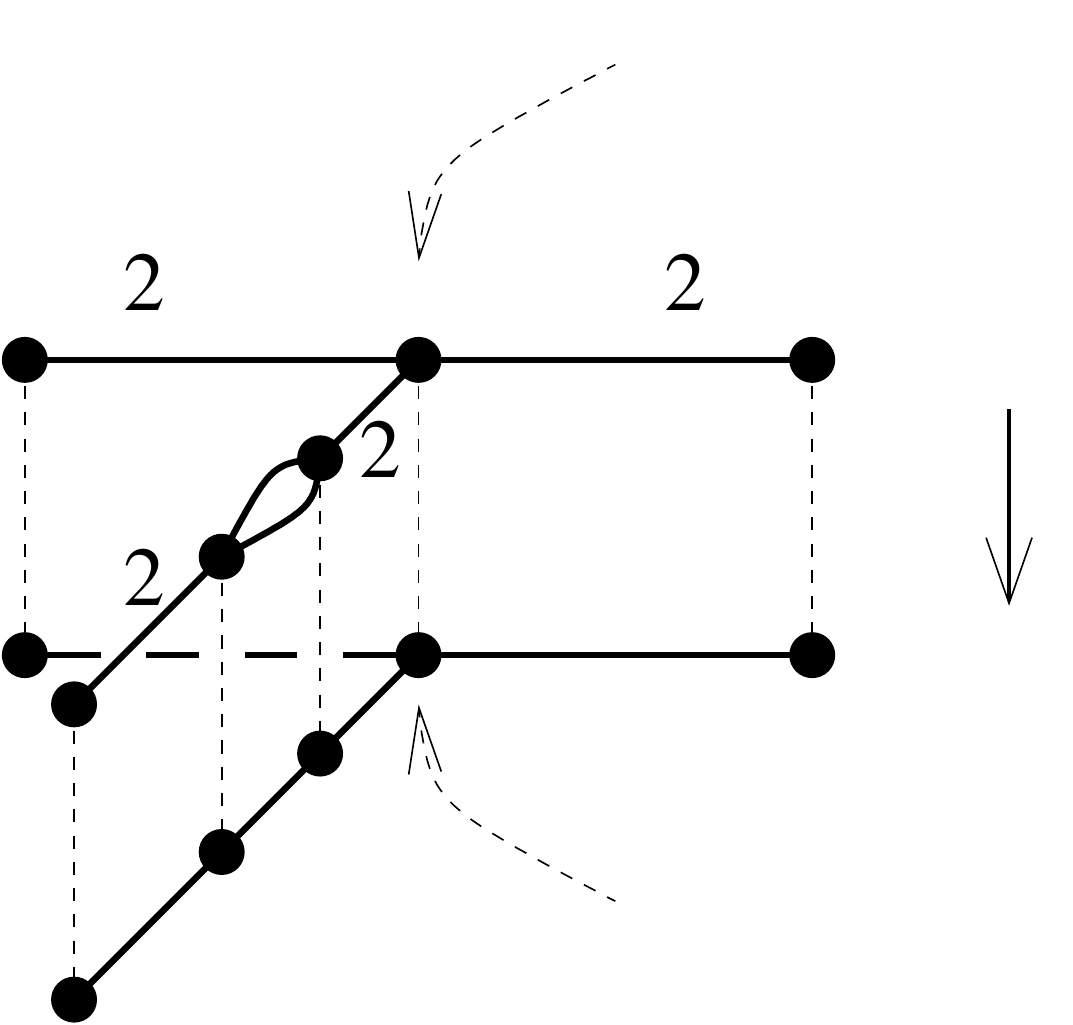}&  \hspace{5ex} &
\scalebox{.31}{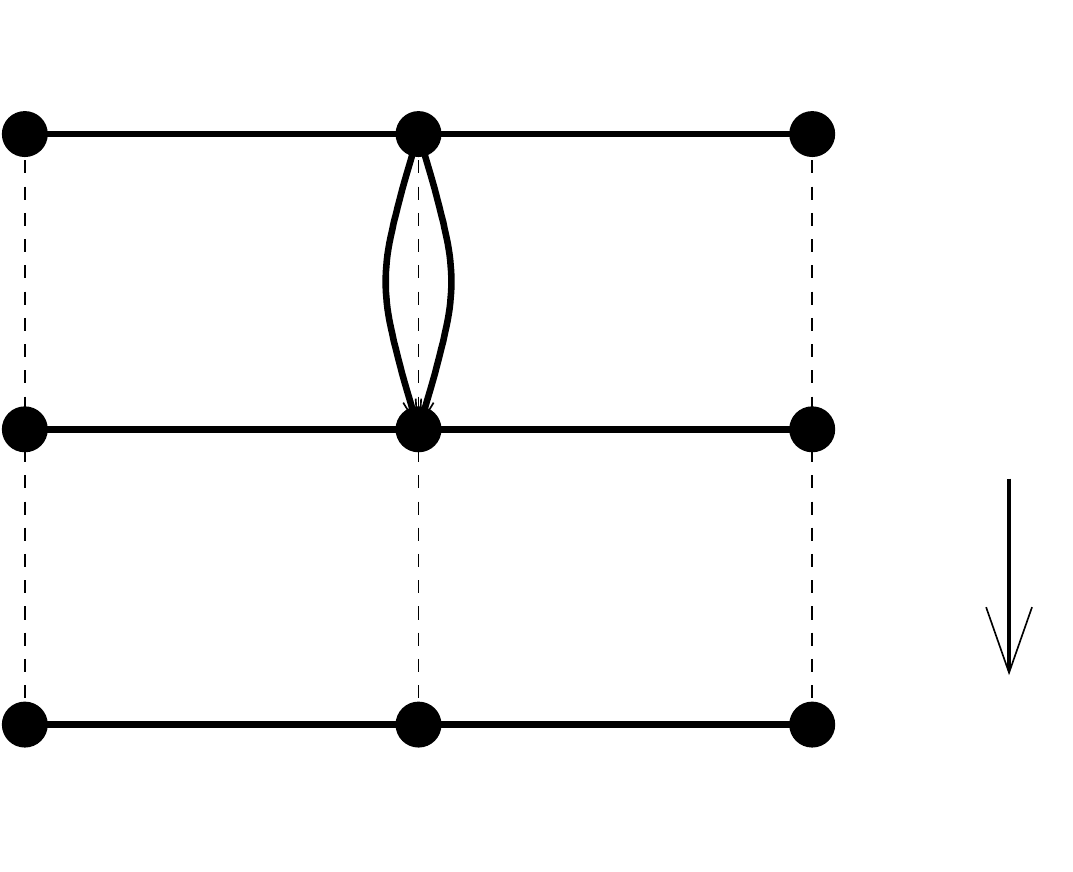}
\\ \\d) $g(p)=0$   && e) $g(p)=1$
\end{tabular}
\caption{}\label{fig:ex morph augmented}
\end{figure}
\end{eg}

\paragraph[Metrized complexes of curves]
Metrized complexes of curves and harmonic morphisms between them are
defined in~\parref{par:metrized.complex.defn}.  We recall part of the
definitions here.  A \emph{$\Lambda$-metrized complex of
$k$-curves} $\cC$ is the data of an underlying augmented $\Lambda$-metric graph
$\Gamma$ with a distinguished vertex set, and for each finite vertex
$p\in\Gamma$ a smooth proper connected $k$-curve $C_p$ of genus $g(p)$,
called the  \emph{residue curve}, and an injective reduction map
$\red_p:T_p(\Gamma)\inject C_p(k)$.  A \emph{harmonic morphism}
$\phi:\cC'\to\cC$ is a harmonic morphism of underlying augmented metric
graphs $\phi:\Gamma'\to\Gamma$, taking finite vertices of $\Gamma'$ to
finite vertices of $\Gamma$, along with a finite morphism 
$\phi:C_{p'}\to C_{\phi(p')}$ for every finite vertex $p'$ of $\Gamma'$
such that $d_{p'}(\phi)\neq 0$,
satisfying the following compatibility conditions:
\begin{enumerate}
\item For every finite vertex $p'\in V(\Gamma')$ 
  and every tangent direction $v'\in T_{p'}(\Gamma')$
  with $d_{v'}(\phi)>0$, we have 
  $\phi_{p'}(\red_{p'}(v')) = \red_{\phi(p')}(d\phi(p')(v'))$, and
  the ramification degree of $\phi_{p'}$ at
  $\red_{p'}(v')$ is equal to $d_{v'}(\phi)$.

\item For every finite vertex $p'\in V(\Gamma')$ with 
  $d_{p'}(\phi)>0$, every tangent direction
  $v\in T_{\phi(p')}(\Gamma)$,
  and every point 
  $x'\in\phi_{p'}\inv(\red_{\phi(p')}(v))\subset C'_{p'}(k)$,
  there exists $v'\in T_{p'}(\Gamma')$ such that
  $\red_{p'}(v') = x'$.
  
\item For every finite vertex $p'\in V(\Gamma')$ with 
  $d_{p'}(\phi)>0$ we have
  $d_{p'}(\phi) = \deg(\phi_{p'})$.
\end{enumerate}
Let $\phi: \cC' \to \cC$ be a finite harmonic
morphism of metrized complexes of curves.  We say that $\phi$ is a
\emph{tame harmonic morphism} if $\phi_{p'}$ is tamely ramified for all
finite vertices $p'\in \Gamma'$.  
We call $\phi$ a {\em tame covering} if in addition it is
a generically \'etale finite morphism of augmented metric graphs. 

\begin{rem}\label{rem:RH} 
  It follows from the Riemann--Hurwitz formula applied to the maps
  $\phi_{p'}:C'_{p'}\to C_{\phi(p')}$ that a harmonic morphism of metrized
  complexes of curves gives rise to an \emph{effective} harmonic morphism of
  augmented metric graphs when each $\phi_{p'}$ is a separable morphism of
  curves;
the integer $r_{p'}$ is then the sum of ramification indices over
all ramification points of $\phi_{p'}$ not contained in
$\red_{p'}(T_{p'}(\Gamma')) $.  
  In particular, tame harmonic morphisms of metrized complexes of
  curves give rise to effective harmonic morphisms of augmented metric
  graphs. 
\end{rem}

\paragraph[Triangulated punctured curves and skeleta]
Let $X$ be a smooth, connected, proper algebraic $K$-curve
and let $D\subset X(K)$ be a finite set of \emph{punctures}. 
Recall from Definitions~\ref{def:semistable vertex set}
and~\ref{def:triangulated curve} that 
a \emph{semistable vertex set of $(X,D)$} is a finite set $V$ of
type-$2$ points of $X^\an$ such that $X^\an\setminus(V\cup D)$ is a 
disjoint union of open balls and finitely many once-punctured open balls
and open annuli. If $V$ is a semistable vertex set of $(X,D)$, 
then $(X,V\cup D)$ is called a \emph{triangulated punctured curve}.
The semistable vertex sets of $(X,D)$ are in bijective correspondence with
the semistable $R$-models of $(X,D)$.   

To a triangulated punctured curve $(X,V\cup D)$ one associates a canonical
$\Lambda$-metrized complex of curves $\Sigma(X,V\cup D)$ called its
\emph{skeleton}.  The genus of the underlying augmented metric graph
$\Gamma$ is equal to genus $g(X)$ of $X$. 
There is a canonical closed embedding $\Gamma\inject X^\an$ and a
retraction map $\tau:X^\an\to\Gamma$.  

A \emph{finite morphism} of triangulated punctured $K$-curves
$\phi:(X',V'\cup D')\to(X,V\cup D)$ consists of a finite morphsim
$\phi:X'\to X$ such that $\phi\inv(V) = V'$, 
$\phi\inv(D) = D'$, and 
$\phi\inv(\Sigma(X,V\cup D)) = \Sigma(X',V'\cup D')$ as sets.
Here we restate Corollary~\ref{cor:morphism.to.harmonic}:

\begin{prop*}
  Let $\phi:(X',V'\cup D')\to(X,V\cup D)$ be a finite morphism of
  triangulated punctured curves.  Then $\phi$ naturally
  induces a finite harmonic morphism of $\Lambda$-metrized complexes of
  curves 
  \[ \Sigma(X',\, V'\cup D') \To \Sigma(X,\, V\cup D). \]
\end{prop*}

\begin{defn}
  A finite harmonic morphism $\bar\phi:\Gamma'\to\Gamma$ of metrized complexes
  of curves (resp.\ augmented metric graphs, resp.\ metric graphs) is said
  to be \emph{liftable} provided that there exists a finite morphism of
  triangulated punctured $K$-curves 
  $\phi:(X',V'\cup D')\to(X,V\cup D)$ and an isomorphism of $\bar\phi$ with the induced
  finite harmonic morphism of skeleta $\Sigma(X',V'\cup D')\to\Sigma(X,V\cup D)$ 
  (resp.\ of augmented metric graphs underlying the skeleta, 
  resp.\ of metric graphs underlying the skeleta).
\end{defn}

\begin{rem}\label{rem:lift effecitve}
  Among all finite harmonic morphisms of augmented metric graphs, only the
  effective ones have a chance to be liftable to a finite morphism of
  triangulated punctured $K$-curves.  Since
  the induced morphism of skeleta is a finite harmonic morphism of
  metrized complexes of curves, this follows from Remark~\ref{rem:RH}.
\end{rem}

\paragraph[Tropical modifications and tropical curves]
\label{par:def trop curves}
Here we introduce an equivalence relation among metric
graphs; an equivalence class for this relation will be called a \emph{tropical curve}.

\begin{defn}\label{defn:modification}
An \emph{elementary tropical modification} of a $\Lambda$-metric graph $\Gamma_0$ is a
$\Lambda$-metric graph $\Gamma = [0,+\infty] \cup \Gamma_0$ obtained from $\Gamma_0$ by  
attaching the segment $[0,+\infty]$ to $\Gamma_0$ in such a way that $0 \in [0,+\infty]$ gets identified with 
a finite $\Lambda$-point $p \in \Gamma_0$. 
If $\Gamma_0$ is augmented, then
$\Gamma$ naturally inherits a genus function from $\Gamma_0$ by
declaring that every point of $(0,+\infty]$ has genus $0$.

An (augmented) $\Lambda$-metric graph $\Gamma$ obtained from an (augmented) $\Lambda$-metric
graph $\Gamma_0$ by a finite sequence of elementary tropical modifications is called a {\em tropical 
modification} of $\Gamma_0$.
\end{defn}

\smallskip

If $\Gamma$ is a tropical modification of $\Gamma_0$, then there is a natural retraction map
$\tau: \Gamma \rightarrow \Gamma_0$ which is the identity on $\Gamma_0$ and contracts
each connected component of $\Gamma \setminus \Gamma_0$ to the unique point in $\Gamma_0$ lying in the
topological closure of that component.    
The map $\tau$ is a (non-finite) harmonic morphism of (augmented) metric graphs. 

\smallskip

\begin{eg}
We depict an elementary tropical modification in Figure \ref{fig:ex modification}(a), 
and a tropical modification which is the sequence of two elementary tropical
modifications in  Figure \ref{fig:ex modification}(b).
\begin{figure}[h]
\begin{tabular}{ccc}

\scalebox{.32}{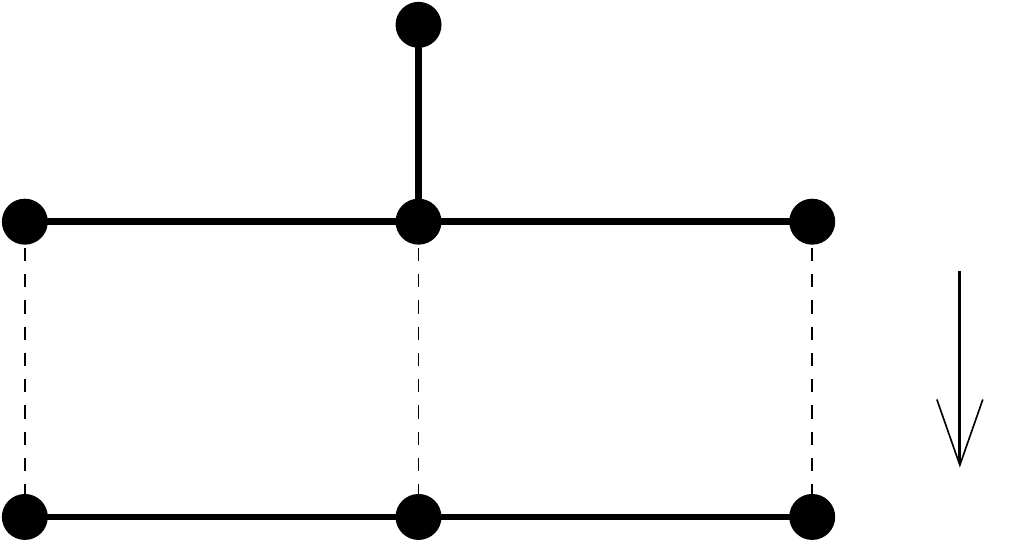}&  \hspace{5ex} &
\scalebox{.32}{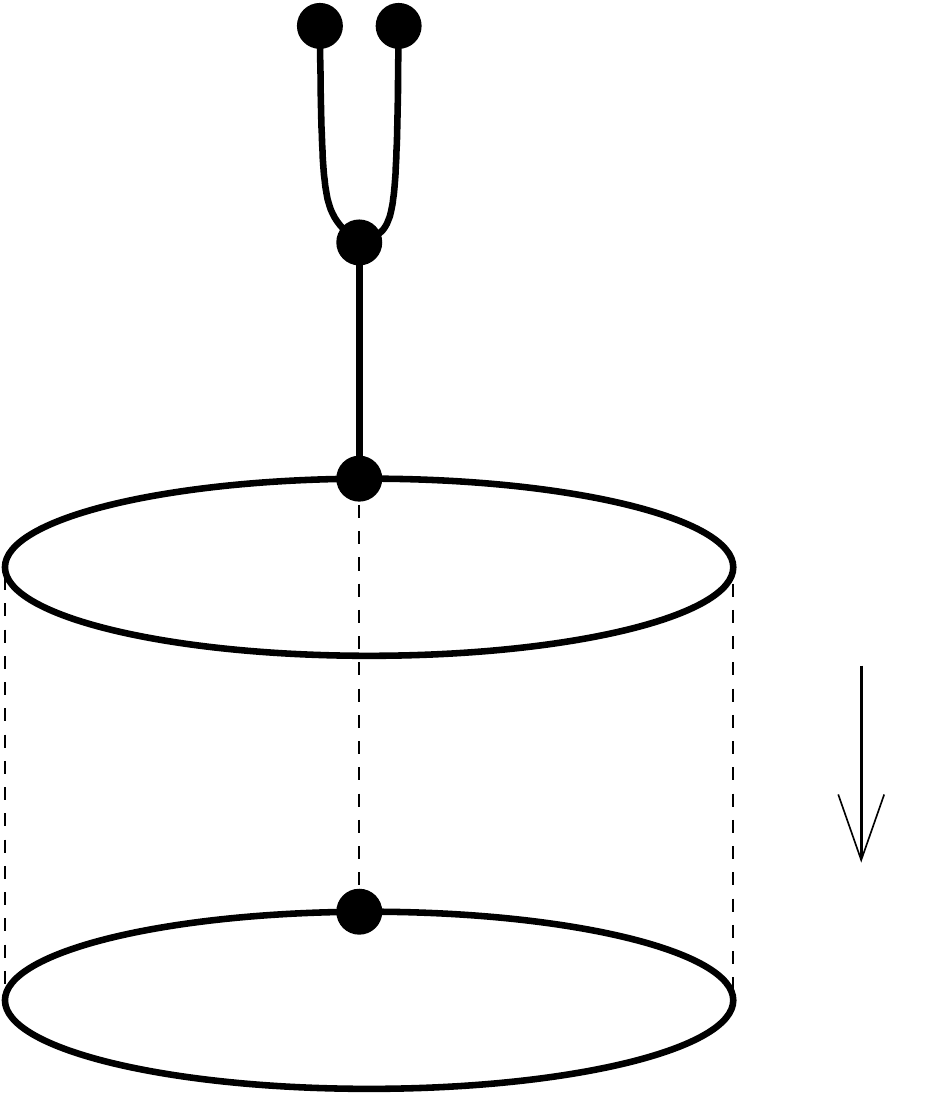}
\\ \\a) && b) 

\end{tabular}
\caption{Two tropical modifications}\label{fig:ex modification}
\end{figure}
\end{eg}

Tropical modifications generate an equivalence relation $\sim$ on the set of
(augmented) $\Lambda$-metric
graphs.
\begin{defn}
A \emph{$\Lambda$-tropical curve} (resp.\ an \emph{augmented $\Lambda$-tropical curve})
is an equivalence class of $\Lambda$-metric graphs (resp.\ augmented $\Lambda$-metric
graphs) with respect to $\sim$.
\end{defn}

\smallskip
In other words, a $\Lambda$-tropical curve is a 
$\Lambda$-metric graph considered up to tropical modifications and
their inverses (and similarly for augmented tropical curves).
By abuse of terminology, we will often refer to a tropical curve in terms of one
of its metric graph representatives.

\begin{eg}\label{def:rational tropical}
There exists a unique rational (augmented) tropical curve, which we denote by
$\T\PP^1$. 
Any rational (augmented) metric graph whose $1$-valent vertices are
all infinite is obtained by a sequence of tropical modifications from
the metric graph consisting of a unique finite vertex (of genus $0$).
\end{eg}

\begin{eg}\label{ex:add finite edge}
Let $\Gamma_0$ be a $\Lambda$-metric graph,  $p\in\Gamma_0$ a finite
$\Lambda$-point, and $l\in\Lambda\setminus\{0\}$.
We can construct a new $\Lambda$-metric graph $\Gamma$
by attaching the segment $[0,l]$ to $\Gamma_0$ via the identification of
 $0 \in [0,l]$  with $p$. 
Then $\Gamma_0$ and $\Gamma$ represent the same tropical curve, since the
elementary tropical modification of $\Gamma_0$ at $p$ and  the
elementary tropical modification of $\Gamma$ at the right-hand endpoint of $[0,l]$ are the
same metric graph.
\end{eg}

\begin{defn}\label{def:trop morph}
Let  $ \Gamma$ (resp.\ $\Gamma'$) be a representative
of a $\Lambda$-tropical curve $C$ (resp.\ $C'$),  
and assume we are given a harmonic morphism of $\Lambda$-metric graphs
$\phi:  \Gamma' \rightarrow \Gamma$.

An \emph{elementary tropical modification of $\phi$} is 
a harmonic morphism $\phi_1:  \Gamma'_1 \rightarrow \Gamma_1$ of $\Lambda$-metric graphs, 
where $\tau:\Gamma_1\to \Gamma$ is an
elementary tropical modification,  $\tau':\Gamma'_1\to \Gamma'$ is a
tropical modification, and such that $\phi\circ\tau' = \tau\circ\phi_1$.

A \emph{tropical modification of $\phi$} is a finite sequence of
elementary tropical modifications of $\phi$.

Two harmonic morphisms $\phi_1$ and $\phi_2$ of $\Lambda$-metric graphs 
are said to be \emph{tropically equivalent} if there exists a 
harmonic morphism which is a tropical modification of both
$\phi_1$ and $\phi_2$.

A \emph{tropical morphism of tropical curves $\phi:C'\to C$} is a harmonic morphism
of $\Lambda$-metric graphs between some representatives of $C'$ and  $C$,
considered up to (the equivalence relation generated by) tropical
equivalence, and which has a finite representative. 

One makes similar definitions for morphisms of 
augmented tropical curves,
with the additional condition that all harmonic morphisms should be effective.
\end{defn}

Note that it might happen that two non-equivalent  
 morphisms of augmented metric graphs represent the same tropical
morphisms of non-augmented tropical curves.

\begin{rem}
The collection of $\Lambda$-metric graphs (resp.\ augmented $\Lambda$-metric graphs), together with harmonic morphisms (resp.\ effective harmonic morphisms) between them, 
forms a category.  Except for the condition of having a finite representative, one could try to think of tropical curves, together with tropical morphisms between them, as the localization of this category with respect to tropical modifications.
However, there are some technical problems which arise when one tries to make this rigorous (at least if we demand that the localized category admit a calculus of fractions): 
as we will see in Example~\ref{ex:trop morph}, tropical equivalence is
\emph{not}  a transitive relation between morphisms of $\Lambda$-metric graphs.
On the other hand, the restriction of tropical equivalence of
morphisms 
(resp.\ of augmented morphisms) to the collection of finite morphisms
(resp.\ of generically \'etale morphisms) \emph{is} transitive (and hence an equivalence relation).
This is one reason why we include the condition that $\phi$ has a finite representative in our definition of a morphism of tropical curves; another reason is that all morphisms of tropical curves which arise from
algebraic geometry automatically satisfy this condition.
See~\parref{par:alg.trop.curves.functor}.
\end{rem}

\begin{eg}\label{ex:trop morph}
The  morphism of (totally degenerate augmented) 
metric graphs depicted in Figure 
\ref{fig:ex morph augmented}(b) (resp.\ \ref{fig:ex trop morphism}(b))
is an elementary tropical modification of
the one depicted in \ref{fig:ex trop morphism}(a) 
(resp.\ \ref{fig:ex morph augmented}(b)).

The tropical morphisms $\phi_1$ and $\phi_2$
of totally degenerate augmented tropical
curves depicted in Figure \ref{fig:ex trop morphism}(c) and (d) are both
elementary tropical modifications of the morphism $\phi$ 
depicted in Figure 
\ref{fig:ex trop morphism}(e).
\begin{figure}[h]
\begin{tabular}{lll}

\includegraphics[width=2.5cm,
  angle=0]{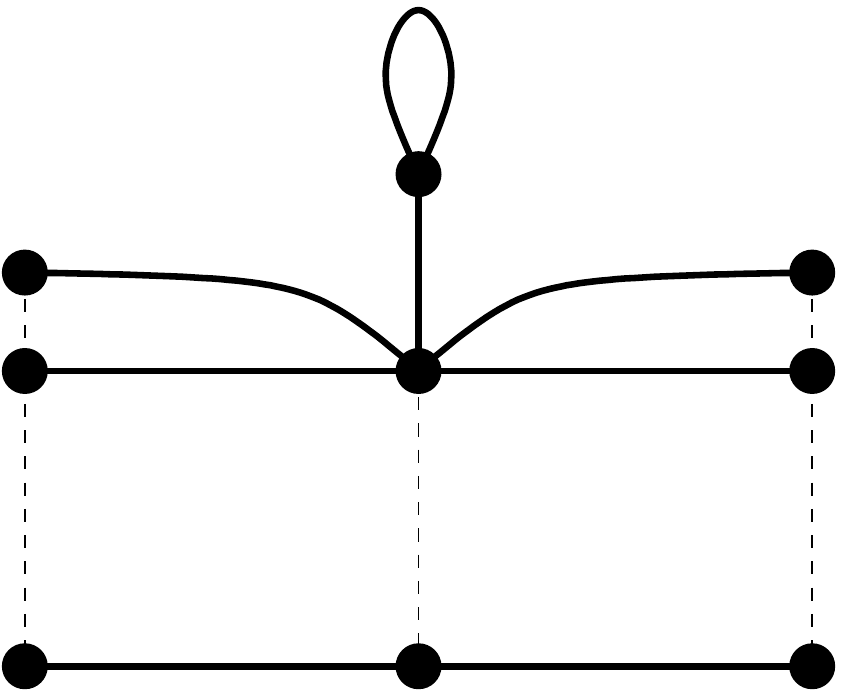} \hspace{3em} &
\includegraphics[width=2.9cm, angle=0]{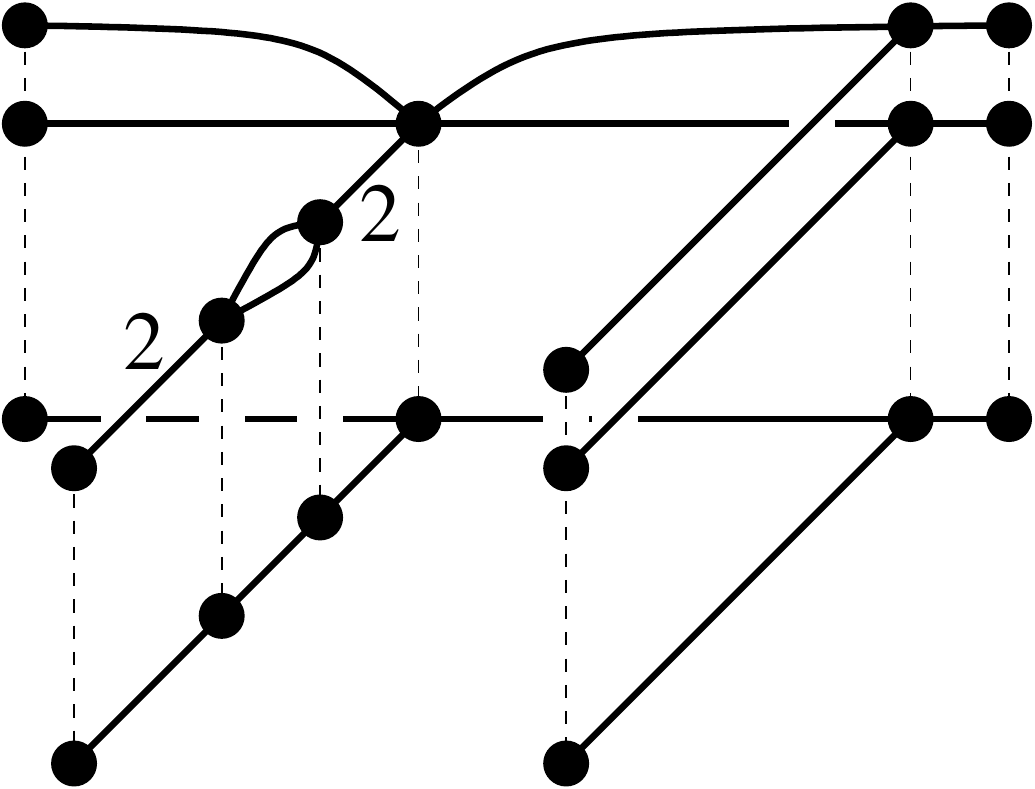}
\hspace{3em} &
\scalebox{.32}{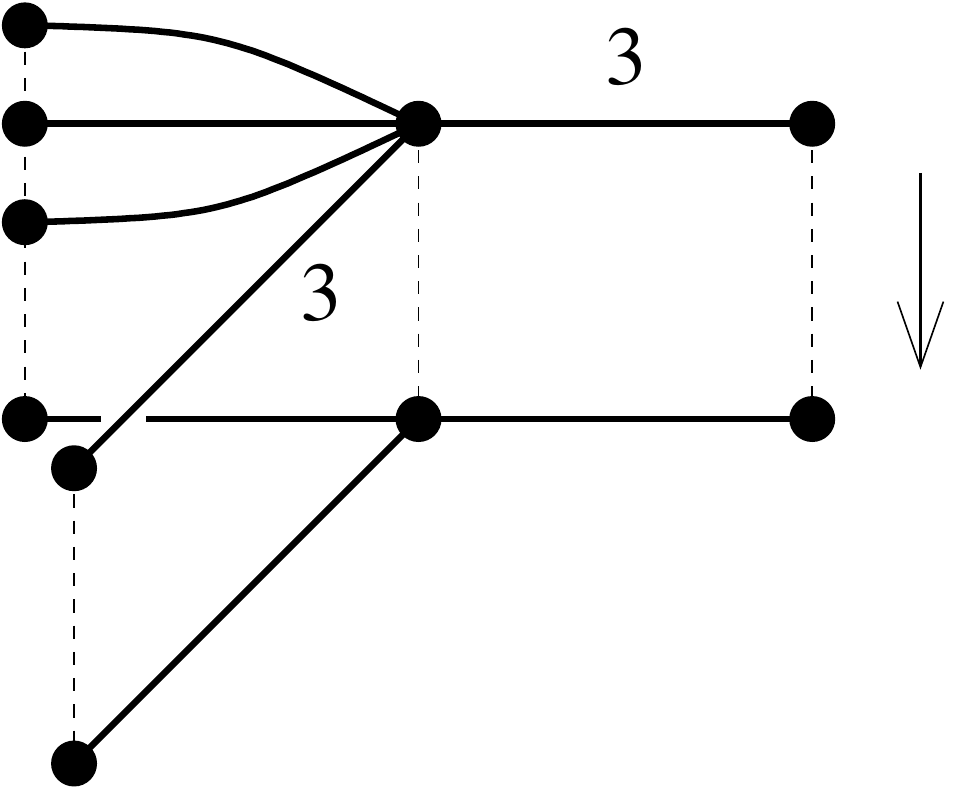}
\\ \\a) & b) &c) $g(p')=0$, $R_{p'}=0$ 
\end{tabular}

\bigskip

\begin{tabular}{ll}
\scalebox{.32}{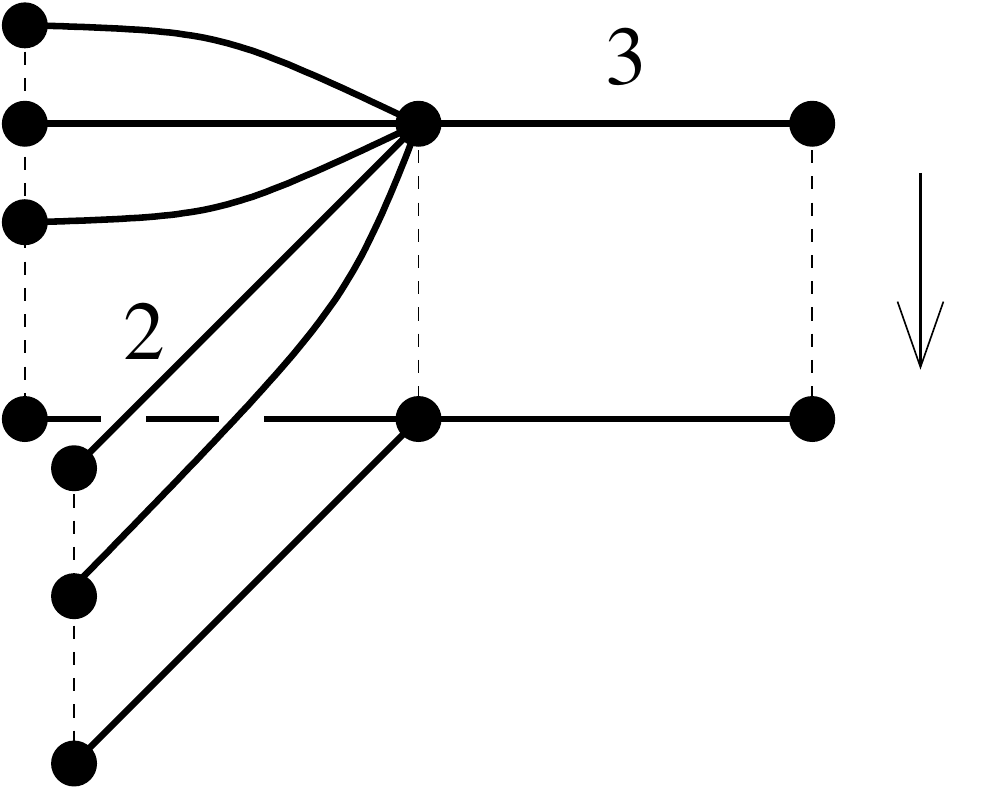} \hspace{3em} & 
\scalebox{.32}{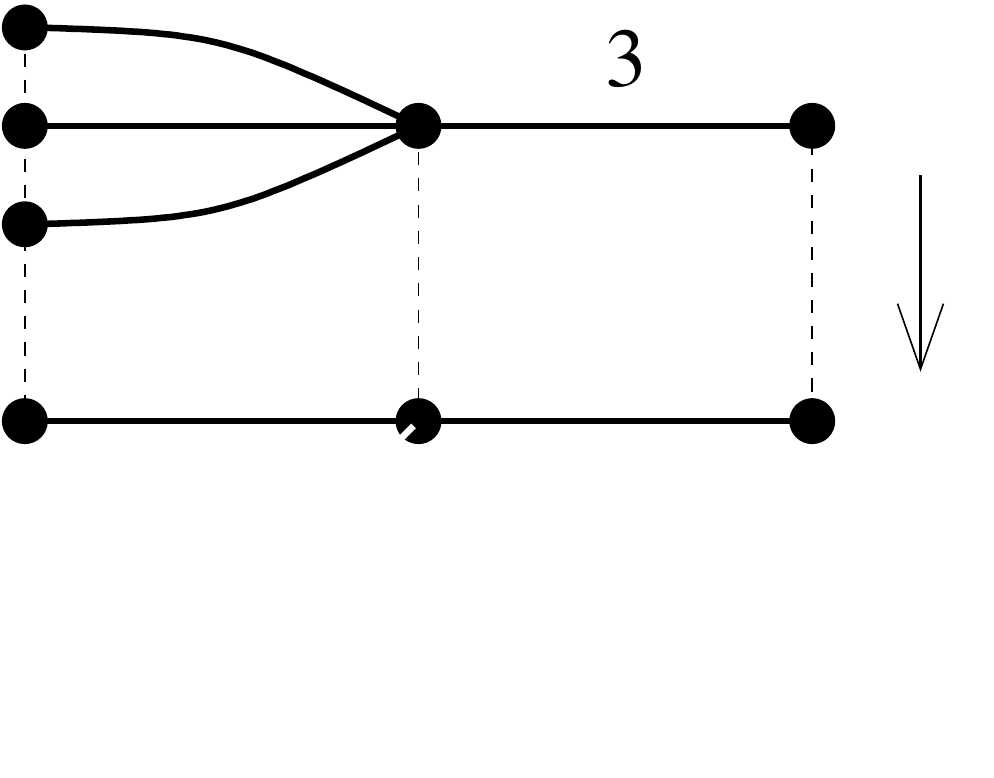}
\\ \\  d) $g(p')=0$, $R_{p'}=1$
& e) $g(p')=0$, $R_{p'}=2$
\end{tabular}
\caption{}\label{fig:ex trop morphism}
\end{figure}

The tropical morphisms $\phi_1$ and $\phi_2$
 depicted in Figure \ref{fig:ex trop morphism2}(a) and (b) are both
elementary tropical modifications of the morphism $\phi$ 
depicted in Figure 
\ref{fig:ex trop morphism2}(c).
\begin{figure}[h]
\begin{tabular}{ccc}

\scalebox{.32}{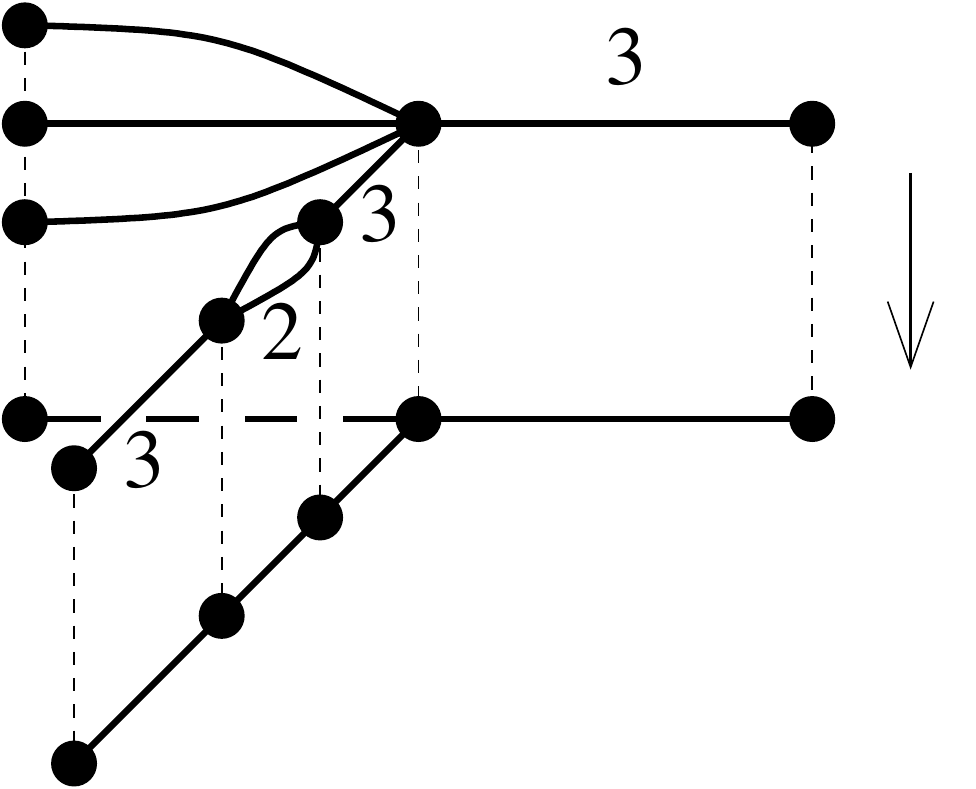}& \hspace{.4cm}
\scalebox{.32}{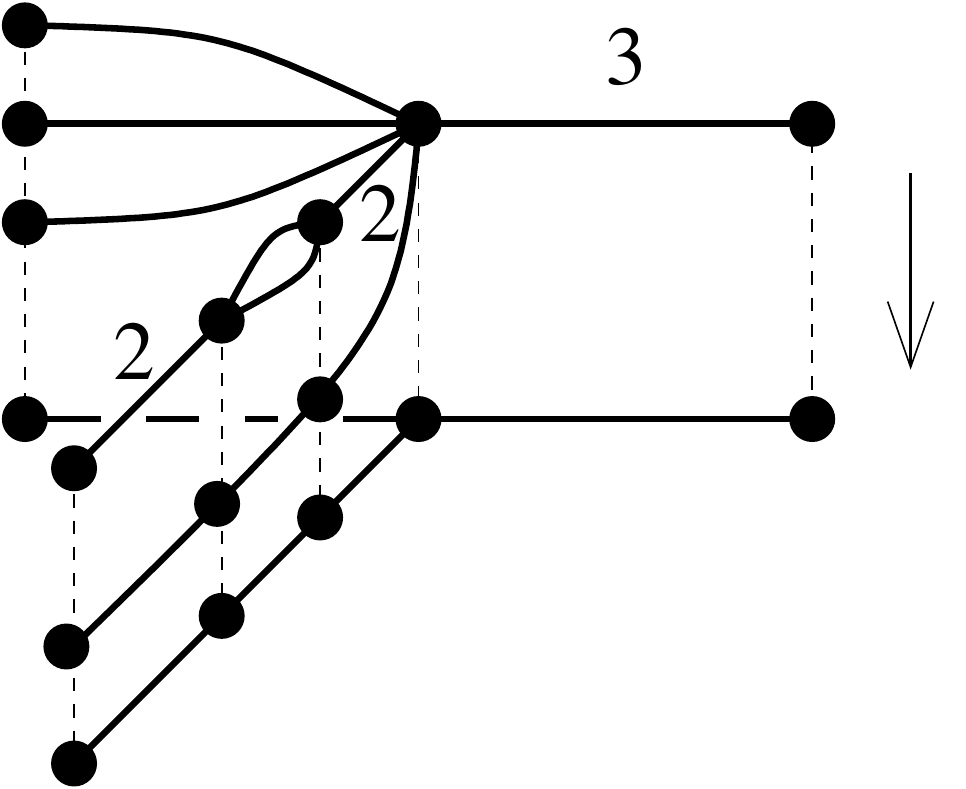}& \hspace{.5cm}
\scalebox{.32}{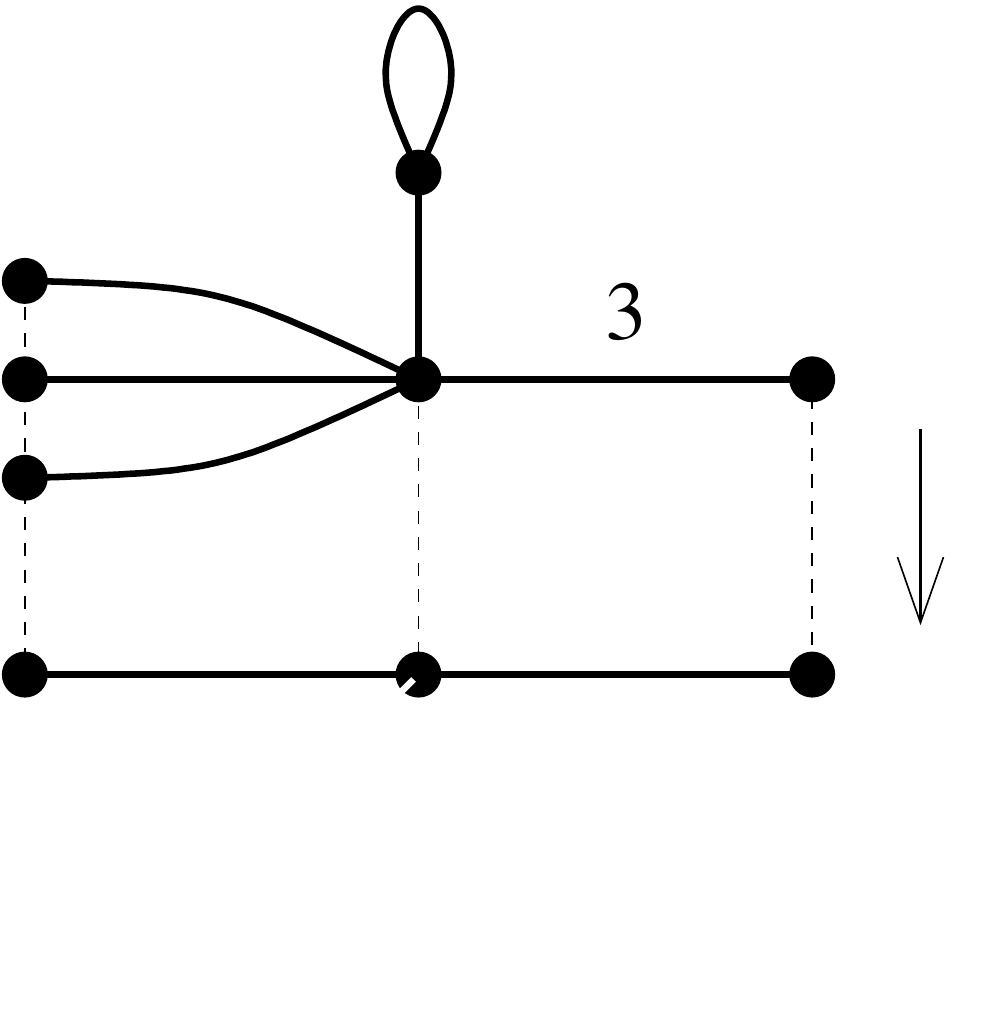}
\\ \\a) & b) &c)

\end{tabular}
\caption{}\label{fig:ex trop morphism2}
\end{figure}

On the other hand, the harmonic morphism $\phi:\Gamma'\to\Gamma$ depicted in
Figure \ref{fig:ex morph augmented}(c) 
with $d=1$ is not tropically equivalent
to any finite morphism: since $\phi$ has degree $1$,
the cycle of the source graph will  be contracted to a point by
any harmonic morphism of metric graphs tropically equivalent to
$\phi$.
In particular, $\phi$ does not give rise to a tropical morphism.
\end{eg}

\smallskip
As mentioned above, tropical equivalence is not transitive among morphisms of
metric graphs (resp.\ of  augmented metric graphs). 
For example, the two morphisms $\phi_1$ and $\phi_2$ depicted 
in Figure \ref{fig:ex trop morphism}(c) and (d)  are not
tropically equivalent as augmented morphisms:
since $R_{p'}=0$ in Figure
\ref{fig:ex trop morphism}(c), any edge appearing in a tropical
modification of $\phi_1$ will have degree $1$.

Note that the preceding harmonic morphisms $\phi_1$ and $\phi_2$ {\em are}
tropically equivalent as morphisms of metric graphs (i.e.  forgetting
the genus function).
However, tropical equivalence is not transitive for tropical morphisms either, for essentially the same reason: 
the two tropical morphisms $\phi_1$ and $\phi_2$
 depicted in Figure \ref{fig:ex trop morphism2}(a) and (b)  are not
tropically equivalent.

Nevertheless, the restriction of tropical equivalence of morphisms
to the set of finite morphisms (resp.\ generically \'etale morphisms)
 is an
equivalence relation. Hence a  tropical morphism (resp.\ an augmented
tropical morphism) can also
be thought as an equivalence class of finite
harmonic morphisms (resp.\ generically \'etale morphisms). 
In particular there exists a natural
composition rule for tropical morhisms (resp.\ augmented tropical
morphisms), turning 
tropical
curves (resp.\ augmented tropical curves)
 equipped with tropical morphisms into a category.

\begin{rem}
In the definition of a tropical morphism of augmented
tropical curves, in
addition to the condition of being a harmonic morphism and the ``up
to tropical modifications'' considerations, we
imposed two rather strong conditions, namely being effective
and having a finite representative.
We already saw in Example \ref{ex:trop morph} that the finiteness condition is
non-trivial. 
The effectiveness condition is also non-trivial: for example,
the harmonic morphism $\phi:\Gamma'\to\Gamma$ of totally degenerate augmented
metric graphs
depicted in Figure \ref{fig:ex morph augmented}(c) 
with $d=2$ is not tropically equivalent
to any finite effective morphism of totally degenerate augmented
metric graphs.  Indeed, for any tropical modification of $\phi$ which is effective, at most two edges 
adjacent to $p'$  can have degree $2$; since $\Gamma'$ already has
two such edges for $\phi$, any tropical modification of $\phi$ which is finite and effective will
contract the cycle of $\Gamma'$  to a point.

 We refer to~\cite{Br12} for a general definition of a tropical morphism
$\phi:C\to X$ from an augmented tropical curve to a non-singular
tropical variety, including Definition \ref{def:trop morph} as a
particular case. 
\end{rem}

\paragraph[Algebraic and tropical curves] 
\label{par:alg.trop.curves.functor}
Restating Lemma~\ref{lem:larger.vertex.set} and
Remark~\ref{rem:skeleta.modifications}, we have:

\begin{prop*}
  Let $(X,V\cup D)$ be a triangulated punctured $K$-curve.  
  Let $D'\subset X(K)$ be a finite set and let $V'$ be a semistable vertex
  set of $(X,D')$, so $(X,V'\cup D')$ is another triangulated punctured
  $K$-curve with underlying curve $X$.  Then the augmented metric graphs
  underlying $\Sigma(X,V'\cup D')$ and $\Sigma(X,V\cup D)$ represent the
  same tropical curve.
\end{prop*}

\smallskip
The above Proposition implies that one can associate a canonical
(augmented) tropical curve to any smooth proper connected $K$-curve $X$. 
This association is \emph{functorial} by
Corollary~\ref{cor:skel.compat.exists}: 

\begin{prop*}
  Let $\phi:X'\to X$ be a finite morphism of smooth proper connected
  $K$-curves, let $D\subset X(K)$ be a finite set, and let
  $D' = \phi\inv(D)$.  Then there exist semistable vertex sets 
  $V,V'$ of $(X,D)$ and $(X',D')$, respectively, such that
  $\phi$ induces a finite morphism of triangulated punctured curves
  $\phi:(X',V'\cup D')\to(X,V\cup D)$.  In particular, $\phi$ induces a
  finite harmonic morphism on suitable choices of skeleta.
\end{prop*}

\smallskip
Again we emphasize that a tropical morphism of tropical curves functorially induced by a
finite morphism of algebraic curves is \emph{effective} and has a
\emph{finite representative}.

\begin{defn}
  We say that a tropical morphism of tropical curves $\bar\phi:C'\to C$ 
  is \emph{liftable} provided that there exists a finite morphism of
  smooth proper connected $K$-curves $\phi:X'\to X$ functorially inducing
  $\bar\phi$ on skeleta in the above sense.
\end{defn}

\smallskip
We will also make use in the text of the notion of tropical modifications of
metrized complexes of curves. 

\begin{defn}
Let
$\cC_0$ be a $\Lambda$-metrized complex of $k$-curves. 
\begin{itemize}
\item A \emph{refinement}  of $\cC_0$ is any
$\Lambda$-metrized complex of $k$-curves
 $\cC$ obtained from $\cC_0$ by adding a
finite set of $\Lambda$-points $S$ of $\cC_0 \setminus V(\cC_0) $ to
the set $V(\cC_0)$ of vertices of $\cC_0$ 
(see Definition \ref{def:metr cplx}),
setting $C_p=\PP^1_k$ for all $p\in S$, and defining the map $\red_p$
by choosing any two distinct closed points of $C_p$. 
\item An \emph{elementary tropical modification} of $\cC_0$ is a
$\Lambda$-metrized complex of $k$-curves $\cC$ obtained from $\cC_0$ 
by an elementary tropical modification of the underlying metric
graph at a vertex $p$ of $\cC$, 
with the map $\red_p$ extended to $e$
by choosing any closed point of $C_p$ not in the image of the reduction map
for $\cC_0$. 
\item Any metrized complex of curves
 $\cC$ obtained from a metrized complex of curves
 $\cC_0$ by a finite sequence of 
 refinements and elementary tropical modifications is called a {\em tropical 
modification} of $\cC_0$.
\end{itemize}
\end{defn}

\section{Lifting harmonic morphisms of metric graphs to morphisms of
  metrized complexes}\label{sec:lifting morphisms metric graphs}

There is an obvious forgetful functor which sends 
metrized complexes of curves to
(augmented) metric graphs, and harmonic morphisms of metrized complexes to 
harmonic morphisms of (augmented) metric graphs. 
A harmonic morphism of (augmented) metric graphs is said to be
\emph{liftable to  a harmonic morphism of metrized  complexes of $k$-curves} 
if it lies in the image of the forgetful functor.
 
We proved in Theorem \ref{thm:lifting} that every tame covering
of metrized complexes of curves can be lifted to a tame covering of
algebraic curves. 
In this section we study the problem of lifting harmonic morphisms
of (augmented) metric graphs to finite morphisms of metrized complexes
(and thus to tame coverings of proper smooth curves, thanks to 
Proposition~\ref{prop:lifting}).

\paragraph[Lifting finite augmented morphisms] 
Recall that $k$ is an algebraically closed field  of
characteristic $p\geq 0$. A finite harmonic morphism $\phi$ of (augmented)
metric graphs is called a {\it tame harmonic morphism} 
if either $p=0$ or all the local degrees of $\phi$ along edges  are prime to $p$.
Lifting of tame harmonic morphisms of augmented metric graphs to 
tame harmonic morphisms of metrized complexes of $k$-curves is equivalent to the
existence of tamely ramified covers of 
$k$-curves of given genus with some given prescribed ramification profile.

 \medskip

\subparagraph \label{par:def A}
A partition $\mu$ of an integer $d$ is a multiset of natural numbers $d_1,\dots, d_l \geq 1$ with $\sum_i d_i =d$. 
 The integer $l$, called the length of $\mu$, will be denoted by $l(\mu)$.

Let $g', g\ge 0$ and $d>0$ be integers, and let  
 $M=\{\mu_1,\ldots,\mu_s\}$ be a collection of
$s$ partitions of $d$. 
Assume that the integer $R$ defined by
\begin{equation}\label{eq:R}
R:=d(2-2g) + 2g' -2 -sd + \sum_{i=1}^sl(\mu_i)
\end{equation}
is non-negative. Denote by 
$\mathcal A_{g',g}^d(\mu_1,\dots,\mu_s)$ 
the set of all tame
coverings $\phi: C ' \to C$ of smooth proper curves over $k$, with the following properties:

\begin{itemize}
\item[(i)] The curves $C$ and $C'$ are irreducible of genus $g$ and $g'$, respectively;

\item[(ii)] The degree of $\phi$ is equal to $d$;  

\item[(iii)] The branch locus of $\phi$ contains (at least) $s$ distinct points
 $x_1,\dots, x_s\in C$, and the ramification profile of $\phi$ at the points 
 $\phi^{-1}(x_i)$ is given by $\mu_i$, for 
$1\le i\le s$.
\end{itemize}

As we will explain now, the lifting problem for morphisms of augmented metric graphs to  morphisms of
metrized complexes over a field $k$ reduces to the emptiness or non-emptiness of certain sets 
$\mathcal A_{g' ,g}^d(\mu_1,\dots,\mu_s)$. This latter problem is quite
subtle, and no complete satisfactory answer is yet known (see also~\parref{par:hurwitz}). In some simple cases, however, one can
ensure that $\mathcal A_{g',g}^d(\mu_1,\dots,\mu_s)$ is non-empty. 
For example, if all the partitions $\mu_i$ are trivial (i.e., they each consist of $d$ 1's), then $\mathcal A_{g',g}^d(\mu_1,\ldots,\mu_s)$ 
is non-empty.  
Here is another simple example. 

\begin{eg}\label{ex:non-vanishing h}
For an integer $d$ prime to characteristic $p$ of $k$, the set
$\mathcal A_{0,0}^d((d),(d))$ is non-empty since it contains the map 
$z\mapsto z^d$. This is in fact the only map in $\mathcal
A_{0,0}^d((d),(d))$ up to the action of the group $\mathrm{PGL}(2,k)$ on the target curve and
${\mathbf P}^1$-isomorphisms of coverings. 
\end{eg}

\subparagraph
 Let $\phi: \Gamma' \rightarrow \Gamma$ be a finite harmonic morphism
 of augmented metric graphs. 
 Using the definition of a harmonic morphism, one can associate to any point $p'$ of
$\Gamma'$ a collection 
$\mu_1(p'),\dots,\mu_s(p')$ of $s$
partitions of the integer $d_{p'}(\phi)$, where $s=\val(\phi(p'))$, as follows: if
$T_{\phi(p)}(\Gamma) = \{ v_1,\dots,v_s \}$
denotes all the tangent directions to $\Gamma$ at $\phi(p')$,
then $\mu_i(p')$ is the partition of
$d_{p'}(\phi)$ which consists of the various local degrees of $\phi$ in
all tangent directions $v'\in T_{p'}(\Gamma')$ mapping to $v_i$.

\medskip
 
The next proposition is an immediate consequence of the various definitions involved
once we note that, by Example \ref{ex:non-vanishing h}, 
there are only finitely points $p' \in \Gamma'$ for which the question of
non-emptiness of the sets $\mathcal A_{g(p'), g(\phi(p'))}^{d_{p'}(\phi)}$
arises.  It provides a ``numerical criterion'' for a tame harmonic morphism of
augmented metric graphs to be liftable to a tame harmonic morphism of metrized
complexes of curves.

\begin{prop}\label{prop:lift augmented to complex} 
Let $\phi:\Gamma' \to\Gamma$ be a tame harmonic morphism of augmented
metric graphs.  
Then $\phi$ can be lifted to a tame harmonic morphism of metrized complexes over $k$
if and only if 
for every point $p'$ in $\Gamma'$, the set 
$\mathcal A_{g(p'),g(\phi(p'))}^{d_{p'}(\phi)}(\mu_1(p'),\ldots,\mu_{\val( \phi( p' ))}(p'))$ is non-empty. 
\end{prop}

\subparagraph\label{par:hurwitz}
In characteristic $0$, the lifting problem for finite augmented
morphisms of metric graphs can be further reduced to a vanishing question for certain Hurwitz
numbers.

Fix an irreducible smooth proper
curve $C$ of genus $g$ over $k$, and let 
$x_1,\dots, x_s,y_1,\ldots,y_R$  be a set of distinct points on $C$. 
The Hurwitz set  $\mathcal H_{g',g}^d(\mu_1,\ldots,\mu_s)$ is the
set of $C$-isomorphism classes of all coverings in $\mathcal A_{g',g}^d(\mu_1,\ldots,\mu_s)$ satisfying 
(i), (ii), and (iii)
 in~\parref{par:def A} for the curve $C$ and the points $x_1,\dots, x_s$ that we have fixed, and
which in addition satisfy: 
\begin{itemize}
\item[(iv)] The integer $R$ is given by~\eqref{eq:R}, and for each $1\le i\le R $, $\phi$ 
has a unique simple ramification point $y'_i$ lying above $y_i$.
\end{itemize}
Note that by the above condition, the branch locus of $\phi$ consists precisely of the points $x_i, y_j$.
The  Hurwitz number $H_{g',g}^d(\mu_1,\ldots,\mu_s)$ is defined as
$$H_{g',g}^d(\mu_1,\ldots,\mu_s):=\sum_{\phi\in \mathcal H_{g',g}^d(\mu_1,\ldots,\mu_s)} \frac{1}{|\mathrm{Aut_{C}}(\phi)|},$$
 and  does not depend on the choice of $C$ and the closed points
$x_1,\dots, x_s,y_1,\ldots,y_R$ in $C$.

\begin{eg}\label{rem:Hurwitznumbers}
It is known, see for example \cite{AKS84}, that
$$H_{g,0}^2=\frac{1}{2}, \quad  H_{g,0}^3((3),\ldots(3))>0, \quad
H_{0,0}^4((2,2),(2,2),(3,1))=0.$$
For the reader's convenience, and since we will use it several times in the sequel, we sketch a proof
of the fact that $H_{0,0}^4((2,2),(2,2),(3,1))=0$.
By the Riemann--Hurwitz formula and the Riemann Existence Theorem,
$H_{0,0}^4((2,2),(2,2),(3,1))\neq 0$ if and only if there exist elements $\sigma_1,\sigma_2,\sigma_3$ in the symmetric group ${\mathfrak S}_4$ having cycle decompositions of type $(2,2),(2,2),(3,1)$, respectively, such that $\sigma_1\sigma_2\sigma_3=1$ and such that the $\sigma_i$ generate a transitive subgroup of ${\mathfrak S}_4$.
However, elementary group theory shows that the product $\sigma_1\sigma_2$ cannot be of type $(3,1)$
(the transitivity condition does not intervene here).
For a proof which works in any characteristic, see Lemma~\ref{lem:non deg 4} below.
\end{eg}

\medskip

All  Hurwitz numbers can be theoretically  computed, for example using
 Frobenius Formula (see \cite[Theorem A.1.9]{LZZ}). Nevertheless, 
the problem of understanding their vanishing 
is wide open. The above example shows that Hurwitz numbers in degree
at most three are all positive, which is not the case 
in degree four. Some families of (non-)vanishing Hurwitz
numbers are known (see Example \ref{ex:h non-van}). However,
in general one has to explicitly compute a given Hurwitz number
to decide if this latter vanishes or not. We refer the reader to
\cite{AKS84}, \cite{PP06}, and \cite{PP08}, along with the references therein, 
for an account of what is known about this subject.
We will use the vanishing of
$H_{0,0}^4((2,2),(2,2),(3,1))$ in Section~\ref{sec:examples} to
construct a $4$-gonal augmented graph (see Section~\ref{sec:examples} for the definition) which cannot be lifted 
to any $4$-gonal proper smooth algebraic curve over $K$.
\begin{eg}\label{ex:h non-van}
Some partial results are known concerning the (non-)vanishing of Hurwitz numbers. For
example, it is known that double Hurwitz numbers (i.e., when $s=2$) are
all positive (this can be seen for example from the presentation of the cut-join equation given 
in~\cite{CJM10}), 
as well as all the Hurwitz numbers $H_{g',g}^d(\mu_1,\dots,\mu_s)$ when
$g\ge 1$ and $R\geq 0$ (\cite{Hus62, AKS84}). On the other hand, it is proved in \cite{PP08}
that
$$H_{0,0}^{d}\left((d-2,2),(2,\ldots,2),(\frac{d}{2}+1,1,\ldots,1)
\right)=0 \quad \mbox{for all }d \geq 4\mbox{ even.} $$
\end{eg}

\begin{eg} \label{eg:DM}
As another example of non-vanishing Hurwitz numbers, one has $H_{0,0}^{d'}(\mu_1,\ldots,\mu_s,(d')) >0$ 
for all integers $d' \geq 1$ when the integer $R$ defined in \eqref{eq:R} is zero (i.e., if the combinatorial Riemann--Hurwitz formula holds);
see \cite[Proposition 5.2]{AKS84} or \cite[Proposition 7.2]{DeMarco-McMullen}.
\end{eg}

\medskip

The non-emptiness of $\mathcal A_{g',g}^d(\mu_1,\ldots,\mu_s)$ can 
 be reduced to the non-emptiness of the Hurwitz set $\mathcal H_{g',g}^d(\mu_1,\ldots,\mu_s)$. 
\begin{lem}\label{lem:AH} Suppose that $k$ has characteristic $0$. Then  $\mathcal A_{g',g}^d(\mu_1,\ldots,\mu_s)$
  is non-empty if and only if $H_{g',g}^d(\mu_1,\ldots,\mu_s)\ne 0$.
\end{lem}

\pf
Since $ \mathcal H_{g',g}^d(\mu_1,\ldots,\mu_s)$ is a subset of  $\mathcal A_{g',g}^d(\mu_1,\ldots,\mu_s)$, 
obviously we only need to prove that if $\mathcal A_{g',g}^d(\mu_1,\ldots,\mu_s)\neq \emptyset$, 
then the Hurwitz set is also non-empty. Let $\phi: C' \rightarrow C$ be an element of 
$\mathcal A_{g',g}^d(\mu_1,\ldots,\mu_s)$, branched over $x_i \in C$ with ramification profile 
$\mu_i$ for $i=1,\dots, s$, 
and let $z_1, \dots, z_t$ be all the other points in the branch locus of $\phi$. 
Denote by $\nu_i$ the ramification profile of $\phi$ above the point $z_i$. 
Fix a closed point $\star$ of $C \setminus\{x_1,\dots,x_s, z_1,\dots, z_t\}$.  
The \'etale fundamental group $\pi_1(C \setminus\{x_1,\dots,x_s, z_1,\dots, z_t\}, \star)$ 
is the profinite completion of 
the group generated by a system of generators $a_1,b_1,\dots,a_{g}, b_{g}, c_1,\dots, c_{s+t}$ satisying 
 the relation $[a_1,b_1]\dots[a_{g},b_{g}] c_1\dots c_{s+t}=1$, where $[a,b] = aba^{-1}b^{-1}$ 
(see~\cite{SGA1}). 
In addition, the data of $\phi$ is equivalent to the data of a surjective morphism $\rho$ from 
$\pi_1(C \setminus\{x_1,\dots,x_s, z_1,\dots, z_t\}, \star)$ to a transitive subgroup of the symmetric group 
$\mathfrak{S}_{d}$ of degree $d$ such that the partition $\mu_i$ (resp.\ $\nu_i$) of $d$ corresponds to the 
lengths of the cyclic permutations in the decomposition of $\rho(c_i)$ (resp.\ $\rho(c_{s+i})$) 
in $\mathfrak{S}_{d}$ into products of cycles, for $1\leq i\leq s$ (resp.\ $1\leq i\leq t$).
By Riemann--Hurwitz formula, we have $R = \sum_{i=1}^t (d -l(\nu_i))$.  

Now note that each $\rho(c_{s+i})$ can be written as a product of $d-l(\nu_i)$ 
transpositions $\tau_i^1,\dots, \tau^{d-l(\nu_i)}_i$ in $\frak{S}_d$, i.e., 
$\rho(c_{s+i}) = \tau_i^1\dots \tau^{d-l(\nu_i)}_i$. 
Rename the set of $R$ distinct points $y_1,\dots, y_R$ of $C \setminus \{x_1,\dots, x_s, \star\}$ 
as $z_i^1,\dots, z_{i}^{d-l(\nu_i)}$ for $1\leq i\leq t$. 

 The \'etale fundamental group $\pi_1(C \setminus\{x_1,\dots,x_s, 
z_1^1,\dots, z^{d-l(\nu_1)}_1,\dots, z_t^{d-l(\nu_t)}\}, \star)$ has, as a profinite group, 
a system of generators  $a_1,b_1, \dots, a_g, b_g, c_1,\dots, c_s, c_{s+1}^{1},\dots, c_{s+1}^{d-l(\nu_1)},\dots, c_{s+t}^{d-l(\nu_t)}$ 
verifying the relation  
$$[a_1,b_1]\dots[a_{g},b_{g}] c_1\dots c_s c_{s+1}^1\dots c_{s+1}^{d-l(\nu_1)}\dots c_{s+t}^{1}\dots 
c_{s+t}^{d-l(\nu_t)}=1,$$ and admits a surjective morphism to $\frak{S}_d$  
which coincides with $\rho$ on $a_1,b_1,\dots, a_g,b_g$, and which sends 
$c_{s+i}^{j}$ to $\tau_i^j$ for each $1\leq i\leq t $ and $1\leq j\leq d-l(\nu_i)$. 
The corresponding cover $C''\rightarrow C$ obviously belongs to $\mathcal  A_{g',g}^d(\mu_1,\ldots,\mu_s)$ and 
in addition has simple ramification profile $(2)$ above each $y_i$, i.e., 
it verifies condition (iv) above. This shows that $\mathcal  H_{g',g}^d(\mu_1,\ldots,\mu_s)$ is non-empty. 
\qed

\begin{cor}\label{cor:lift iff H not 0}
Suppose that $k$ has characteristic $0$.
Let $\phi:\Gamma' \to\Gamma$ be a 
 finite morphism  of augmented 
metric graphs, and let $\cC$ be a metrized complex over $k$ lifting  $\Gamma$.
There exists a lifting of $\phi$ to a finite harmonic morphism of metrized complexes $\cC' \to \cC$ over $k$ (and thus to a morphism of smooth proper curves over $K$)
if and only if 
$$\prod_{p'\in V(\Gamma')} H_{g(p'),
  g(\phi(p'))}^{d_{p'}(\phi)}\left(\mu_1(p'),\ldots,\mu_{\val( \phi( p'
  )) }\right)\ne 0.$$ 
In particular, if $\phi$ is effective and $g(p) \geq 1$ for all the points of valency at least three in $\Gamma$, then
$\phi$ lifts to a finite harmonic morphism of metrized complexes over $k$.
\end{cor}

\begin{rem} \label{rem:large characteristic}
If $k$ has positive characteristic $p > d$, then $\mathcal A_{g',g}^d(\mu_1,\ldots,\mu_s)$ has the same cardinality as in 
characteristic zero.  (This follows from \cite{SGA1}, which provides an isomorphism between the tame fundamental group in positive
characteristic $p$ and the prime-to-$p$ part of the {\'e}tale fundamental group in characteristic zero.)
In particular, Lemma~\ref{lem:AH} also holds under the assumption that $p>d$.
\end{rem}

\paragraph[Lifting finite harmonic morphisms]
\label{par:liftingMGMC}
Now we turn to the lifting problem for finite morphisms  of 
non-augmented metric graphs  to  morphisms of metrized complexes of
$k$-curves. In this case 
there are no obstructions to the existence of such a lift.
\begin{thm}\label{thm:lifting harm}
Let $\phi:\Gamma'\rightarrow \Gamma$ be a tame harmonic morphism of 
metric graphs, and suppose that $\Gamma$ is augmented. 
There exists an enrichment of $\Gamma'$ to an augmented metric graph
$(\Gamma',g')$ such that $\phi: (\Gamma',g') \rightarrow (\Gamma,g)$
lifts to a tame harmonic morphism of metrized complexes of curves over $k$ 
(and thus to a morphism of smooth proper curves over $K$).
\end{thm}

\smallskip

Theorem~\ref{thm:lifting harm} is an immediate consequence of
 Proposition~\ref{prop:lift augmented to complex}  and the following theorem. 
(For the statement, we say that a partition $\mu$ of $d$ is {\em tame} if 
either $\mathrm{char}(k)=0$ or all the integers appearing in $\mu$ are prime to $p$.)

\begin{thm}\label{Hurwitz:nonvanishing}
 Let $g\ge 0,d\ge 2, s\ge 1$ be integers. Let $\mu_1,\dots, \mu_s$ be a collection of $s$ 
tame partitions of $d$. 
Then there exists a sufficiently large non-negative integer $g'$ such that $\mathcal A_{g',g}^d(\mu_1,\ldots,\mu_s)$ 
is non-empty.
\end{thm}

\medskip

\pf
We first give a simple proof which works in characteristic zero, and more generally, 
in the case of a tame monodromy group. 
The proof in characteristic $p>0$ is based on our lifting theorem and a deformation argument. 

Suppose first that the characteristic of $k$ is zero. By
Lemma~\ref{lem:AH}, we need to show that  for large enough $g'$ the set
$\mathcal H_{g',g}^d(\mu_1,\ldots,\mu_s)$ is non-empty. 

If $g\geq 1$, for any large enough $g'$
giving $R \geq 0$, we have $\mathcal H_{g',g}^d(\mu_1,\ldots,\mu_s) \neq \emptyset$~\cite{Hus62}. So suppose $g=0$.
 Consider $s+R+1$
distinct points $x_1,\ldots,x_s,z_1,\dots, z_R,\star$ in $C$.  
The \'etale fundamental group $\pi_1( R ):=
\pi_1(C \setminus\{x_1,\dots,x_s, z_1,\dots, z_R\}, \star)$ has, as a profinite group, a
system of generators 
$c_1,\dots, c_s, c_{s+1},\dots, c_{s+R}$ verifying the relation   
$$c_1\dots c_r c_{s+1} \dots c_{s+R}=1.$$
It will be enough to show that for a large enough $R$, there exists a
surjective morphism $\rho$ from $\pi_1( R )$ to $\frak S_d$ so that
$\rho(c_{s+i})$ is a transposition for any $i=1,\dots, R$, and that
for any $i=1,\dots, s$, the partition of $d$ given by the lengths of
the cyclic permutations  in the decomposition of $\rho(c_i)$ is equal to $\mu_i$. 
In this case, the genus $g'$ of the corresponding cover $C'$ of $C$ in
$\mathcal H_{g',0}^d(\mu_1,\dots, \mu_s)$ will be given by 
$$g' = 1 -d + \frac 12\bigl[sd + R - \sum_{i=1}^sl(\mu_i)\bigr].$$

Consider an arbitrary map $\rho$ from 
$\{c_1,\dots,c_s\}$ to $\frak S_d$ verifying the
ramification profile condition for $\rho(c_1),\dots,
\rho(c_r)$. Choose a system of $d$ transpositions $\tau_1,\dots,
\tau_d$ generating $\frak S_d$, and consider a set of transpositions
$\tau_{d+1},\dots, \tau_R$ such that  
$$\rho(c_1)\dots \rho(c_s)\, \tau_1 \dots \tau_d= \tau_{R}\dots
\tau_{d+1}\,.$$ 
This proves Theorem \ref{Hurwitz:nonvanishing} when $k$ has characteristic $0$.

\medskip
Consider now the case of a base field $k$ of positive characteristic
$p>0$. Note that since the prime to $p$ part of the tame fundamental group has
the same representation as in the case of characteristic zero,  the group theoretic
method we used in the previous case can be applied if
the monodromy group is tame, i.e., has size prime to $p$. 
However, in general it is impossible to impose such a
condition on the monodromy group.  For example in the case when 
$p$ divides $d$, the size of the monodromy group is always divisible by $p$. 

\begin{figure}[h]
\begin{tabular}{ccccc}
\scalebox{0.32}{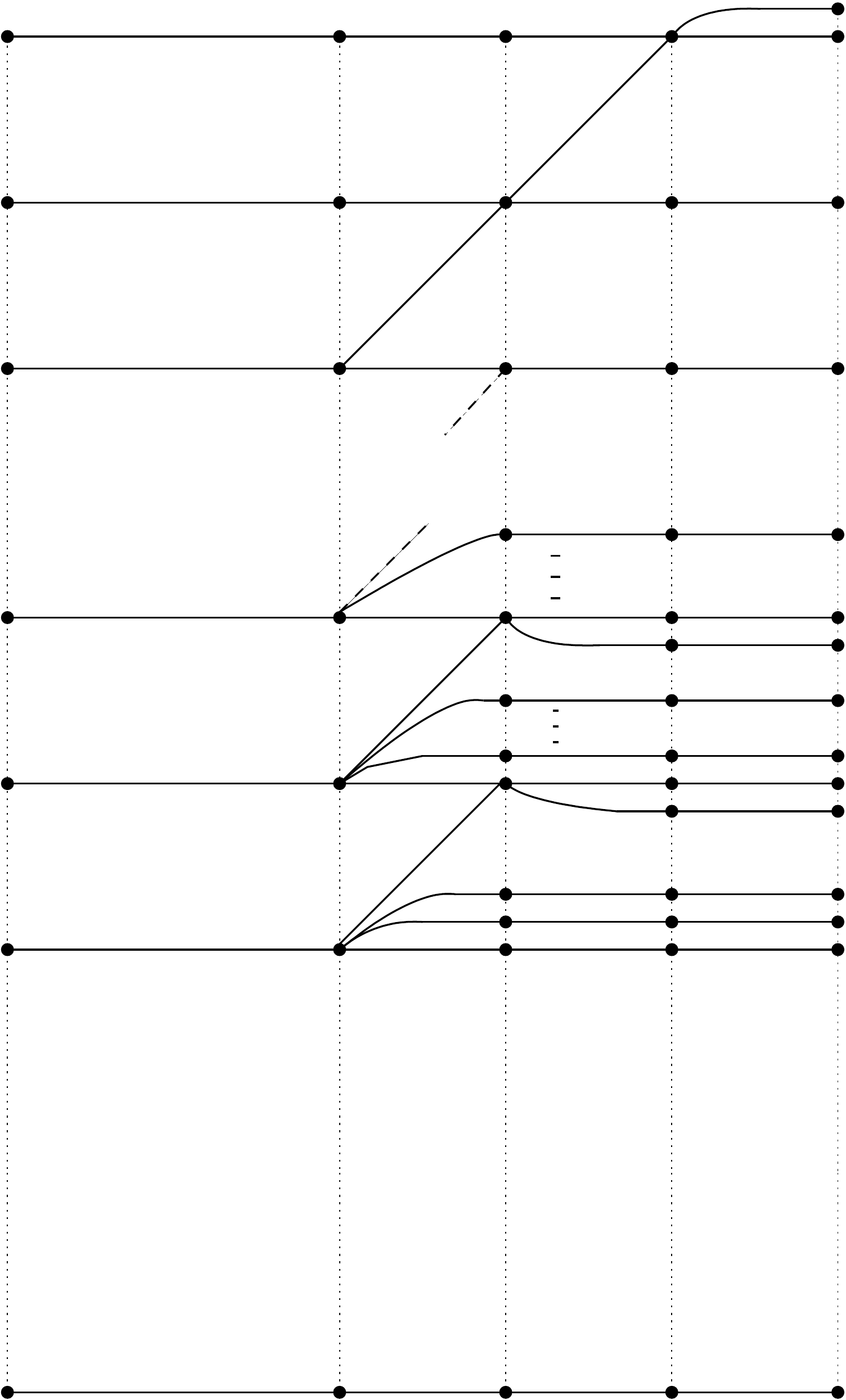} && \hspace{3cm} &&
\scalebox{0.32}{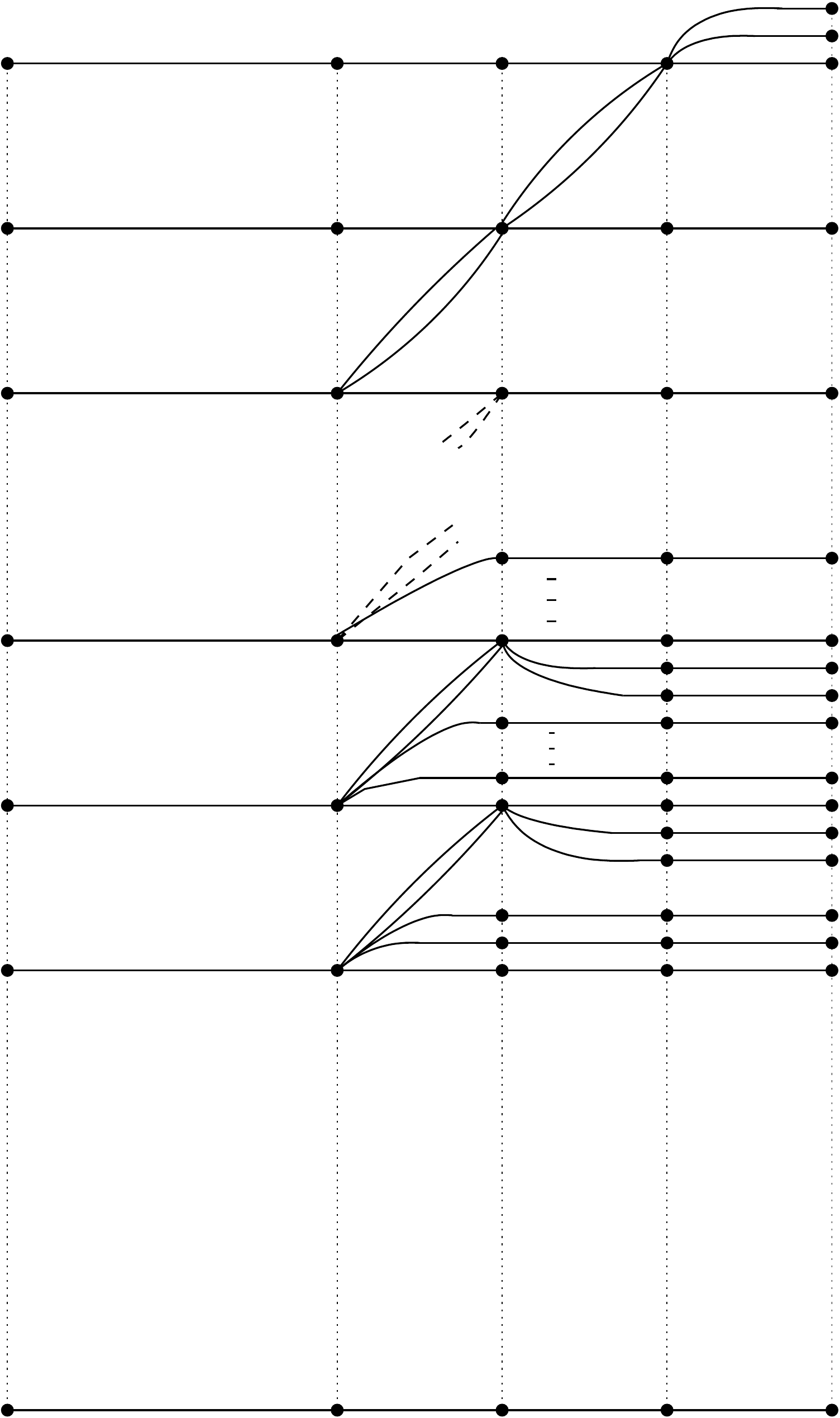}\\
\\ \\a) Reduction in the case $p\ne 2$. && && b) Reduction in the case $p=2$. 
 \end{tabular}
\caption{ }\label{fig:ReductionToCyclic}
\end{figure}

We first describe how to reduce the proof of Theorem~\ref{Hurwitz:nonvanishing} to the case $s=1$ and $g=0$.
Suppose that for each $\mu_i$, $1\leq i\leq s$, there exists a large
enough $g_i$ such that $\mathcal A_{g_i,0}^d(\mu_i)$ is non-empty, and
consider a tame cover $\phi_i: C_i \rightarrow \P^1_k$ in
$\mathcal A_{g_i,0}(\mu_i)$ such that the ramification profile over $0\in \P^1$ is given by $\mu_i$, 
and choose two regular points $x_i,y_i \in \P^1$ (i.e. $x_i,y_i$ are
outside the branch locus of $\phi_i$).  
Choose also a smooth proper curve $C_0$ of genus $g$ which admits a tame cover 
$\phi_0 : C'_0 \to C_0$ of degree $d$ from a smooth proper curve
$C'_0$ of large enough genus $g'_0$.  (The existence of such a cover can be deduced by a similar trick as that discussed at the end of the proof below and depicted in 
Figure~\ref{fig:ReductionToCyclic}.) 
Let $y_0 \in C_0$ be a regular point of $\phi_0$. 

\begin{figure}[h]
\scalebox{0.38}{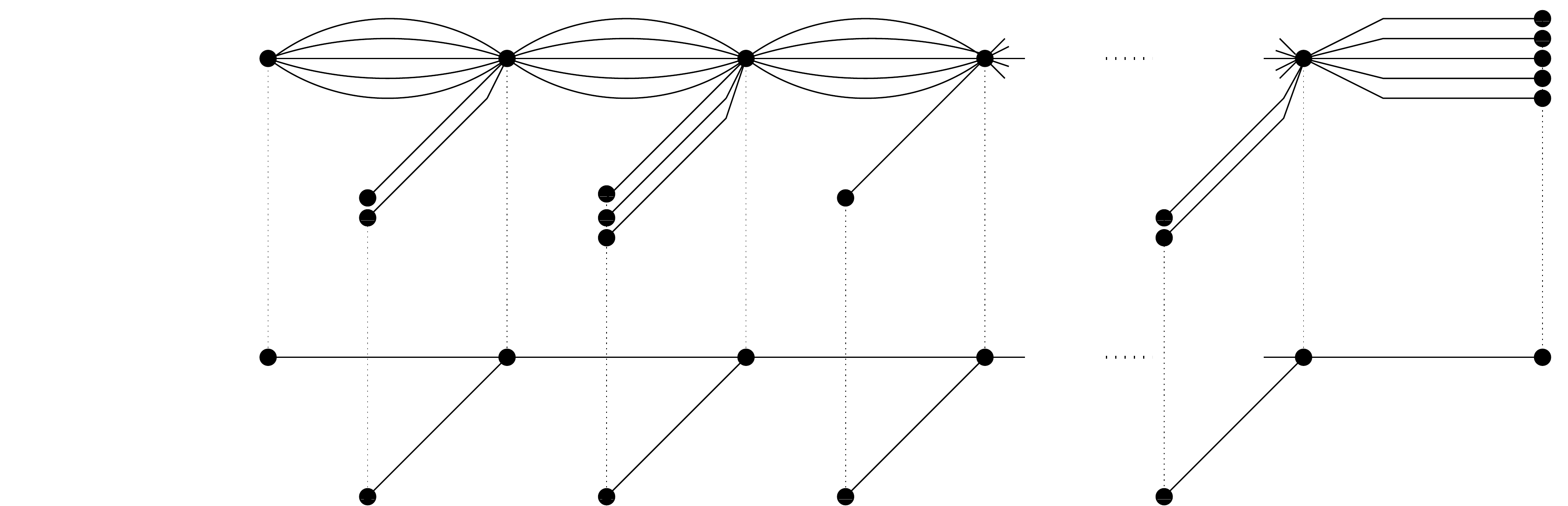}
\caption{ }\label{fig:banana}
\end{figure}

Let $\cC_0$ be the metrized complex over $k$ whose underlying 
 metric graph  is $[0, +\infty]$, with one finite vertex $v_0$ and one infinite vertex $v_\infty$,
 equipped with the metric
induced by $\RR$, and with $C_{v_0}=C_0$ and $\red_{v_0}(\{v_0,v_\infty\}) = y_0$. 
Denote by   $\cC$  the modification of
$\cC_0$ obtained by taking a refinement at $r$ distinct points $0<v_1<\dots< v_s<\infty$, as depicted in Figure~\ref{fig:banana}, 
and by setting $C_{v_i}= \P^1$ and $\red_{v_i}(\{v_i,v_{i-1}\}) = x_i$ and $\red_{v_i}(\{v_i,v_{i+1}\}) =y_{i}$ 
(here $v_{s+1} = v_{\infty}$), and by adding an infinite edge $e_i$ to each $v_i$, and defining $\red_{v_i}(e_i) = 0\in \P^1$.  
Denote by $\Gamma$ the underlying metric graph of $\cC$.  See Figure~\ref{fig:banana}.

Define now the  metric graph $B_{s,d}$ as the chain of $s$ banana
graphs of size $d$:  $B_{s,d}$ has $s+1$ finite vertices 
$u_0,\dots, u_s$ and $u'_1,\ldots,u'_d$ infinite vertices adjacent to $u_s$
 such that $u_i$ is connected to $u_{i+1}$ with $d$
 edges of length $\ell_{\Gamma}(\{v_{i+1}-v_i\})$.
 We denote by $\widetilde B_{s,d}$ the tropical modification of
 $B_{r,d}$ at $u_1,\dots, u_s$, obtained by adding $l(\mu_i)$ infinite edges to $u_i$. 
Eventually we turn $\widetilde  B_{s,d}$ into a metrized complex $\cC_{s,d}$ over $k$ by setting $C_{u_i}=C_i$, and defining 
$\red_{u_i}$ on the $d$ edges between $u_i$ and $u_{i+1}$ by a bijection to the $d$ points in $\phi_i^{-1}(y_i)$, 
$\red_{u_i}$ on the edges between $u_i$ and $u_{i-1}$ by a bijection to the $d$ points in $\phi_i^{-1}(x_i)$, 
and $\red_{u_i}$ on the 
$l(\mu_i)$ infinite edges adjacent to $u_i$ by a bijection to the $l(\mu_i)$ points in $\phi_i^{-1}(0)$.

Obviously, there exists a degree $d$ tame morphism $\phi: \cC_{s,d} \rightarrow \cC$ of curve complexes over $k$
which sends $u_i$ to $v_i$, and has
degrees given by integers in $\mu_i$ above the infinite edge of
$\Gamma$ adjacent to $v_i$, for $i=1,\dots, s$, and 
$\phi_{u_i}=\phi_i$ (see Figure \ref{fig:banana}). According to Proposition \ref{prop:lifting}, the map $\phi$
lifts to a tame morphism of smooth proper curves  $\phi_K : X \to X'$ over $K$ the completion of the algebraic 
closure of  $k[[t]]$. The map $\phi_K$ has partial
ramification profile $\mu_1,\dots, \mu_s$. To deduce now the
non-emptiness of $\mathcal A_{g',g}^d(\mu_1,\dots, \mu_s)$, we note that
there exists a subring $R$ of $K$ finitely presented over $k$ such
that  the map $\phi_K$ descends to a finite morphism $\phi_R : \frak
X \to \frak X'$ between smooth proper curves over $\Spec( R )$. In
addition, over a non-empty open subset $U$ of $\Spec( R )$, $\phi_R$
specializes to a tame cover with the same ramification profile as
$\phi_K$.  Since $U$ contains a $k$-rational point, we infer the
existence of a  large  enough $g'$ such that $\mathcal
A_{g',g}^d(\mu_1,\dots, \mu_s) \neq \emptyset$. 

\medskip

We are thus led to consider the case where $s=1$, $g=0$, $\mu = (d_1,d_2, \dots, d_l)$ with $\sum_i d_i =d$, $d_1, \dots, d_t >1$, and $d_{t+1} = \dots = d_l=1$. 
Figure~\ref{fig:ReductionToCyclic} shows that, just as in the previous reduction, one can reduce to the case where $s=1$ and 
 $\mu_1=\{d\}$ with $(d,p)=1$.   (Note that in Figure~\ref{fig:ReductionToCyclic}(a) the degree of the morphism at some of the middle vertices is two;
 Figure~\ref{fig:ReductionToCyclic}(b) is arranged so that the degrees are all odd.)
 But this is just non-emptiness of $\mathcal A_{0,0}((d))$
(see Example \ref{ex:non-vanishing h}).
\qed

\begin{rem} 
As the above proof shows, when $k$ has characteristic zero one can get an explicit upper bound on the smallest positive integer
$g'$ with $\mathcal H_{g',0}^d(\mu_1,\dots, \mu_s) 
\neq \emptyset$. 
Indeed, the permutation $\rho(c_1)\dots\rho(c_s)\tau_1\dots\tau_d$ can be written as the product of 
$d+\sum_{i=1}^s (d-l(\mu_i))$ transpositions.  So without loss of generality we have
$R -d = d+\sum_{i=1}^s (d-l(\mu_i))$, which means that one can take $g'$ to be $1 +  \sum_{i=1}^r(d-l(\mu_i))$.
For $g\geq 1$, $\mathcal H_{g',g}$ is non-empty as soon as  $R$ is non-negative,
which means in this case that one can take $g'$ to be $1+(g-1)d+\frac 12\sum_i (d-l(\mu_i))$.
\end{rem}

\paragraph[Lifting polynomial-like harmonic morphisms of trees]\label{par:polynomial like}
There is a special case of Theorem~\ref{Hurwitz:nonvanishing} in which one does not need to increase the genus of the source curve.
To state the result, we say (following \cite{DeMarco-McMullen}) that a degree $d$ finite harmonic morphism $\bar\phi : T' \to T$ of metric trees is {\em polynomial-like} if there exists an infinite vertex of $T'$ with local degree equal to $d$. 

\begin{thm} \label{thm:polynomial like}
Assume that the residue characteristic of $K$ is zero or bigger than $d$.
Let $\bar\phi : T' \to T$ be a generically {\'e}tale polynomial-like harmonic morphism of metric trees.
Then there exists a degree $d$ polynomial map $\phi : \PP^1 \to \PP^1$ over $K$ lifting $\bar\phi$.
\end{thm}

\pf
It suffices to prove that $\bar\phi$ can be extended to a degree $d$ harmonic morphism of genus zero metrized complexes of curves.
By Theorem~\ref{thm:lifting}, Proposition~\ref{prop:lift augmented to complex}, and Remark~\ref{rem:large characteristic}, 
this reduces to showing that the Hurwitz numbers given by the ramification profiles around
each finite vertex of $T'$ are all non-zero.  Fix an infinite vertex $\infty$ of $T'$ with local degree $d$.  Then it is easy to see that for any such vertex $v'$, the local degree of $\bar\phi$ at $v'$ is equal to the local degree of $\bar\phi$ in the tangent direction corresponding to the unique path from $v'$ to $\infty$.
(This is analogous to \cite[Lemma 2.3]{DeMarco-McMullen}.)
The result now follows from Example~\ref{eg:DM}.
\qed

\paragraph[Lifting of harmonic morphisms in the case the base has
  genus zero] \label{par:target genus zero}
We now consider the special case where $\Gamma$
has genus zero and present more refined lifting results in this case.  
As explained in~\parref{par:def trop curves}, a given
 harmonic morphism
of (augmented) metric graphs does not necessarily have a tropical
modification which is finite. 
We present below a weakened notion of finiteness of a harmonic morphism,
and prove that any harmonic morphism from an (augmented)
metric graph to an (augmented) rational metric graph
 satisfies this weak finiteness property. We discuss in Section~\ref{sec:applications}
some consequences concerning linear equivalence of divisors on metric graphs.

\begin{defn}
A  harmonic morphism 
$\phi :\Gamma\to T$ from an augmented
metric graph $\Gamma$ to a metric tree $T$ 
is said to admit a \emph{weak resolution} if there
exists
a tropical modification $\tau:\widetilde \Gamma\to\Gamma$ and an augmented
harmonic morphism $\widetilde \phi:\widetilde\Gamma\to T$ such that the restriction 
$\widetilde \phi_{|\Gamma}$ is equal to $\phi$, and some tropical
modification of $\widetilde \phi$ is finite.
\end{defn}

In other words, the  morphism $\phi$ has a weak resolution if it
can be extended, up to increasing the degree of $\phi$ using the
modification $\tau$, to a tropical morphism $\widetilde \phi : \widetilde \Gamma \rightarrow T$.

\begin{eg}\label{ex:weakres}
The harmonic morphism depicted in Figure \ref{fig:ex morph augmented}(c)
with $d=1$ can be
weakly resolved by the harmonic morphisms depicted in Figures~\ref{fig:ex trop morphism}(b) and 
\ref{fig:ex morph augmented}(b).
Another example of a weak resolution is depicted in Figure~\ref{fig:res2}.
\begin{figure}[h]
\begin{tabular}{ccc}
\includegraphics[width=5cm,  angle=0]{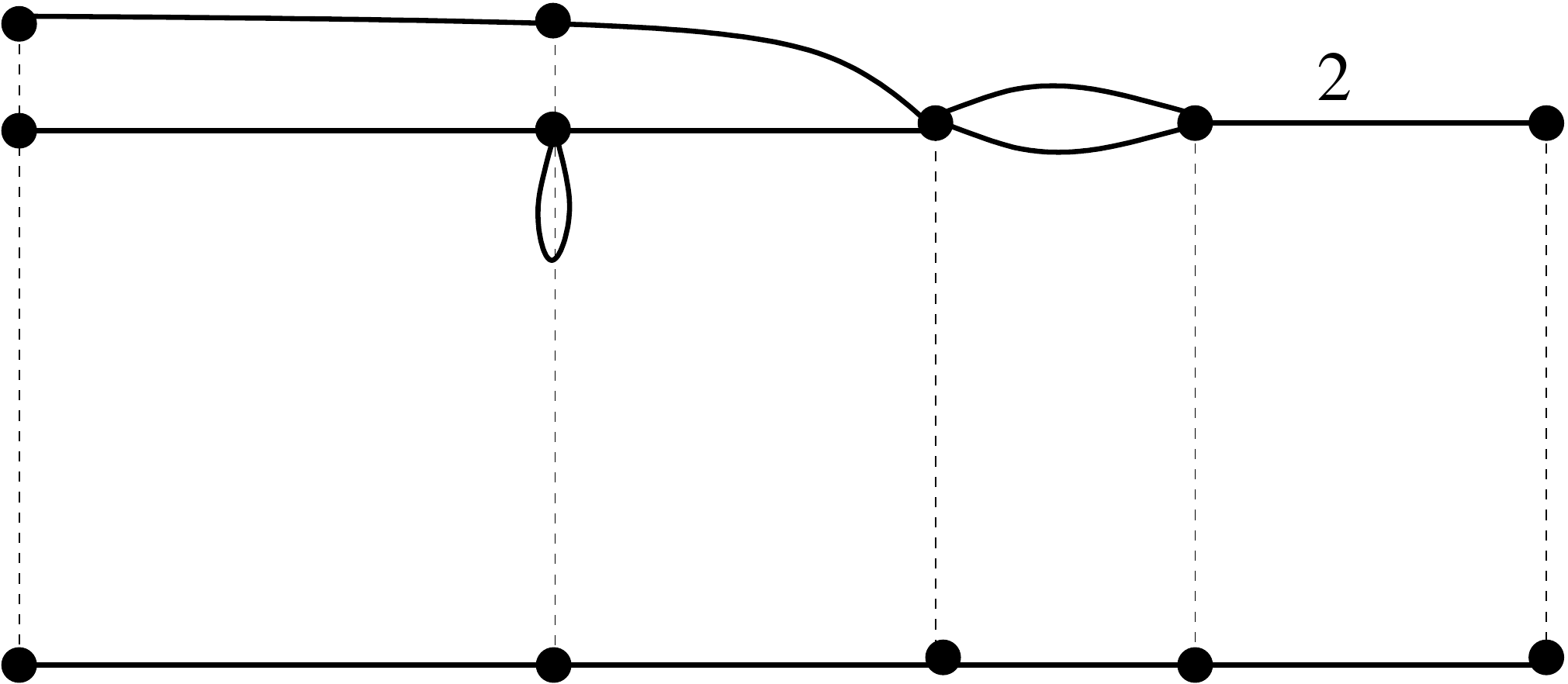}&  \hspace{3ex}&
\includegraphics[width=5cm, angle=0]{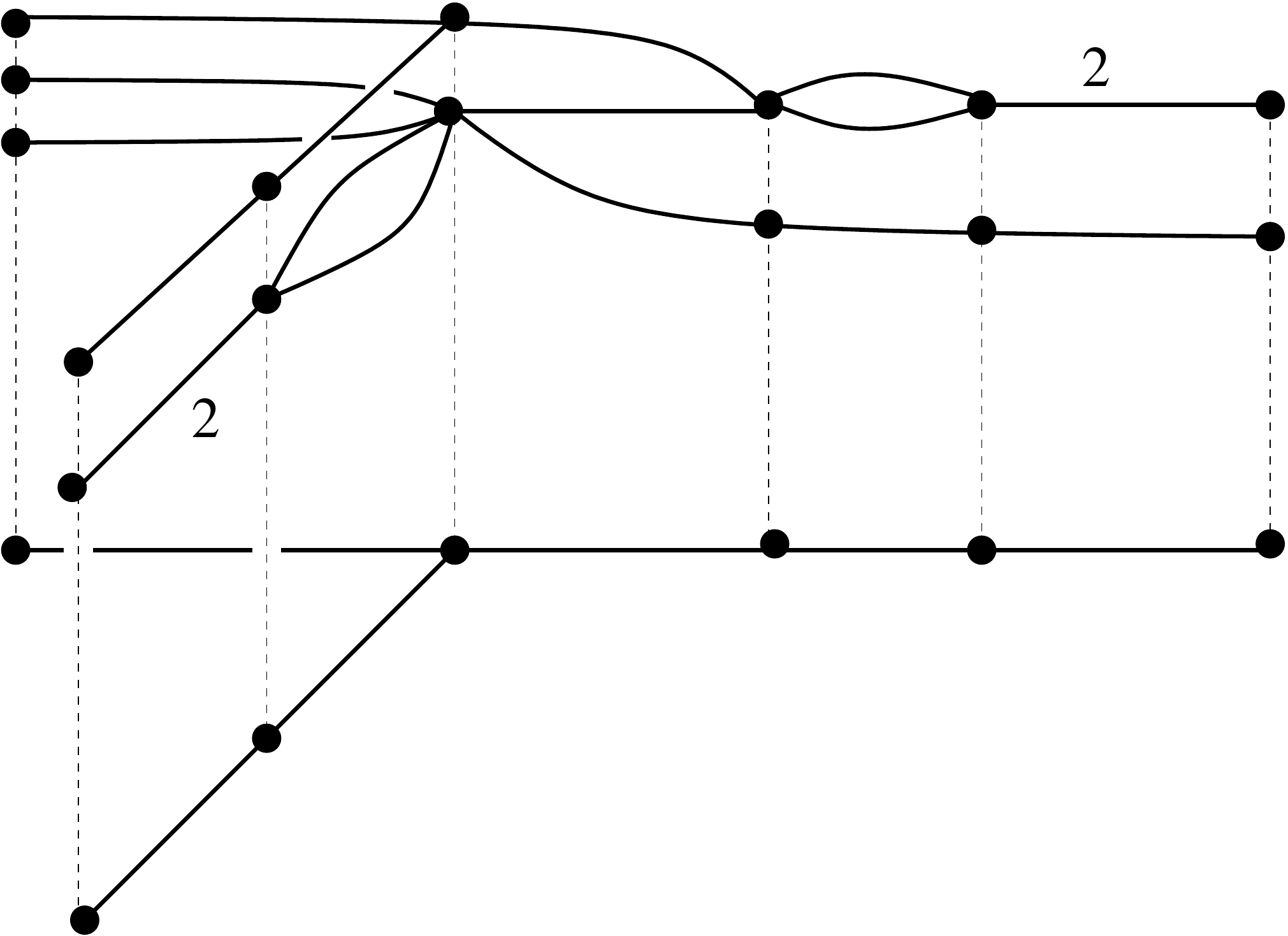} 
 \\a) A  harmonic morphism not tropically && b) A weak resolution of 
\\equivalent to any finite harmonic morphism && the morphism in Figure~\ref{fig:res2}(a)
\end{tabular}
\caption{}\label{fig:res2} 
\end{figure}
\end{eg}

\begin{defn}
Let $\phi :\Gamma\to T$ be a  harmonic morphism from a
metric graph $\Gamma$ to a metric tree $T$.
A point $p\in\Gamma$ is \emph{regular} if $\phi$ is non-constant on all
neighborhoods of $p$.

 The \emph{contracted set} of $\phi$, denoted by $\mathcal E(\phi)$, is
the set of all non-regular points of $\phi$. A \emph{contracted
  component} of $\phi$ is a connected component of $\mathcal E(\phi)$.
\end{defn}

\medskip
The next proposition, together with Proposition
\ref{prop:lifting}, allows us to conclude that any harmonic morphism from an augmented metric 
graph to a metric tree can be realized, up to weak resolutions, as the induced morphism on skeleta of a finite morphism of
triangulated punctured curves.  
Recall that $\Lambda = \val(K^\times)$ is
divisible since $K$ is algebraically closed.

 \begin{prop}[Weak resolution of contractions]\label{prop:liftinggenuszero}
Let $\phi : \Gamma \rightarrow T$ be a harmonic morphism of degree $d$ from a metric
    graph $\Gamma$ to a metric tree $T$.
 \begin{enumerate}
  \item  There exist tropical modifications $\tau: \widetilde\Gamma\to\Gamma$ and $\tau' :\widetilde T\to T$, 
and a harmonic morphism of metric graphs (of degree $\widetilde d\geq d$) 
$\widetilde \phi :\widetilde\Gamma \rightarrow \widetilde T$, such that 
$\widetilde\phi_{|\Gamma\setminus\mathcal E(\phi)} = \phi$, where $\mathcal E(\phi)$ is the contracted part of $\Gamma$.

\item Suppose in addition that $\Gamma$ is augmented, and if $p>0$ that all the non-zero degrees of $\phi$ 
along tangent directions at $\Gamma$ are prime to $p$.  Then there exist tropical modifications of 
$\Gamma$, $T$, and $\phi$ as above such that $\widetilde \phi$ is tame and, in addition, 
there exists a tame harmonic morphism of metrized complexes of $k$-curves
with  $\widetilde\phi$ as 
the underlying finite harmonic
  morphism of augmented metric graphs.
 \end{enumerate}

\end{prop} 
  
\pf
Up to tropical modifications, we may assume that all $1$-valent vertices
of $T$ are infinite vertices.

The proof of (1) goes by giving an algorithm to exhibit a weak resolution of
$\phi$. Note that this algorithm does not produce the weak resolutions
presented in Example \ref{ex:weakres}, since in these cases we could find simpler ones.

\begin{figure}[h]
\begin{tabular}{ccc}
\includegraphics[width=3cm,  angle=0]{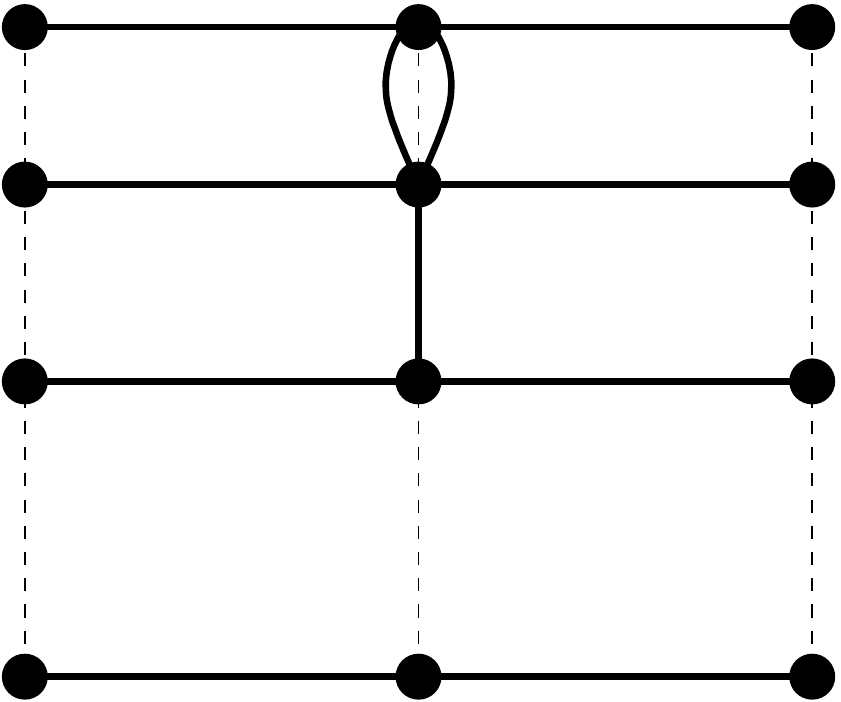}&  \hspace{5ex}
&
\includegraphics[width=3cm, angle=0]{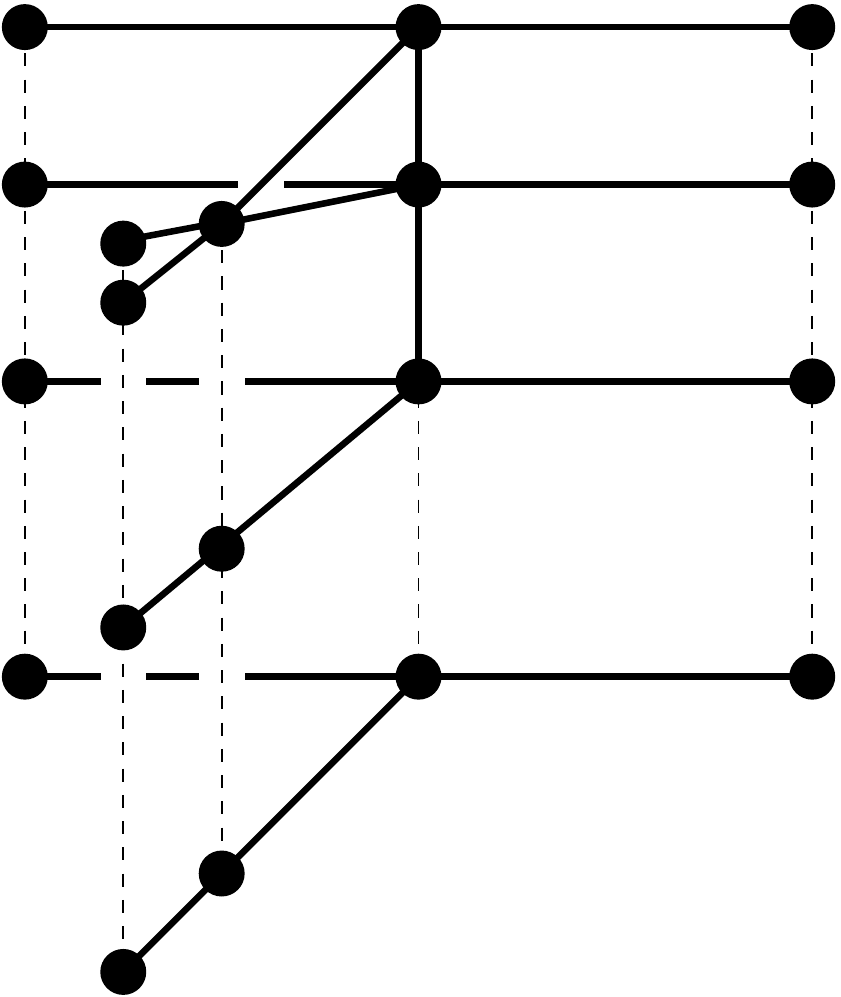}  
 \\a) The morphism $\widetilde \phi$  && b)
 The morphism $\psi_e$ 
\end{tabular}
\caption{The harmonic morphisms $\widetilde \phi$ and $\psi_e$ in the case of Figure
 \ref{fig:ex morph augmented}(c) with $d=1$}\label{fig:resproof}
\end{figure}

Let $V(\Gamma)$ be any strongly semi-stable vertex set of $\Gamma$.
We denote by $d$ the degree of $\phi$, and by $\alpha$ the number of
non-regular vertices of $\phi$.  Given $v$  a non-regular vertex of $\Gamma$, we consider the
tropical 
modification $\tau_v: \widetilde \Gamma_v\to \Gamma $ such that
$\left( \widetilde \Gamma_v\setminus \Gamma \right)\cup\{v\}$ is
isomorphic to $T$ as a metric graph. Considering all those modifications for all
non-regular vertices of $\phi$, we obtain a modification 
$\tau: \widetilde \Gamma\to \Gamma $. We can naturally extend
$\phi$ to a harmonic morphism $\widetilde \phi :\widetilde \Gamma\to
T$ of degree $d+\alpha$ such that $\widetilde \phi_{|\Gamma}=\phi$ and
all degrees of $\widetilde \phi$ on edges not in $\Gamma$
are equal to $1$ (see Figure \ref{fig:resproof}(a) in the case of the
harmonic morphism depicted in Figure \ref{fig:ex morph augmented}(c)
with $d=1$).

By construction, any contracted component of $\widetilde \phi$ is an open edge
of $\Gamma$, and this can be easily resolved. 
Indeed, if $e$ is a finite  contracted edge of 
 $\widetilde \phi$, we  do the following (see Figure
 \ref{fig:resproof}(b)):    
\begin{itemize}
\item consider the tropical modification $\tau_T:\widetilde T\to T$ of $T$ at
  $\widetilde \phi(e)$; denote by $e_1$ the new end of $\widetilde T$;

\item consider $\tau_e:\widetilde
  \Gamma_e\to\widetilde  \Gamma$ the composition of two elementary
tropical modifications  of $\widetilde  \Gamma$ at the middle
  of the edge $e$; denote by $e_2$ and $e_3$ the two new infinite edges of 
$\widetilde  \Gamma_e$, and by $e_4$ and $e_5$ the two new finite
  edges of $\widetilde  \Gamma_e$;
\item subdivide $e_1$ into a finite edge $e_1^0$ of length equal to the lengths of $e_4$ and $e_5$, 
and an infinite edge $e_1^\infty$;
\item consider the morphism of metric graphs $\widetilde \phi_e:
  \widetilde \Gamma_e\to \widetilde T$ defined by  
\begin{itemize}
\item $\widetilde
  \phi_{e|\widetilde \Gamma\setminus \{e_2,e_3,e_4,e_5\}}=\widetilde
  \phi$,
\item $ \widetilde  \phi_e(e_2) = \widetilde \phi_e(e_3) = e_1^\infty$, and $\widetilde \phi_e(e_4) = \widetilde \phi_e(e_5)=e_1^0$,
\item $d_{e_i}(\widetilde  \phi_e)=1$ for $i=2,3,4,5$.
\end{itemize}
\item extend $\widetilde  \phi_e$ to a harmonic morphism of metric
  graphs $\psi_e:\Gamma'\to \widetilde T$, where $\Gamma'$ is a
  modification of $\widetilde \Gamma_e$ at regular vertices in 
$\widetilde  \phi_e^{-1}(\widetilde  \phi(e))$, with 
all degrees of $\widetilde \phi$ on edges not in $\widetilde
\Gamma_e$
 equal to $1$.
\end{itemize}
We resolve in the same way 
  a contracted infinite end of
 $\widetilde \Gamma$.
By applying this process to all contracted edges of $\widetilde \phi$,
we end up with a finite harmonic morphism of metric graphs which is a tropical
modification of $\widetilde \phi$. 

\medskip

Note that in the proof of (1) we increased some of the local degrees by one, but we could have increased these local 
degrees by any amount by inserting an arbitrary number of copies of $T$ in the construction of $\tilde \Gamma$. 
Based on this remark, the proof of (2) now follows the same steps as the proof of (1), using in
addition the following claim: 

\medskip

\noindent{\bf Claim.} 
Let $g'\ge 0$ and $d,s>0$ be integers. 
Let $\mu_1,\dots, \mu_s$ be a collection of $s$ tame partitions of $d$. 
Then there exist arbitrarily large non-negative integers $d'$
 such that $\mathcal A_{g',0}^{d'}(\mu'_1,\ldots,\mu'_s)$ is non-empty,
 where
$\mu'_i$ is the partition of $d'$ obtained by adding a sequence of
 $d'-d$ numbers $1$ to each partition $\mu_i$.
 
\medskip

\begin{figure}[h]
\scalebox{0.32}{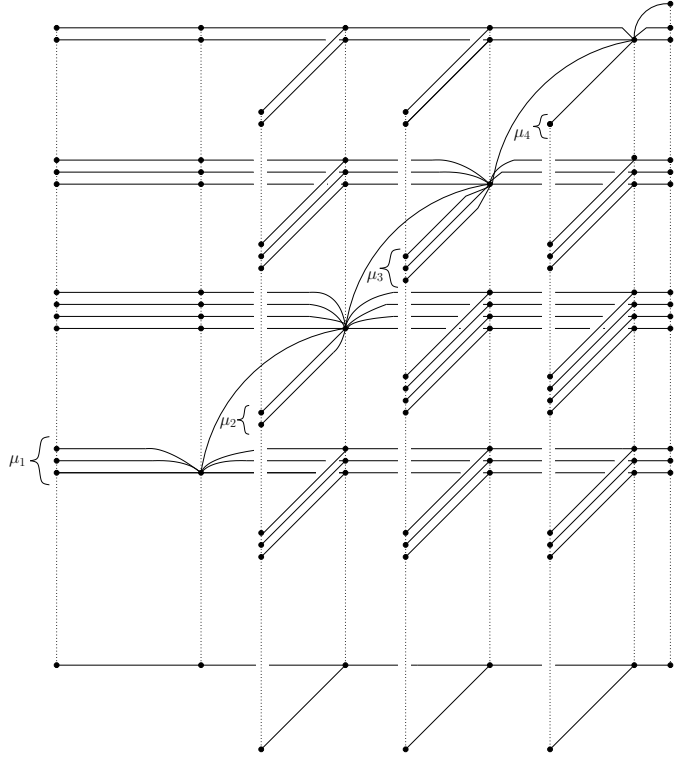}
\caption{Reduction to the case $s=1$ in the proof of (2) in Proposition~\ref{prop:liftinggenuszero}.
  Degrees on (infinite) edges related to $\mu_i$ are exactly the
  integers appearing in $\mu_i$.  All the other degrees are
  one. Degrees over each infinite edge consist of a
  $\mu_i$ and precisely $(s-1)d$ numbers
  1.}\label{fig:reductiontoonemu2} 
\end{figure}

 Figure~\ref{fig:reductiontoonemu2},
  Figure~\ref{fig:ReductionToCyclic}(a), our resolution procedure, and
  the argument used for the positive characteristic case of the 
  proof of Theorem~\ref{Hurwitz:nonvanishing} reduce the proof of the
  claim to the case $s=1$ and $\mu_1 = \{d\}$ with $(d,p)=1$. But in
  this case, for any $g'\geq 0$, by the group theoretic method we used in the proof of 
  Theorem~\ref{Hurwitz:nonvanishing},
  there exists a (tame) covering of $\P^1$ by a curve of genus
  $g'$ having (tame) monodromy group the cyclic group $\Z/d\Z$, and with the property that the ramification 
  profile above the point 1 of $\P^1$ is given by $\mu= \{d\}$. This
  finishes the proof of the claim, and the proposition follows. 
\qed

\section{Applications}\label{sec:applications}

\paragraph[Linear equivalence of divisors]\label{par:rank} 
A {\em (tropical) rational function} on a metric graph $\Gamma$
is a continuous piecewise affine function $F:\Gamma\rightarrow \R$ with integer slopes. 
If $F$ is a rational function on $\Gamma$, $\mathrm{div}(F)$ is the divisor on $\Gamma$ 
whose coefficient at a point $x$ of $\Gamma$ is given by $\sum_{v\in T_x } d_v F$, 
where the sum is over all tangent directions to $\Gamma$ at $x$ and $d_v F$ is the outgoing 
slope of $F$ at $x$ in the direction $v$. 
Two divisors $D$ and $D'$ on a metric graph $\Gamma$ are called {\em
  linearly equivalent} if there exists a rational function $F$ on
$\Gamma$ such that $D - D' = \mathrm{div}(F)$, in which case we write $D\sim D'$.
 For a divisor $D$ on $\Gamma$, the {\em complete linear system of $D$}, denoted $|D|$, is the set of all effective divisors $E$ linearly equivalent to $D$. 
The {\em rank} of a divisor $D \in \mathrm{Div}(\Gamma)$ is defined to be
 \[r_\Gamma (D):= \min_{\substack{E\::\:E\geq 0\\|D-E|=\emptyset}}\deg(E)-1.\]

\medskip

Let $\phi: \Gamma \rightarrow T$ be a finite harmonic morphism from $\Gamma$ to a metric tree  $T$ of  degree $d$. 
For any point $x\in T$, the (local degree of $\phi$ at the points of the) fiber $\phi^{-1}(x)$ defines a divisor of degree $d$ in $\mathrm{Div}(\Gamma)$ that we denote by $D_x(\phi)$. We have
\[D_x(\phi) := \sum_{y\in \phi^{-1}(x)} d_y(\phi) (y),\]
where $d_y(\phi)$ denotes the local degree of $\phi$ at $y$.

\begin{prop}\label{prop:rank-finite} Let $\phi: \Gamma \rightarrow T$ be a finite harmonic morphism of degree $d$ from $\Gamma$ to a metric tree. 
Then for any two points $x_1$ and $x_2$ in $T$, we have $D_{x_1}(\phi) \sim D_{x_2}(\phi)$ in $\Gamma$.  Moreover, for every $x \in T$ the rank of the 
divisor $D_x(\phi)$ is at least one. 
\end{prop}

\pf 
Since $T$ is connected, we may assume that $x_1$ and $x_2$ are sufficiently close; more precisely,
we can suppose that $x_2$ lies on the same edge as $x_1$ with respect to some model $G$ for $\Gamma$.
Removing the open segment $(x_1,x_2)$ from $T$ leaves two 
connected components $T_{x_1}$ and $T_{x_2}$ which contain $x_1$ and $x_2$, respectively.
 Identifying the segment $[x_1,x_2]$ with the interval $[0, \ell]$ by a linear map (where $\ell=\ell([x_1,x_2])$ denotes the length in $T$ of the segment $[x_1,x_2]$) 
 gives a rational function $F: \Gamma \rightarrow [0, \ell]$ by sending $\phi^{-1}(T_{x_1})$ and $\phi^{-1}(T_{x_2})$ to $0$ and $\ell$, respectively. 
It is easy to verify that $D_{x_1}(\phi) - D_{x_2}(\phi) = \mathrm{div}(F)$, which establishes the first part.
 
The second part follows from the first, since $y$ belongs to the support of the divisor $D_{\phi(y)}(\phi) \sim D_x(\phi)$ for all $y \in \Gamma$, 
which shows that $r_\Gamma(D_x(\phi)) \geq 1$.
\qed

\medskip

By Theorem~\ref{thm:lifting harm},
any finite morphism $\phi:\Gamma \rightarrow T$ 
can be lifted to a morphism $\phi: X \rightarrow \P^1$ of smooth proper curves, possibly with $g(X) > g(\Gamma)$. 
This shows that any effective divisor 
on $\Gamma$ which appears as a fiber of a finite morphism to a metric tree
can be lifted to a divisor of rank
at least
 one on a smooth proper curve {\em of possibly higher genus}.

  We are now going to show that the (additive) equivalence relation generated by fibers of
  ``tropicalization'' of 
finite morphisms $X \rightarrow \P^1$ 
coincides with tropical linear equivalence of divisors. 
To give a more precise statement, let $\Gamma$ be a metric graph with first Betti number $h_1(\Gamma)$,
 and consider the family of all smooth proper curves of genus $h_1(\Gamma)$ 
 over $K$ which admit a semistable vertex set $V$ and a finite set of $K$-points $D$ such that 
 the metric graph 
 $\Sigma(X, V\cup D)$ is a modification of $\Gamma$. Given such a curve $X$ 
 and a finite morphism $\phi: X \rightarrow \P^1$, there is a corresponding finite harmonic  morphism 
 $\phi: \Sigma(X, V\cup D) \rightarrow T$ from a modification of $\Gamma$ to a metric tree $T$.
   Two effective divisors $D_0$ and $D_1$ on $\Gamma$ are called {\it
     strongly effectively linearly equivalent} if there exists a
   morphism 
   $\phi: \Sigma(X, V \cup D) \rightarrow T$ 
as above such that $D_0  = \tau_* (D_{x_0}(\phi))$ and
$D_1=\tau_*(D_{x_1}(\phi))$ for two points $x_0$ and $x_1$ in
$T$. Here 
$\tau_*: \mathrm{Div} (\Sigma(X, V \cup D)) \rightarrow
\mathrm{Div}(\Gamma)$ 
is the extension 
 by linearity  of the retraction map $\tau:\Sigma(X, V\cup D)  \rightarrow \Gamma$.
 The equivalence relation on the abelian group $\Div(\Gamma)$ generated by this relation is called {\it effective linear equivalence} of divisors.  
In other words, two divisors $D_0$ and $D_1$ on $\Gamma$ are effectively linearly equivalent
if and only if there exists  an effective divisor $E$ on $\Gamma$ such that $D_0+E$ and $D_1+E$ are strongly effectively linearly equivalent. This can be summarized as follows:
$D_0$ and $D_1$ on $\Gamma$ are effectively linearly equivalent if and only if
there exists a lifting of $\Gamma$ to a smooth proper curve $X/K$ of genus $h_1(\Gamma)$, 
and a finite morphism $\phi: X \rightarrow \P^1$ such that $\tau_*(\phi^{-1}(0)) = D_0 +E$ and 
$\tau_*(\phi^{-1}(\infty))=D_1+E$ for some effective divisor $E$, where
$\tau_*$ is the natural retraction map from $\Div(X)$ to $\Div(\Gamma)$.

  \begin{thm}\label{thm:equ effective/linear}
  The two notions of linear equivalence and effective linear equivalence of divisors on a metric graph $\Gamma$ coincide. 
As a consequence, linear equivalence of divisors is the additive equivalence relation 
generated by (the retraction to $\Gamma$ of) fibers of finite harmonic morphisms from a tropical modification of $\Gamma$ 
to a metric graph of genus zero.
  \end{thm}

\pf Consider two divisors $D_0$ and $D_1$ which are effectively linearly equivalent. There exists an effective divisor $E$
and a finite harmonic morphism 
$\phi: \widetilde \Gamma \to T$, from a tropical modification of
$\Gamma$
to a metric tree, 
such that $D_0+E = D_{x_0}(\phi)$ and $D_1+E = D_{x_1}(\phi)$ for two points $x_0,x_1\in T$. By Proposition~\ref{prop:rank-finite}
we have $D_0+E \sim D_1+E$, which implies that $D_0$ and $D_1$ are
linearly equivalent in 
in $\widetilde \Gamma$, and hence in $\Gamma$.

To prove the other direction,  it will be enough to show that if $D$ is linearly equivalent to zero,
then there exists an effective divisor $E$ such that $D+E$ and $E$ are fibers of a finite harmonic morphism $\phi$ from 
a modification of $\Gamma$ to a metric tree $T$, and such that $\phi$ 
can be lifted to a morphism $X \rightarrow \P^1$.   

 By assumption, there exists a rational function $f: \Gamma \to \R\cup\{\pm\infty\}$ such that 
$D+\mathrm{div}(f) =0$.  We claim that there is a tropical
modification $\widetilde \Gamma$ of $\Gamma$ together with an extension of $f$ to 
a (not necessarily finite) harmonic morphism
 $\phi_0: \widetilde \Gamma \rightarrow \R\cup \{\pm\infty\}$. 
The tropical modification $\widetilde \Gamma$ is obtained from $\Gamma$ by choosing a vertex set which contains all the points in the support of $D$, 
adding an infinite edge to any finite vertex
in $\Gamma$ with $\ord_v(f) \neq 0$, and 
extending $f$ as an affine linear function of slope $-\ord_v(f)$ along this infinite edge. 
It is clear that the resulting map $\phi_0$ is harmonic.

Consider now the retraction map $\tau: \widetilde \Gamma \to \Gamma$, and note that for the two divisors 
$D_{\pm \infty}(\phi_0)$, we have 
$\tau_*\bigl(D_{\pm \infty}(\phi_0)\bigr) = D_\pm$, where $D_+$ and $D_-$ denote the positive and negative part of $D$, respectively.
By Proposition~\ref{prop:liftinggenuszero}, 
there exist tropical modifications $\overline \Gamma$ of 
$\widetilde \Gamma$ and $T$ of $\R\cup\{\pm\infty\}$ such that $\phi_0$ extends to a finite harmonic morphism
 $\phi: \overline \Gamma \to T$ which can be lifted to a finite morphism $X \to \P^1$.
If we denote (again) the retraction map $\overline \Gamma \to \Gamma$ by $\tau$, then
$\tau_*\bigl(D_{\pm \infty}(\phi)\bigr) = D_\pm +E_0$ for some effective divisor $E_0$ in $\Gamma$. 
Setting $E = D_{-} + E_0$, the divisors $D+E$ and $E$ are strongly effectively linearly equivalent, and the theorem follows. 
\qed

\medskip

\begin{eg} 
Here is an example which illustrates the distinction between the notions of (effective) linear equivalence 
and strongly effective linear equivalence of divisors, as introduced above.

Let $\Gamma$ be the metric graph depicted in Figure \ref{glasses}(a), with arbitrary lengths, 
and $K_\Gamma = ( p )+ ( q )$ the 
canonical divisor
on $\Gamma$.

\begin{figure}[h]
\begin{tabular}{ccc}
\scalebox{.32}{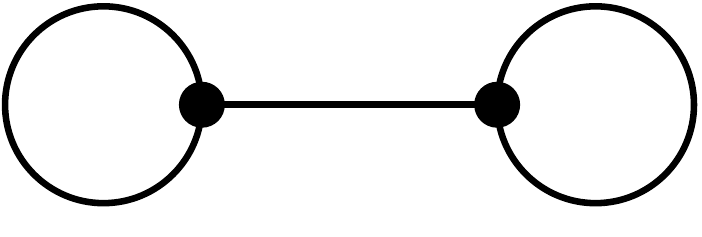}& \hspace{.4cm}
\scalebox{.28}{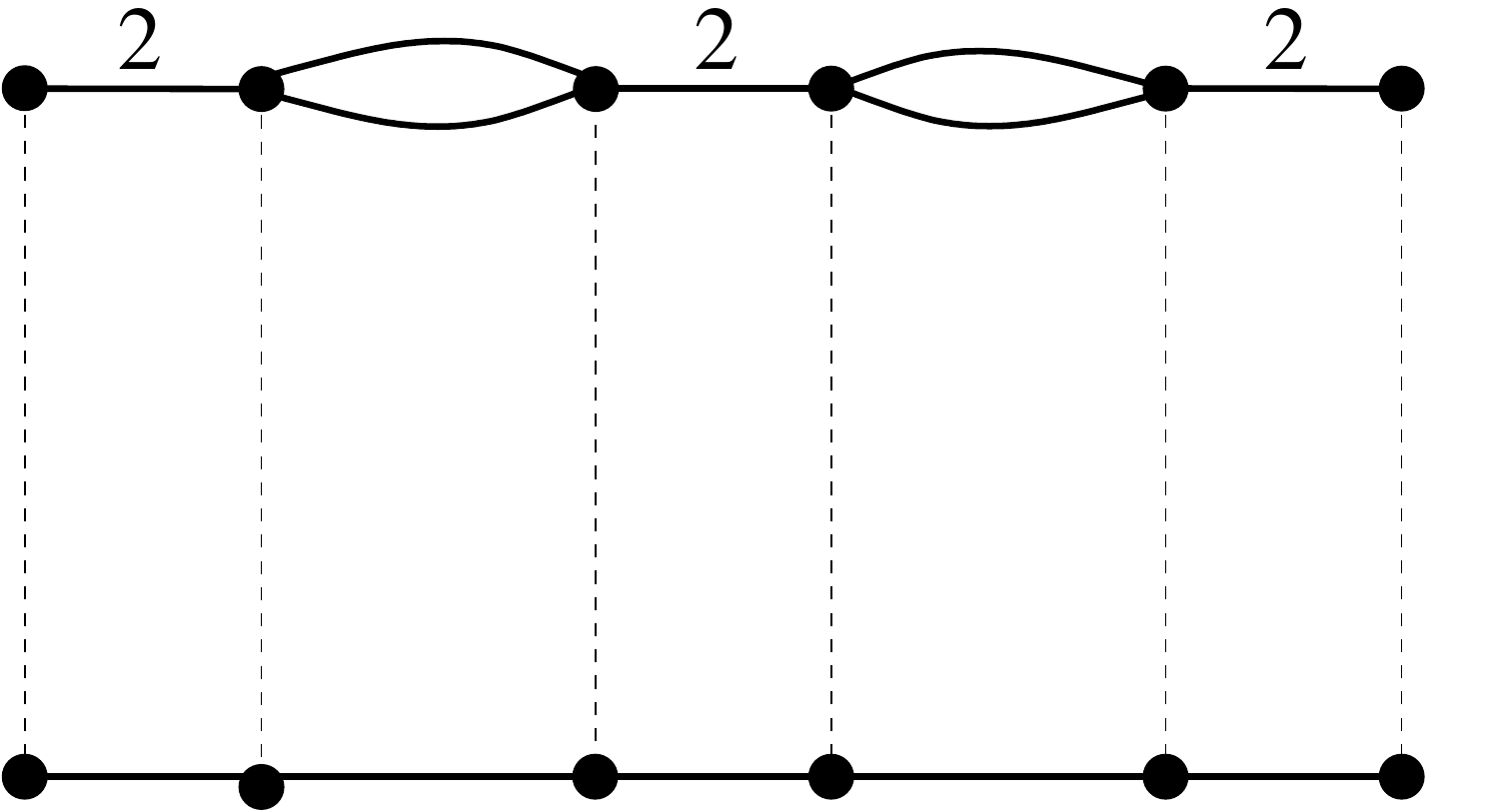}&  \hspace{.4cm}
\scalebox{.24}{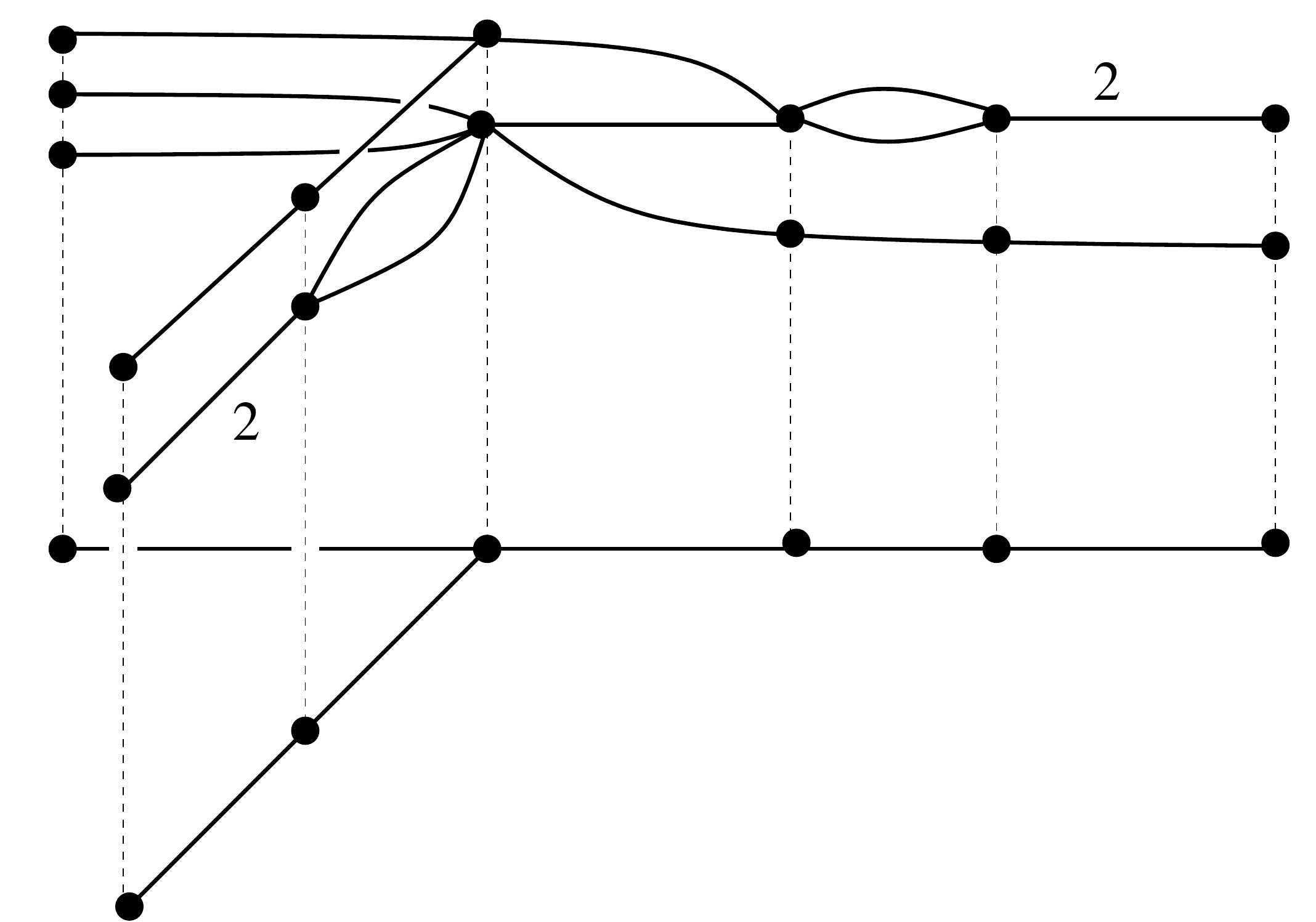}
\\ a) $K_\Gamma=( p )+ ( q )$ & b) An effective lift of $2( t )$
& c) A non-effective lift of $K_\Gamma$
\end{tabular}
\caption{}\label{glasses}
\end{figure}

We claim that $K_\Gamma$ is not the specialization of any effective divisor 
of degree two representing the 
canonical class of a smooth proper curve of genus two over $K$. 
More precisely, we claim that for any triangulated punctured curve $(X, V \cup D)$ over $K$
such that 
$\Sigma(X, V\cup D)$ is a tropical modification of $\Gamma$, and for
any effective divisor $\mathcal D$ in $\Div(X)$ with $K_\Gamma =
\tau_*(\mathcal D)$, we must have $r_X(\mathcal D) = 0$.
(Here $\tau_*$ denotes the specialization map from
$\Div(X)$ to $\Div(\Gamma)$ and $r_X(D) = \dim_K(H^0(X, \mathcal O(\mathcal D)))-1$.)
Indeed, otherwise there would exist a degree $2$ finite harmonic
morphism $\pi : \widetilde \Gamma \to T$ from some tropical
modification of $\Gamma$ to a metric tree
with the property that $\pi(p)=\pi(q)$.  Restricting such a harmonic morphism to the preimage in $\widetilde \Gamma$ of the loop containing $p$ 
would imply, by Proposition~\ref{prop:rank-finite}, that the divisor $(p)$ has rank one in a genus-one metric graph, which is impossible.
On the other hand, Figure \ref{glasses}(b) shows that the divisor $2( t )\sim ( p )+(q ) $ can
be lifted to an effective representative of the canonical class $K_X$, where $t$ is the
middle point of the loop edge with vertex $q$. This shows that the two linearly equivalent divisors $D_0= ( p )+( q )$ and $D_1= 2( t )$
are not strongly effectively linearly equivalent. 

However, $D_0$ and $D_1$ are effectively linearly equivalent.  Indeed,
adding $E=( p )$ to $D_0$ and $D_1$, respectively, gives the two divisors 
$2(p)+(q)$ and $2(t)+(p)$ which are retractions of fibers of a degree $3$ finite
harmonic morphism from a tropical modification of $\Gamma$ to a tree,
as shown in Figure  \ref{glasses}(c).
Consequently, $D_0+(p)$ and $D_1+(p)$ can be lifted to linearly equivalent effective divisors on a smooth proper curve $X$.

Note also that  Figure \ref{glasses}(c) shows that 
since $(p_1)+(p_2) + (q) -(p_3)$ 
can be lifted to  a non-effective representative of the canonical class $K_X$, 
there exists a non-effective divisor $\mathcal D$ in the canonical class $K_X$ of $X$ such that 
$\tau_*( \mathcal D) = ( p ) + (q)$. 
\end{eg}

\paragraph[Tame actions and quotients]\label{par:quotient} 
Let $\cC$ be a metrized
complex of $k$-curves, and denote by $\Gamma$ the underlying metric graph of $\cC$.
An {\it automorphism} of $\cC$ is a (degree one) finite harmonic morphism of metrized complexes $h : \cC \rightarrow
 \cC$ which has an inverse.
  The group of automorphisms of $\cC$ is denoted by
$\mathrm{Aut}(\cC)$. 

Let $H$ be a finite subgroup of $\mathrm{Aut}(\cC)$.   
The action of  $H$ on $\cC$ is {\it generically free} if for any vertex $v$ of $\Gamma$, 
the inertia (stabilizer) group $H_v$ acts freely on an open subset of $C_v$. 
A finite subgroup $H$ of $\mathrm{Aut}(\cC)$ is called {\it tame} if the action of $H$ on $\cC$ is generically free and all the inertia subgroups $H_x$
for $x$ belonging to some $C_v$ are cyclic of the form $\mathbf{Z} / d \mathbf{Z}$ 
for some positive integer $d$, with $(d,p) =1 $ if $p>0$.   In this case we say that the action of $H$ on $\cC$ is tame. 
 
 \begin{rem}~\label{rem:tameactions}
 The stabilizer condition in the definition of tame actions 
 is equivalent to requiring the cover 
$C_v \rightarrow C_v/H_v$ be tame, where $H_v$ is the stablizer of the vertex $v$. 
To see that this latter condition implies all the stablizers of points on 
$C_v$ are cyclic, consider a uniformizer $\pi$ at a point $x$, and consider the map 
$H_x \rightarrow k^{\times}$ which 
sends an element $h\in H_x$ to $h(\pi)/\pi$. This is independent of the choice of the uniformizer, 
and embeds $H_x$ in the subgroup of roots of unity in $k^{\times}$, from which the assertion follows. The other direction is clear from the definition.

\noindent Note that, more generally,  one has a filtration of $H_v$ with higher ramification groups 
$H_v \supseteq H_0 = H_x  \supseteq H_{1} \supseteq H_{2} \supseteq \dots$, the quotient $H_0/H_1$ 
is a finite cyclic group of order prime to the characteristic $p$, and $H_i/H_{i+1}$ are all $p$-groups. In the case of tame actions, $H_1$ is trivial.
\end{rem}

\medskip
 
In this section, we characterize tame group actions $H$ on $\cC$ which lift to an action of $H$ on some smooth proper curve $X/K$ lifting 
$\cC$.
The main problem to consider is whether there exists a refinement $\widetilde \cC$ of $\cC$
and an extension of the action of $H$
to $\widetilde \cC$ such that the quotient $\widetilde \cC/H$ can be defined, 
and such that the projection map $\pi : \widetilde \cC \rightarrow  
\widetilde\cC/H$ is a tame harmonic morphism. 
The lifting of the action of $H$ to a smooth proper curve $X$ as above will then be a consequence of our lifting theorem. 

\paragraph
Let $H$ be a tame group of automorphisms of a metrized complex $\cC$. 
Let $W_H = W_H(\cC)$ be the set of all $w \in \Gamma$ lying in in the middle of
an edge $e$ such that there is an element $h \in H$ having $w$ as an isolated fixed point.
Denote by $H_w$ the stabilizer of $w \in W_H$.  It is easy to see that
$H_w$ consists of all elements $h$ of $H$ which 
restrict on $e$ either to the identity or to the symmetry 
with center $w$. In particular, if $h|_{e}\ne \rm{id}$, then $h$
permutes the two vertices $p$ and $q$ adjacent to $e$.
For $w \in W_H$, the inertia group 
$H_{\red_p(e)} = H_{\red_q(e)} \cong \Z / d_{e} \Z$ (for some integer $d_e$)
is a normal subgroup of index two in $H_w$:
$$ 0 \To H_{\red_p(e)} \To H_w \To \Z/2 \Z \To 0.$$

We make the following assumption on the groups $H_w$: 

\begin{defn}
  A tame group of automorphisms $H$ of a metrized complex $\cC$ satisfies
  the \emph{dihedral condition} provided that, 
  for all $w \in W_H$, the stabilizer group $H_w$ is isomorphic to
  the dihedral group generated by two elements $\sigma$ and $\zeta$ with the
  relations
  $$\sigma^2 = 1,\quad \zeta^{d} =1, \sptxt{ and } \sigma \zeta \sigma = \zeta\inv$$
  for some integer $d$, such that $H_{\red_p(e)} = \angles{\zeta}$. 
\end{defn}

The dihedral condition means that the above short exact sequence splits,
and the action of $\Z/2\Z\cong\{\pm 1\}$ on $H_{\red_p(e)}$ is given by 
$h \rightarrow h^{\pm 1}$ for $h\in H_{\red_p(e)}$.

\medskip

We can now formulate our main theorem on lifting tame group actions:

\begin{thm}~\label{thm:liftinggroupaction} 
Let $H$ be a finite group with a tame action on a metrized complex $\cC$.
\begin{enumerate}
\item If $W_H \neq \emptyset$, then the dihedral condition and $\mathrm{char}(k)
  \neq 2$ are the necessary and sufficient conditions for the existence of
  a refinement $\widetilde \cC$ of $\cC$ such that the action of $H$ on
  $\cC$ extends to a tame action on $\widetilde\cC$ such that 
  $W_H(\td\cC) = \emptyset$.%
\footnote{See \cite[\S{2.3}]{raynaud:cover_specialization} for a related discussion,
  including remarks on the situation in characteristic $2$.}
\item If $W_H = \emptyset$, then the quotient $\cC/H$ exists in the
  category of metrized complexes.  In addition, the action of $H$ on $\cC$
  can be lifted to an action of $H$ on a triangulated punctured $K$-curve
  $(X,V \cup D)$ such that $\Sigma(X, V\cup D_0) \cong \cC$ with
  $D_0\subset D$, the action of $H$ on $X \setminus D$ is \'etale, and the
  inertia group $H_x$ for $x\in D$ coincides with the inertia group
  $H_{\tau(x)}$ of the point $\tau(x) \in \Sigma(X, V\cup D_0) = \cC$.
\end{enumerate}
\end{thm}

\pf
Suppose that $W_H\neq\emptyset$, that the dihedral condition holds, and
that $\chr(k)\neq 2$. Fix an orientation of the edges of $\Gamma$, and for
an oriented edge $e$, denote by $p_0$ and $p_\infty$ the two vertices of
$\Gamma$ 
which form the tail and the head of $e$, respectively. 
Let $w$ be a point  lying in the middle of an oriented edge $e = (p_0,p_\infty)$ of $\Gamma$ which is an isolated fixed point of some elements of $H$. 
Take the refinement $\widetilde \cC$ of $\cC$ obtained by adding all such points $w$ 
to the vertex set of $\Gamma$ and by setting $C_w = \P^1_k$, $\red_e( \{w, p_0\})= 0$, 
and $\red(\{w, p_\infty\}) =\infty$.
To see that the action of $H$ on $\cC$ extends to $\widetilde \cC$, first note that one can define 
a generically free action of $H_w$ on $\P^1_k$ 
(equivalently, one can embed $H_w$ in $\mathrm{Aut}(\P_k^1)$) in a way compatible with the action 
of $H_w$ on $\Gamma$, i.e., such that all the elements of 
$H_{\red_{p_0}(e)} = H_{\red_{p_\infty}(e)}$ fix the two points $0$ and $\infty$ of $\P^1_k$, 
and such that the other elements of $H_w$ permutes the two points $0, \infty \in \P^1_k$. Indeed, 
the dihedral condition is the necessary and sufficient condition for the existence of such an action. 
Under this condition and upon a choice of a 
$d_{e} = |H_{\red_{p_0}(e)}|$-th root of unity $\zeta_{d_e}\in k$,  
and upon the choice of the point $1\in \P^1_k$ as a fixed point of $\sigma$, 
the actions of the two generators $\sigma$ and $\zeta$ of $H_w$ on $\P^1$ are given 
by $\sigma(z) = 1/z$ and $\zeta(z) = \zeta_{d_e} z$, respectively. 

Fix once for all a $d$-th root of unity $\zeta_d \in k$ for each positive
integer $d$ (with $(d,p)=1$ in the case $p>0$).  Given $h\in H$, we extend
the action of $h$ on $\cC$ to an action on $\td\cC$ 
in the following way.  Let $w\in W_H(\cC)$ and let $e$ be the edge
containing $w$, with the orientation chosen above.  If $h(w) \neq w$, we define 
$h_w: C_w\to C_{h(w)}$ by $h_w ={\rm id}_{\P^1}$ if $h$ is compatible with
the orientations of $e$ and $h(e)$, and we set $h_w(z) = z\inv$ otherwise.
If $h\in H_w$, we define the action of $h$ on $C_w$ as above.  This
defines a generically free 
action of $H$ on $\widetilde \cC$.  The inertia groups of the points
$0,\infty,$ and $\pm 1$ in $C_w$ are $\Z/d_e\Z$, $\Z/d_e\Z$, and $\Z/2
\Z$, respectively.  Since $p \neq 2$, this shows that the action of $H$ on
$\widetilde \cC$ is tame.  By construction we have 
$W_H(\td\cC) = \emptyset$.

Working backward, one recovers the necessity of the dihedral condition and
$\chr(k)\neq 2$.  Indeed, any $\td\cC$ satisfying the conditions of
the theorem must contain each $w\in W_H(\cC)$ as a vertex.  Since $H_w$
acts on $\bP^1_k$ in the manner described above, it must be a dihedral
group; since its action on $C_w$ has stabilizers of order $\pm 2$, we must
have $\chr(k)\neq 2$.

Now we assume that the action of $H$ on $\cC$ is tame and 
that no element of $H$ has an isolated fixed point in the middle of an edge.
We will define the quotient metrized complex $\cC/ H$.
The metric graph underlying $\cC /H$ is the quotient graph
$\Gamma/H$ equipped with the following metric: given an edge $e$ of
$\Gamma$ of length $\ell$ and stabilizer $H_e$, we define the length of its
projection in $\Gamma/H$ to be $\ell \cdot |H_e|$.
The projection map $\Gamma \rightarrow \Gamma/H$ is a tame finite harmonic morphism.

For any vertex $p$ of $\Gamma$, the $k$-curve associated to its image
 in $\cC /H$ is $C_p/H_p$. 
The marked points of $C_p/H_p$ are the different orbits of the marked points of $C_p$, 
and are naturally in bijection with the edges of $\Gamma/H$ adjacent
to the projection of $p$. The projection map
$\cC \rightarrow \cC/H$ is a tame harmonic morphism of metrized complexes.

To see the second part, let $\cC'$ be the (tropical) modification of $\cC$
obtained as follows: for any closed point  
$x \in C_p$ with a non-trivial inertia group and which is not the
reduction $\red_p(e)$  of any edge $e$ adjacent to $p$, 
consider the elementary tropical modification of $\cC$ at $x$. 
Extend the action of $H$ to a tame action on $\cC'$ by defining $h_x: e_x \rightarrow e_{h(x)}$ to be affine with slope one 
for any such point. Let $\pi :\cC' \rightarrow \cC'/H$ be the projection map. 
Let $(X',V' \cup D')$ be a triangulated punctured $K$-curve such that 
$\cC(X', V' \cup D') \cong \cC'/H$.  By Theorem~\ref{thm:lifting2}, 
the tame harmonic morphism $\pi$ lifts to a morphism of 
triangulated punctured $K$-curves $(X, V \cup D) \rightarrow (X', V' \cup D')$.
By Remark~\ref{rem:auto}, we have an injection 
$\iota: \Aut_{X'}(X) \hookrightarrow \Aut_{\cC'/H}(\cC')$. 
By the construction given in the proof of Theorem~\ref{thm:lifting2}, it is easy to see that every $h\in H$ lies 
in the image of $\iota$, and thus $H \subset \Aut_{X'}(X)$. The last part follows formally from 
the definition of the modification $\cC'$
and the choice of $X$ as the lifting of $\pi: \cC' \rightarrow \cC'/H$. 
\qed

\begin{rem}[Compare with Remark~\ref{rem:tameactions}]\label{rem:tameactionsbis}
If the characteristic of $k$ is positive, the lifting of the action of a finite group on a metrized complexes cannot be guaranteed in general
without further assumptions. Indeed, even in the smooth case, i.e., where the metrized complex 
consists of a single vertex $v$ and a single curve $C_v$, 
there are obstructions to the lifting~\cite{Oort, OSS97, GM98, BM00}, e.g., due to the fact that the automorphism group of a smooth proper curve in positive 
characteristic does not respect the Hurwitz upper bound $84(g-1)$. 
However, Pop's recent proof of the Oort conjecture~\cite{Pop14}, based on the results of Obus and Wewer~\cite{OW14}, 
shows that in the smooth case, 
the action can be lifted under the 
assumption that the stablizers of points are all cyclic. A natural question is then to see whether our theorem 
can be extended by 
only requiring  all the stablizers of points to be cyclic (without the tameness assumption).
\end{rem}

\paragraph[Characterization of liftable hyperelliptic augmented metric graphs]\label{par:hyper}
Let $\Gamma$ be an augmented metric graph and denote by $r^\#$ the weighted rank function on divisors introduced in~\cite{AC11}.
Recall that this is the rank function on the non-augmented metric graph $\Gamma^\#$
obtained from $\Gamma$ by attaching $g(p)$ cycles, called {\it virtual cycles}, of (arbitrary) positive lengths to each $p \in \Gamma$ with $g(p)>0$.
We say that an augmented metric graph $\Gamma$ is {\it hyperelliptic} if $g(\Gamma) \geq 2$ and there exists a divisor $D$ in $\Gamma$ of degree two 
such that $r_\Gamma^\#(D)=1$. 
An augmented metric graph is said to be \emph{minimal} if it contains neither 
infinite vertices nor $1$-valent vertices of genus $0$. 
Every augmented metric graph $\Gamma$ is tropically equivalent to a minimal
 augmented metric graph $\Gamma'$, which is furthermore
 unique if $g(\Gamma)\ge 2$. Since the tropical rank and weighted rank
 functions are invariant under tropical modifications, an augmented
 metric graph $\Gamma$ is hyperelliptic if and only if $
 \Gamma'$ is. Hence we restrict in this section to the case of minimal
 augmented metric graphs.

The following proposition is a refinement of a result from~\cite{Cha12} on vertex-weighted metric graphs 
(itself a strengthening of results from \cite{BN09}):

\begin{prop}~\label{prop:hypercharacterization} For a minimal
 augmented metric graph $\Gamma$ of genus at least two,
the following assertions are equivalent: 
\begin{enumerate}
\item $\Gamma$ is  hyperelliptic;
\item There exists  
an involution $s$ on $\Gamma$ such that:
\begin{enumerate}
\item $s$ fixes all the points $p \in  \Gamma$ with 
 $g(p) > 0$;
\item the quotient 
$ \Gamma/s$ 
is a
  metric tree; 
\end{enumerate}  
 \item There exists 
an effective finite harmonic morphism of degree two $\phi: \Gamma\rightarrow T$ from $\Gamma$ 
to a metric tree $T$ such that the local degree at 
any point $p \in \Gamma$ with $g(p) >0$ is two. 
 \end{enumerate}
Furthermore  if the involution $s$
exists, then it is unique.

\end{prop}

\pf
The implication $(2) \Rightarrow (3)$ is 
obtained by taking $T=\Gamma/s$ and letting $\phi$ be the natural quotient map.

To prove $(3) \Rightarrow (1)$, 
we observe that a finite harmonic morphism of degree two 
$\phi:  \Gamma \rightarrow T$ with local degree two at each vertex $p$
with 
$g(p) >0$
naturally extends to an effective 
finite harmonic morphism of degree two from a tropical modification 
$\Gamma'$
of $ \Gamma^\#$ to a tropical modification $T'$ 
of $T$ as follows: 
$\Gamma'$ is obtained by modifying $\Gamma^\#$ once at the midpoint of
each of its virtual cycles, and $T'$ is obtained by modifying  $T$ precisely $g(p)$
times at each point $\phi(p)$ with $g(p)>0$.  The map $\phi$ extends uniquely to an effective
finite degree two harmonic morphism $\phi': \Gamma' \rightarrow T'$,
since $\phi$ has local degree two at $p$ whenever $g(p)>0$.
By Proposition~\ref{prop:rank-finite}, the linearly equivalent 
degree two divisors $D_x(\phi')$ have rank 
one in $\Gamma'$ as $x$ varies over all points of $T'$,
which shows that $\Gamma$ is hyperelliptic.

It remains to prove $(1) \Rightarrow (2)$.  
A \emph{bridge edge} of $\Gamma$
is an edge
$e$ such that $\Gamma\setminus e$ is not connected. Let  
$\Gamma'$ be  the augmented metric graph obtained by
removing all bridge edges from $\Gamma$.
Since $\Gamma$ is
minimal, any connected component of $\Gamma'$  has  positive
genus. In particular  the involution $s$, if exists, has to restrict
to an involution on each such connected component. This implies that
$s$ has to fix pointwise any bridge edge. Hence we may now assume
without loss of generality that
$\Gamma$ has no bridge edge.
In this case 
$s$ has the following simple definition:  for any point $p \in \Gamma$, 
since $r_{\Gamma^\#}(D) = 1$ and $\Gamma$ is 
two-edge connected, there exists a unique point $q=s(p)$ such that $D
\sim (p ) + ( q )$.  
This also proves the uniqueness of the involution.
\qed

\medskip
From now until the end of the section we assume that $\mathrm{char}(k)\neq 2$. An {\it involution} on a metrized complex $\cC$ is a finite harmonic morphism 
$s:\cC \rightarrow \cC$ with $s^2 =\mathrm{id}_{\cC}$. 
An involution is called {\it tame} if the action of the group
generated by $\angles s\cong \Z /2\Z$ on $\cC$ is tame.

If $X/K$ is a (smooth proper) hyperelliptic curve, 
then the augmented metric graph $\Gamma$ associated to stable model of $X$ is hyperelliptic.
Indeed if $s_X$ is an involution on $X$, then the quotient map 
$X\to X/s$ tropicalizes to an effective tropical morphism $\phi:\Gamma\to T$
of degree $2$. The condition that $\phi$ has local degree $2$ at each
point $p$ with  $g(p)>0$ comes from the fact that any non-constant
algebraic map from a positive genus curve to $\P^1$ has degree at
least two.
The next theorem, combined with Proposition~\ref{prop:hypercharacterization}, 
provides a complete characterization of hyperelliptic 
augmented metric graphs which can be realized as the skeleton of a 
hyperelliptic curve over $K$.

\medskip

\begin{thm}\label{thm:characterization hyperelliptic}
Let $\Gamma$ be a minimal hyperelliptic augmented metric graph, and let
$s:\Gamma\to\Gamma$ be
the involution given by Proposition~\ref{prop:hypercharacterization} $(2)$.  
Then the following assertions are equivalent:
\begin{enumerate}
\item There exists a hyperelliptic smooth proper curve $X$ over $K$ and
an involution $s_X: X \rightarrow X$ such that $\Gamma $ is the
minimal skelton of $X$, and $s$ coincides with the reduction of  $s_X$ to $\Gamma$.

\item For every $p \in \Gamma$ we have
$$2g( p ) \ge \kappa( p )-2,$$
where $\kappa(p)$ denotes the number of tangent directions at $p$ which are
fixed by $s$. 
\end{enumerate}
\end{thm}

\pf
Consider the finite harmonic morphism 
$\pi: \Gamma \to \Gamma /s$. 
We note that the tangent directions at $p$ which are fixed by $s$ are exactly those along which $\pi$ has local degree two. 
Thus the condition $2g(p) \geq \kappa(p) -2$ is equivalent to the ramification index $R_p$ 
being non-negative: see Section~\ref{sec:prelim}. 
This proves $(1) \Rightarrow (2)$. 

To prove $(2) \Rightarrow (1)$, we use
Proposition~\ref{prop:lifting} and Theorem~\ref{thm:liftinggroupaction}.
According to these results, it suffices to prove that the involution $s: \Gamma \rightarrow \Gamma$ lifts to an 
involution $\overline s: \cC \rightarrow \cC$ for some metrized
complex $\cC$
with underlying augmented metric graph $\Gamma$ 
such that $\cC /\overline s$ has genus zero. 
The existence of such a lift follows from the observation that Hurwitz numbers of
degree two are all positive (see~\parref{par:liftingMGMC}, Remark~\ref{rem:Hurwitznumbers}).
\qed

\begin{figure}[h]
\begin{tabular}{c}
\scalebox{.32}{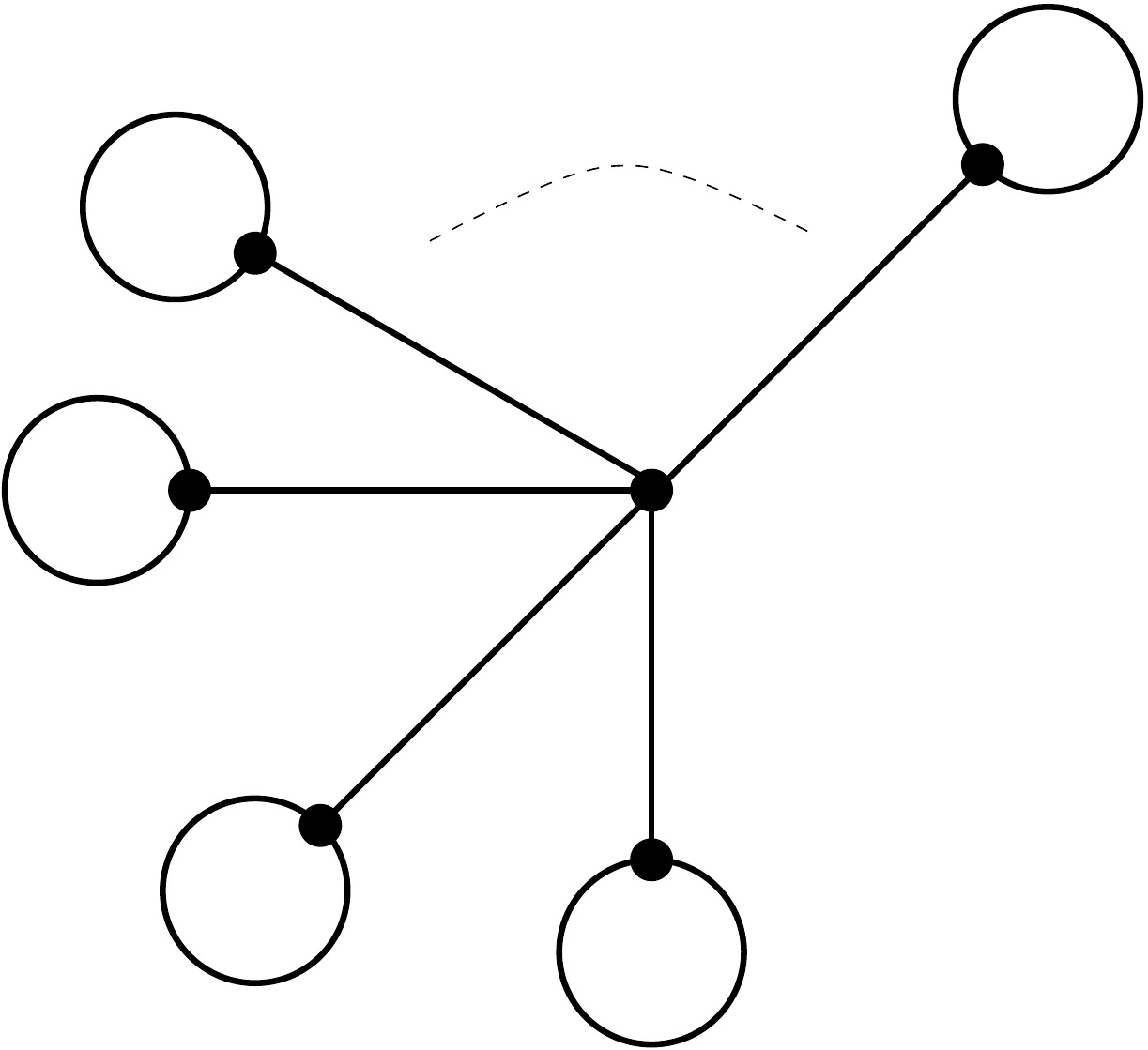}
\end{tabular}
\caption{}\label{hyper pic}
\end{figure}

\begin{eg}\label{ctrex hyper}
Let $\Gamma$ be the augmented 
metric graph of genus $g$ depicted in Figure \ref{hyper pic}
with arbitrary positive lengths. 
It is clearly hyperelliptic, and
since the involution $s$ restricts to the identity on each bridge
edge, it fixes all tangent directions at $p$. 
Then one can lift $\Gamma$ as a hyperellitptic
curve of genus $g$ if and only $2g(p)\ge \kappa-2$. In particular, if
$g(p)=0$ then this metric graph cannot be realized as the skeleton of a
 hyperelliptic curve 
as soon as $\kappa\ge 3$.
\end{eg}

\medskip

Since the hyperelliptic involution is unique for both curves and minimal augmented metric graphs, and since the tangent directions fixed by the hyperelliptic involution on an augmented metric graph correspond to bridge edges, 
we can reformulate Theorem~\ref{thm:characterization hyperelliptic} as follows,
obtaining a metric strengthening of \cite[Theorem 4.8]{caporaso:gonality}:

\begin{cor}
\label{cor:characterization hyperelliptic}
Let $\Gamma$ be a minimal augmented metric graph of genus $g \geq 2$.
Then there is a smooth proper hyperelliptic curve $X$ over $K$ of genus $g$ having $\Gamma$ as its minimal skeleton 
if and only if $\Gamma$ is hyperelliptic and for every $p \in \Gamma$ the number of bridge edges adjacent to $p$
is at most $2g(p)+2$.
\end{cor}

\section{Gonality and rank}\label{sec:examples}

A fundamental (if vaguely formulated) question in tropical geometry is the following:
If $X$ is an algebraic variety and $\T X$ is a tropicalization of $X$
(whatever it means), which properties of $X$ can be read off from $\T X$?
In this section, we discuss more precisely (for curves) the relation between
the classical and tropical notions of gonality, and of the rank of a divisor. 
It is not difficult to prove that the gonality of a tropical curve (resp.\
the rank of a tropical divisor) provides a lower bound
for the gonality (resp.\ an upper bound for the rank) of any lift (this
is a consequence, for example, of Corollary~\ref{cor:morphism.to.harmonic}).
Here we address the question of
sharpness for these inequalities: 

\begin{enumerate}
\item Can a $d$-gonal (augmented or non-augmented) tropical curve $C$ always be lifted to a $d$-gonal
  algebraic curve?
\item Can a divisor $D$ on an (augmented or non-augmented) tropical curve $C$ always be lifted to divisor of the same rank on an algebraic curve lifting $C$?
\end{enumerate}

It follows immediately  from Theorem~\ref{thm:lifting harm} that the
answer to Question $(1)$ is \emph{yes} if $C$ is not
augmented, i.e., if we are allowed to arbitrarily increase 
the genus of finitely many points in $C$. On the other hand, 
we prove in this section that the answer to Question $(1)$ in the case $C$ is
augmented, and the answer to Question 
$(2)$ in both cases, is \emph{no}.
 
\medskip

We refer to \cite{BakerNorine, MikhalkinZharkov, AC11, AB12} for the basic definitions concerning 
ranks of divisors on metric graphs, augmented metric graphs, and metrized complexes of curves.

\paragraph[Gonality of augmented graphs versus gonality of algebraic curves]
\label{par:gonality}
 An augmented tropical curve $C$ is said to have  an augmented
 (non-metric) graph $G$ as its {\it combinatorial type} if $C$ admits
 a representative 
whose underlying augmented graph is $G$.
Given an augmented graph $G$, 
we denote by $\mathcal M(G)$ the set of all augmented metric graphs which define a tropical curve $C$
 with combinatorial type $G$.  
 In other words, $\mathcal M(G)$ consists of all augmented metric graphs which can be obtained
by a finite sequence of tropical modifications (and their inverses) from an augmented metric graph $\Gamma$ 
with underlying augmented graph $G$. 
When  no confusion is possible, we identify an (augmented)
tropical curve with any of its representatives as an (augmented) metric graph: in what follows, 
we deliberately write $C \in \mathcal M(G)$ for a tropical curve $C$ with combinatorial type $G$. 
Note that the spaces $\mathcal M(G)$ appear naturally in the
 stratification of the moduli space of tropical curves of 
 genus $g(G)$, see for example \cite{caporaso:gonality}.

\begin{defn}
An augmented
tropical curve $C$ is called \emph{$d$-gonal} if there exists a
tropical morphism  $C\to\T\P^1$ of degree $d$.

An augmented graph $G$ is called \emph{stably $d$-gonal} if there
exists a $d$-gonal augmented tropical curve $C$ whose combinatorial type is $G$.
\end{defn}

In other words, an augmented graph $G$ is stably $d$-gonal if and only if there is an augmented metric graph $\Gamma \in \mathcal M(G)$ 
which admits an effective finite harmonic morphism of degree $d$ to a metric tree. 

 \begin{rem} \label{rem:gonality.definitions}
Our definition of the stable gonality of a graph is equivalent to the one given in \cite{cornelissen:graph_li-yau}.
See Appendix A of {\em loc. cit.} for a detailed discussion of the relationship between stable gonality and other tropical or graph-theoretic notions of
gonality in the literature, e.g. Caporaso's definition in \cite{caporaso:gonality}.  
 \end{rem}

\medskip

In this section we prove the following theorem, which is an immediate
consequence of Corollary \ref{cor:morphism.to.harmonic} and  
Propositions~\ref{non42} and~\ref{non4} below.

\begin{thm}\label{thm:nonlift.gonality}
There exists an augmented stably $d$-gonal graph $G$ such that for any 
augmented metric graph  $\Gamma \in \mathcal M(G)$ and any smooth proper connected $K$-curve $X$
lifting $\Gamma$, the gonality of $X$ is strictly larger than $d$.
\end{thm}

\medskip

\begin{figure}[h]
\scalebox{.32}{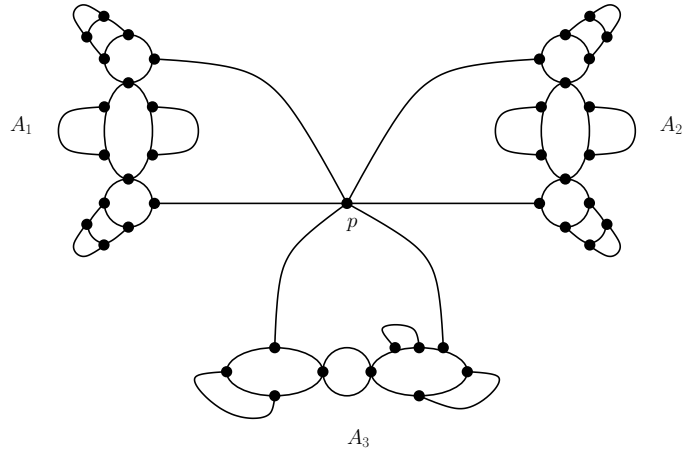}
\caption{The graph $G_{27}$}\label{23}
\end{figure}

Let $G_{27}$  be the graph depicted in Figure~\ref{23}, which we promote to a totally degenerate 
augmented graph by taking the genus function to be identically zero. Note that $g(G_{27})=27$, and that
  $G_{27} \setminus\{p\}$ has three connected components, which we denote by
  $A_1$, $A_2$, and $A_3$ according to Figure~\ref{23}.
  
Given 
an element $\Gamma \in \mathcal M(G_{27})$ and a tropical morphism
$\phi:C\to \T\P^1$ 
 from the tropical curve represented by $\Gamma$ to $\T\P^1$,  we denote by $\phi_i$ the
restriction of $\phi$ to 
(the metric subgraph in $\Gamma$ which corresponds to)
$A_i$, and by $\phi_p$ the restriction of $\phi$ to a
small neighborhood of the point $p$. 

\begin{prop}\label{non42}
The graph $G_{27}$ depicted in Figure~\ref{23} is stably $4$-gonal. 
\end{prop}

\pf
\begin{figure}[h]
\begin{tabular}{ccccc}
\scalebox{.21}{\input{Figures/Mapv.pdf_t}}& \hspace{3ex} &
\includegraphics[width=4cm, angle=0]{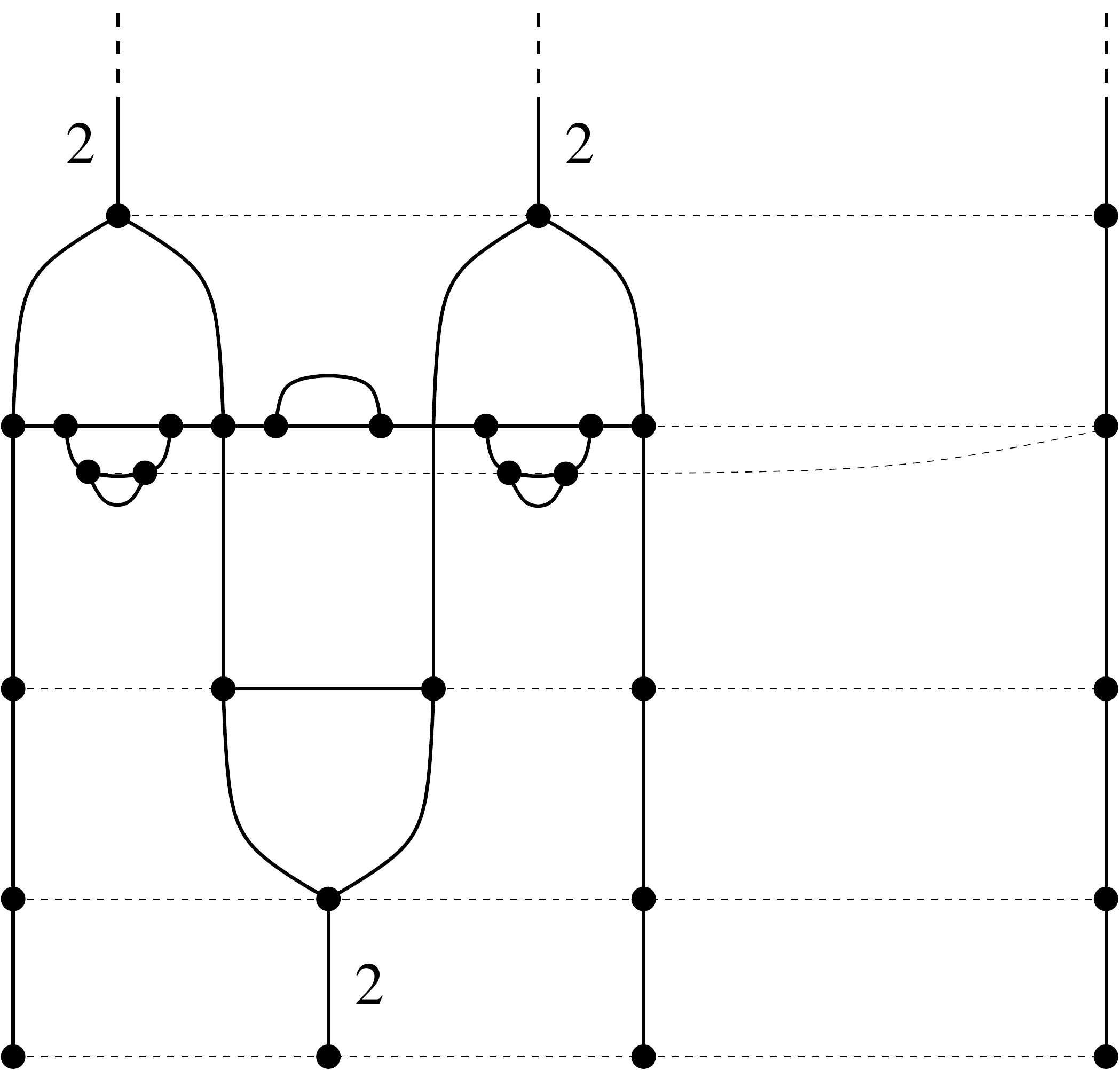}& \hspace{3ex} &
\includegraphics[width=4cm, angle=0]{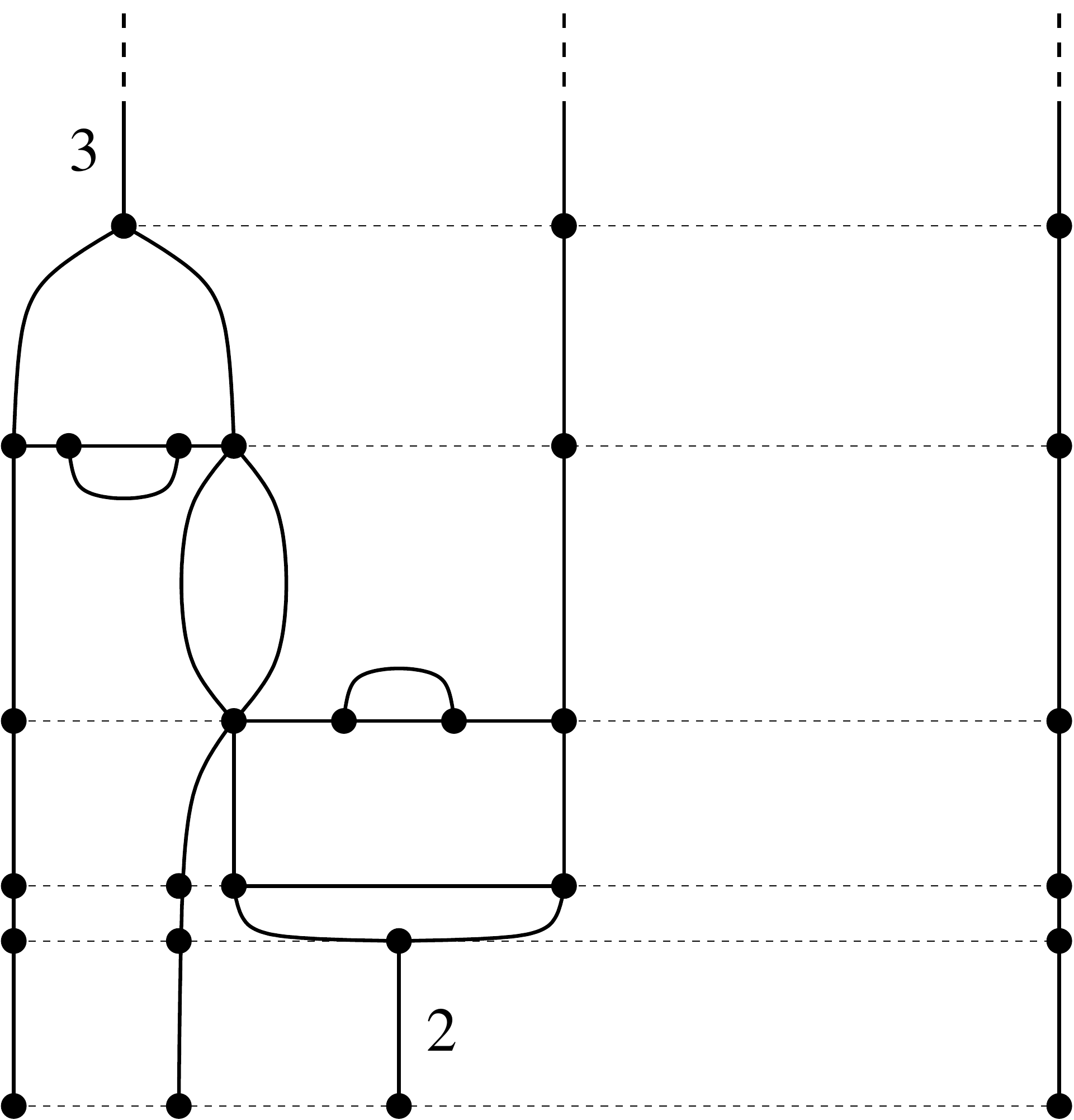}
\\ \\ a) $\phi_v$ && b) $\phi_1=\phi_2$&& c) $\phi_3$
\end{tabular}
\caption{A tropical morphism of degree four.}\label{Map}
\end{figure}
We need to show the existence of a suitable tropical curve $C$ with combinatorial type $G_{27}$ which admits a 
 tropical morphism of degree four to $\T\P^1$. 
For a suitable choice of edge lengths on $G_{27}$, we get an element
$\Gamma \in \mathcal M(G_{27})$ such that there exists a harmonic morphism from $\Gamma$
to a
star-shaped genus zero augmented metric graph with three infinite
edges, which has 
restrictions $\phi_1, \phi_2, \phi_3, \phi_v$ to $A_1, A_2, A_3$, and
a small neighborhood of $p$, respectively, given as in Figure
\ref{Map}.
We claim that $\phi$ induces a tropical morphism,
i.e., that there exists a tropical modification of $\phi$ which is finite and effective. 

\begin{figure}[h]
\begin{tabular}{ccc}
\includegraphics[width=3cm, angle=0]{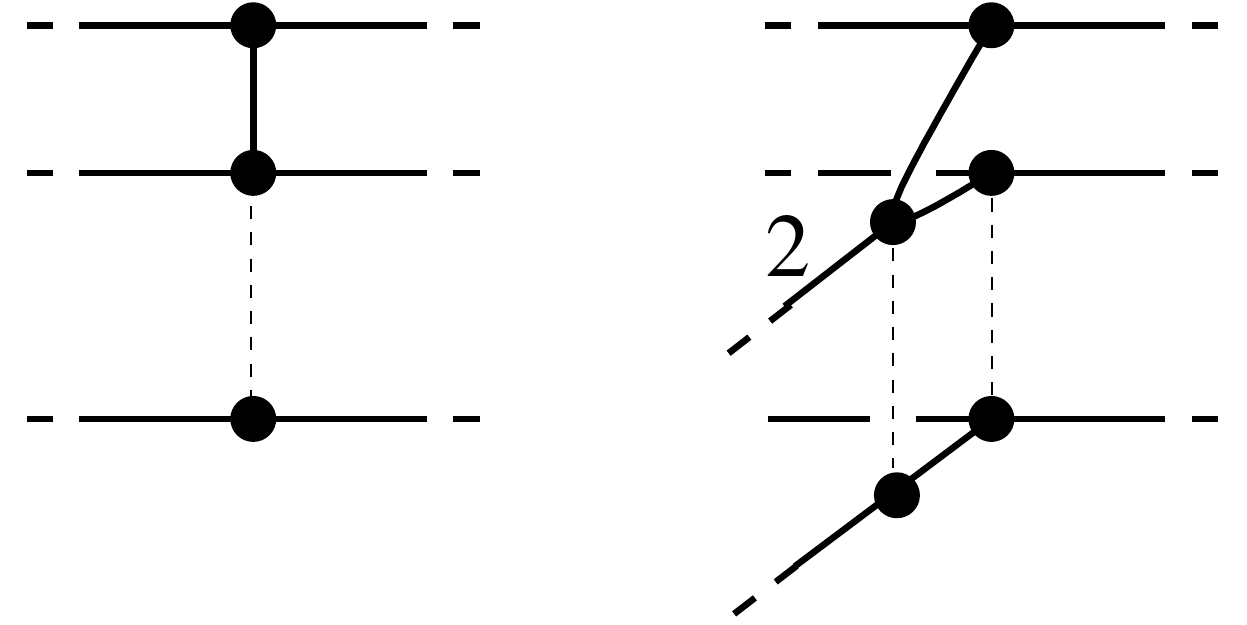}& \hspace{.3cm} 
\includegraphics[width=3.5cm, angle=0]{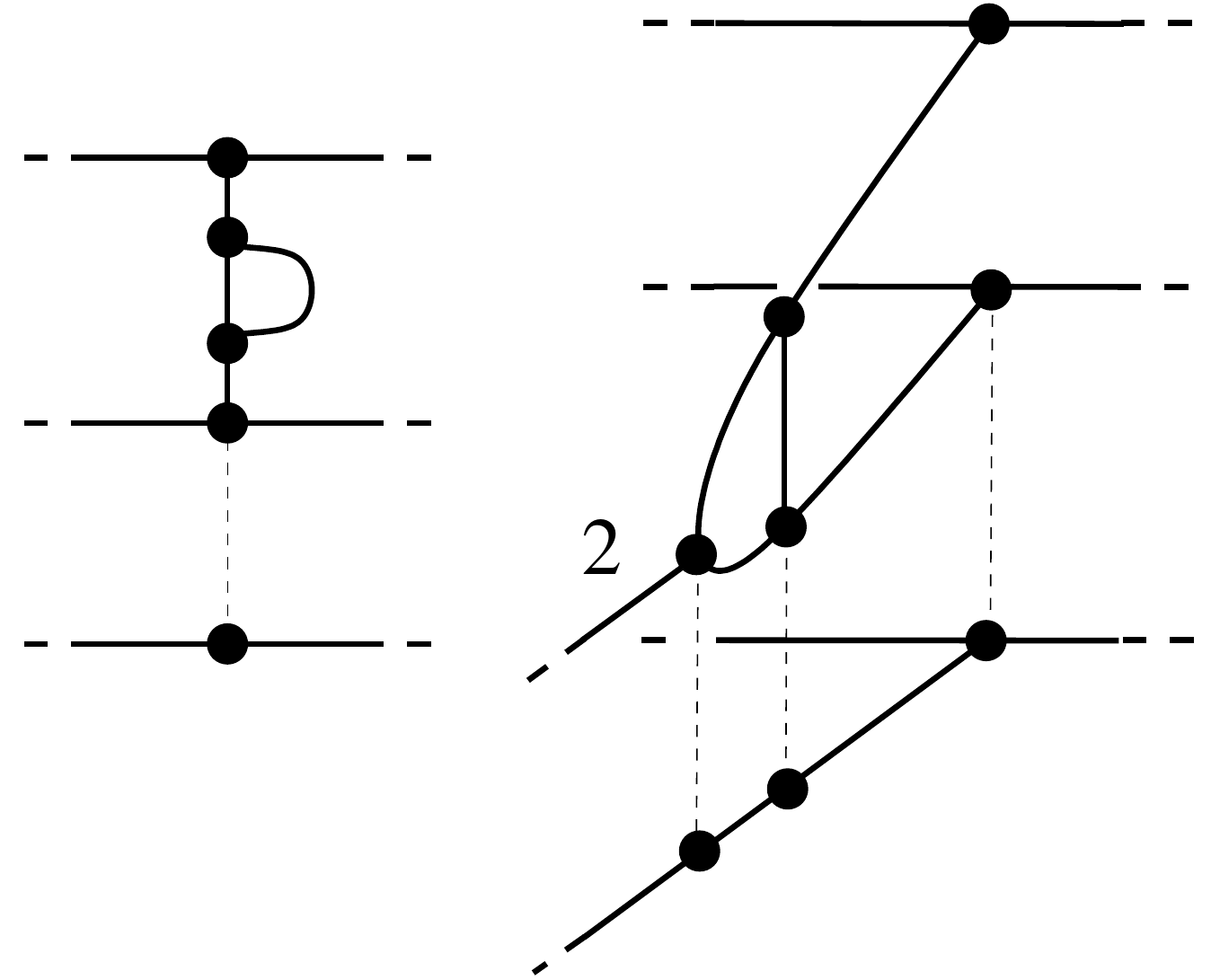}& \hspace{.3cm}
\includegraphics[width=3.5cm, angle=0]{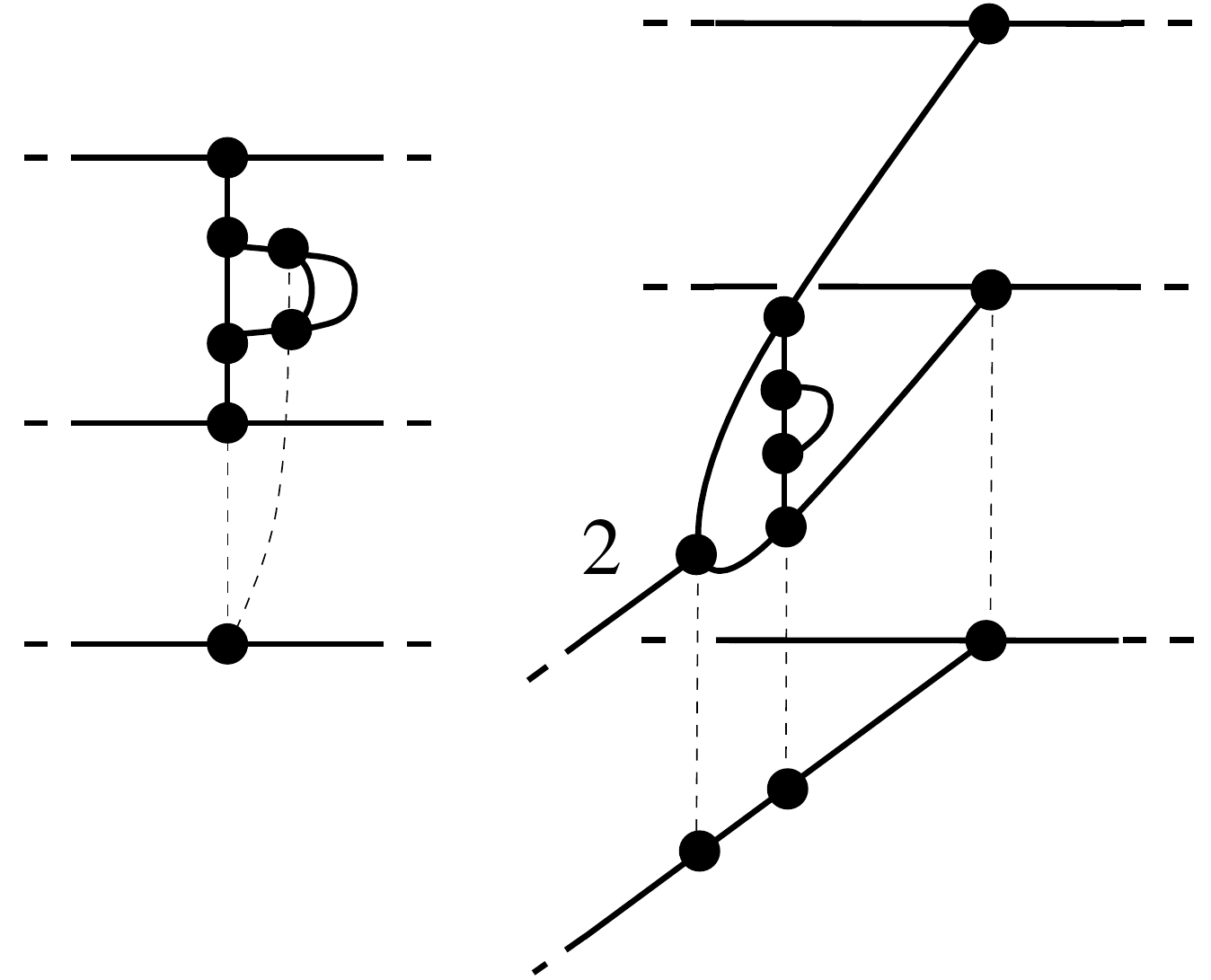}
\\ \\ a) Resolution in one step & \hspace{.3cm}  b) Resolution in two steps  & \hspace{.3cm} c) Resolution in three steps
\\ & \,\,\,\,\,\,\,\,\, (combined with
case a) &  \,\,\,\,\,\,\, (combined with
case b)

\end{tabular}
\caption{Patterns to resolve contractions in the harmonic morphisms $\phi_1,\phi_2,$ and $\phi_3$.}\label{resolve contr}
\end{figure}

Note that each of the morphisms $\phi_1$ and $\phi_2$ contains a fiber of genus five, 
while the morphism $\phi_3$ has two different fibers of genus one. 
All the other fibers of $\phi_1, \phi_2,$ and $\phi_3$ are either finite or connected of genus zero.   
We depict in  Figure \ref{resolve contr} a few patterns which show how to resolve
contractions of $\phi$, turning $\phi$ into an augmented tropical morphism.
Figure~\ref{resolve contr}(a) shows  how to resolve a contracted segment 
(resolving contracted fibers of genus zero). 
Figure~\ref{resolve contr}(b) shows how to resolve a contracted cycle 
(resolving the contracted cycles in $\phi_3$ and the middle contracted cycle in $\phi_1$ and $\phi_2$): 
the idea is to reduce to the case of a contracted segment, 
in which case one can use the resolution given in Figure~\ref{resolve contr}(a) to finish. 
And finally, Figure~\ref{resolve contr}(c) shows how to resolve the two contracted double-cycles in $\phi_1$ and $\phi_2$ by reducing to the case
already treated in Figure~\ref{resolve contr}(b).
 Note that performing these tropical modifications  
 impose conditions on the length of the contracted edges in $\Gamma$,
 e.g., in Figure~\ref{resolve contr}(b), the two edges adjacent to the contracted  
cycle should have the same length. 
Nevertheless, by appropriately choosing the edge lengths, 
we get the existence of a metric graph $\Gamma \in \mathcal M(G_{27})$  which admits a finite morphism of 
degree four to a metric tree. It is easily seen that this morphism is effective; thus we get a tropical 
curve $C$ with combinatorial type $G_{27}$ and a tropical morphism of degree four to $\T \P^1$, 
finishing the proof of the proposition.
\qed

\medskip

To conclude the proof of Theorem~\ref{thm:nonlift.gonality}, we now show the following:

\begin{prop}\label{non4}
There is no  metrized complex of $k$-curves 
with underlying augmented metric graph in $\mathcal M(G_{27})$ and admitting a finite morphism 
of degree four to a metrized complex of $k$-curves of genus zero.
\end{prop}

\medskip

We emphasize that the statement holds for any (algebraically closed) field $k$.
The proof of Proposition \ref{non4} relies on some technical lemmas that we are now going to state.

\medskip

We first recall a formula given in~\cite{AB12} for the rank of divisors on a metric graph $\Gamma = \Gamma_1\vee \Gamma_2$ which is 
obtained as a wedge sum of two metric graphs $\Gamma_1$ and $\Gamma_2$.  Recall that given two metric graphs $\Gamma_1$ 
and $\Gamma_2$ and distinguished points $t_1\in \Gamma_1$ and $t_2\in \Gamma_2$, the {\it wedge sum} or {\it direct sum} of $(\Gamma_i,t_i)$, 
denoted ${\Gamma} = \Gamma_1 \vee \Gamma_2$, is the metric graph obtained by identifying the points $t_1$ and $t_2$ in the disjoint union of 
$\Gamma_1$ and $\Gamma_2$. Denoting by $t \in {\Gamma}$ the image of $t_1$ and $t_2$ in $\Gamma$, one refers to $t\in \Gamma$ 
as a {\it cut-vertex} and to $\Gamma=\Gamma_1\vee \Gamma_2$ as the decomposition corresponding  to the cut-vertex $t$.  
(By abuse of notation, we will use $t$ to denote both $t_1$ in $\Gamma_1$ and $t_2$ in $\Gamma_2$.) 
There is an addition map $\Div(\Gamma_1) \oplus \Div(\Gamma_2) \rightarrow \Div(\Gamma)$ which sends a pair of divisors $D_1$ and $D_2$ in $\Div(\Gamma_1)$ and $\Div(\Gamma_2)$ to the divisor
$D_1+D_2$  on $ \Gamma$ defined by pointwise addition of the coefficients in $D_1$ and $D_2$.

\medskip

Let $D_1\in \Div(\Gamma_1)$ and $D_2 \in \Div(\Gamma_2)$. For any
non-negative $m$, define 
$\eta_{\Gamma_1,D_1}( m )$ as minimum integer $h$ such that
$r_{\Gamma_1}( D_1 + h(t_1)) = m$. 
 Then 
\begin{equation}\label{eq:ranksum}
r_{\Gamma}(D) = \min_{m\geq 0} \big\{m + r_{\Gamma_2}(D_2-\eta_{\Gamma_1,D_1}(m)(t_2))\big\}.
\end{equation}
(see~\cite{AB12} for details). 

\medskip

\begin{figure}[h]
\begin{tabular}{ccc}
\scalebox{.38}{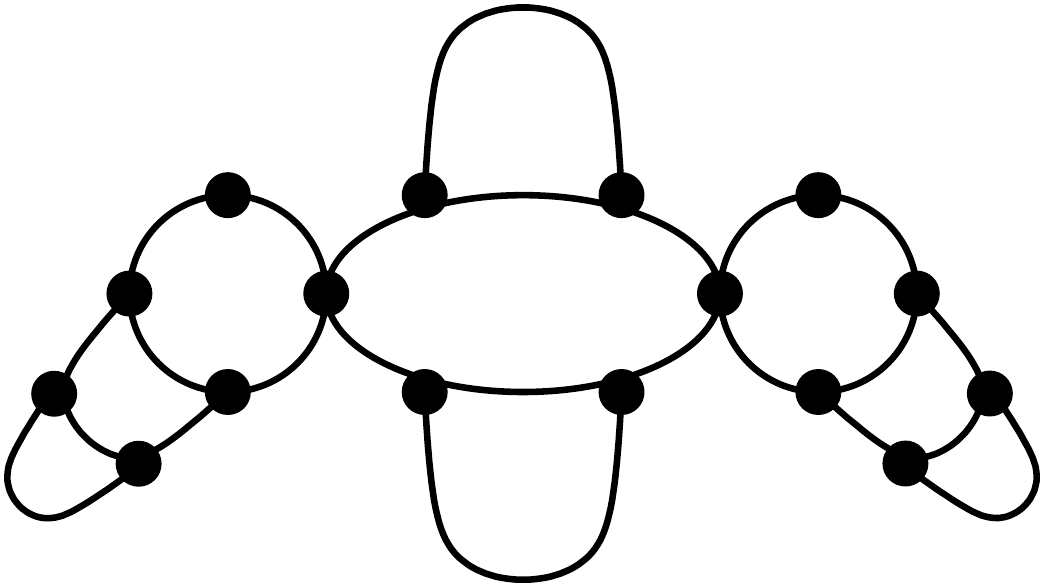}& \hspace{4ex} &
\scalebox{.38}{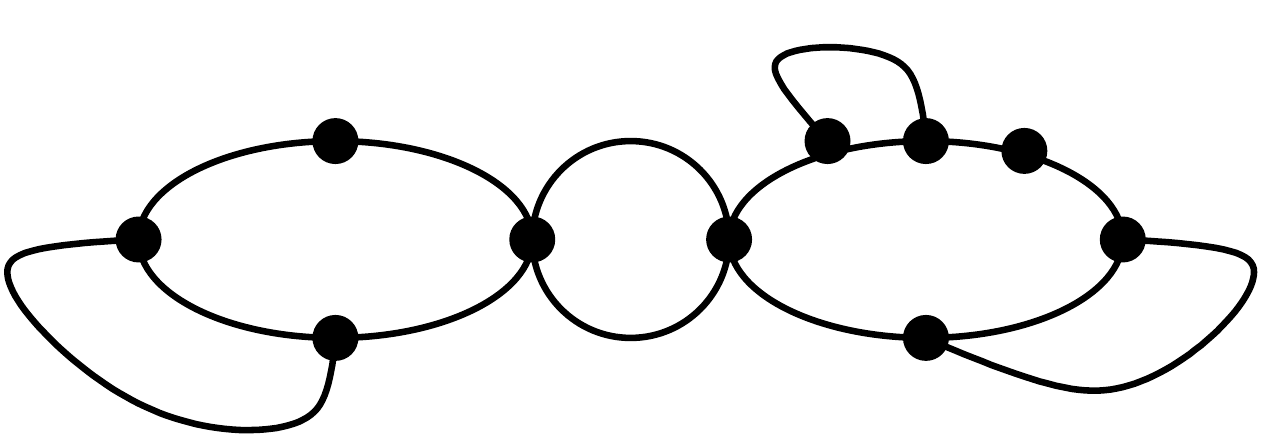}
\\ \\ a) A metric graph $\Gamma$ in $\mathcal M(A_1)=\mathcal M(A_2)$ && b) A metric graph $\Gamma$ in $\mathcal M(A_3)$
\end{tabular}
\caption{}\label{g9}
\end{figure}

In what follows, equation \eqref{eq:ranksum}
will be applied to a metric graph $\Gamma \in \mathcal M(A_1) =
\mathcal M(A_2)$
(see~Figure~\ref{g9}(a) 
and Lemma~\ref{lg9}), 
 to a metric graph $\Gamma \in \mathcal M(A_3)$ 
(see~Figure~\ref{g9}(b) 
and Lemma~\ref{lg6}), and to $\Gamma_{27} \in \mathcal M(G_{27})$ with cut-vertex $p$ 
in the proof of Proposition~\ref{non4}.

\begin{lem}\label{lg9}
Let $\Gamma$ be a metric graph in $\mathcal M(A_1) = \mathcal M(A_2)$ as
depicted in Figure \ref{g9}(a).  For any non-negative integers $a\le 3$ and $b\le 1$, the divisors $a(
p )+b( q )$ and $b( p )+a( q )$ have rank zero in $\Gamma$. 
\end{lem}

\pf
By symmetry it is  enough to prove the lemma for the divisor $D=3( p )+( q )$. 
Consider the decomposition $\Gamma = \Gamma_{p} \vee \Gamma_{q}$ associated to the cut-vertex $t$ in $\Gamma$, where $\Gamma_{p}$ and $\Gamma_{q}$ denote the closure in $\Gamma$ of the 
the two connected components of
$\Gamma\setminus\{t\}$ which contain the points $p$ and $q$, respectively.

We claim that 
$\eta_{\Gamma_q,(q)}(1) = 3$.
 Assume for the moment that this is true. Then by~\eqref{eq:ranksum}, we have 
$$0 \leq r_\Gamma (3 (p ) + ( q))  \leq 1 + r_{\Gamma_{p}} (3 ( p ) - 3 ( t )).$$ 

By Lemma~\ref{lem:g3} below, in $\Gamma_{p}$ we have $r_{\Gamma_{p}} (3 ( p ) - 3 ( t )) = -1$.  We thus infer that  $r_\Gamma (3 (p ) + ( q))=0$. 

\medskip

It remains to prove that 
$\eta_{\Gamma_q,(q)}(1) = 3$.
In other words, we need to show that in $\Gamma_{q}$ we have $r_{\Gamma_{q}}(2 ( t ) + ( q )) = 0$. 
For this, consider the decomposition $\Gamma_{q} = \Gamma_q^t \vee \Gamma_q^q$ 
corresponding to the cut-vertex $s$ in $\Gamma_{q}$, where
$\Gamma_q^t$ and $\Gamma_q^q$ denote the components which contain $t$
and $q$, respectively.
We claim that 
 $\eta_{\Gamma_{q}^t,2(t)} (1) =1$.
Assuming the claim, we 
 have $0 \leq r_{\Gamma_{q}}(2 ( t ) + ( q )) \leq 1 + r_{\Gamma_q^q} ( ( q ) - ( s )) = 0$ 
(since $q$ and $s$ are not linearly equivalent in $\Gamma_q^q$; see Lemma~\ref{lem:g3}). 
So it remains to prove that 
$\eta_{\Gamma_{q}^t,2(t)} (1) =1$.
This is equivalent to $r_{\Gamma_q^t} ( 2 ( t ))=0$, which is obviously the case. 
\qed

\begin{figure}[h!]
\begin{tabular}{ccc}
\scalebox{.36}{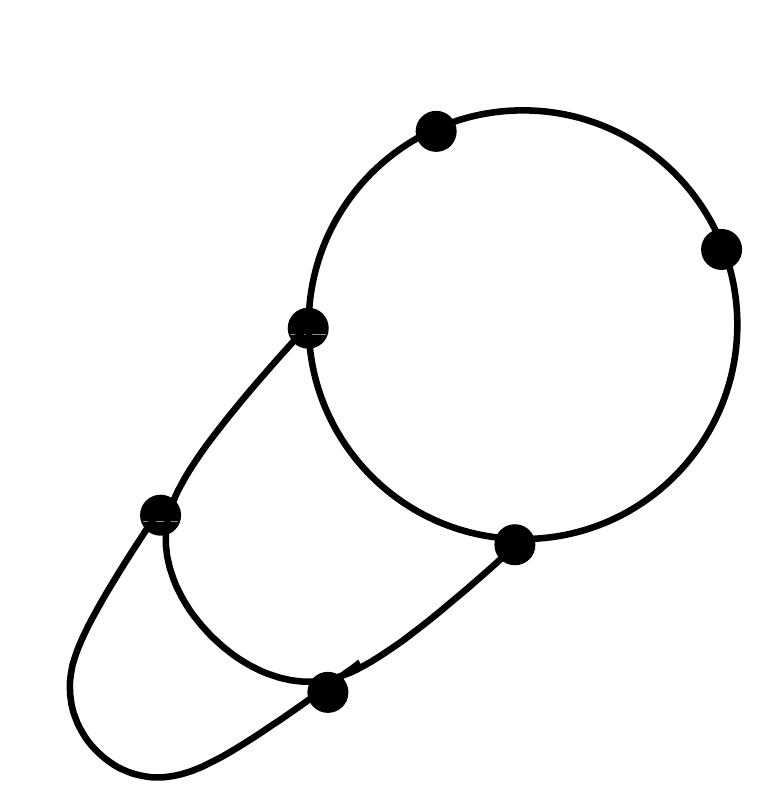}& \hspace{2.5cm} &
\scalebox{.36}{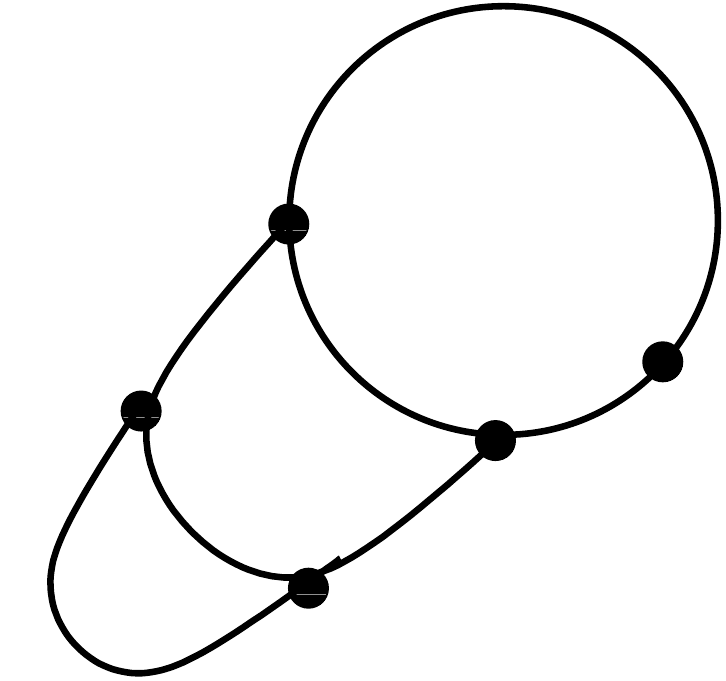}
\\ \\ a) 
the divisor $3( p ) - 3 ( t )$  is not && b) $3( p )- 3 ( t ) \sim 3 ( u ) - 3( t' )$
\\ 
\hspace{0.5cm} 
rationally equivalent to zero. && 
\end{tabular}
\caption{}\label{fig:g3}
\end{figure}

\begin{lem}\label{lem:g3}
Let $\Gamma$ be any metric graph in $\mathcal M(G_3)$, where $G_3$ is
the totally degenerate
graph depicted in Figure~\ref{fig:g3}(a).
Then the two divisors $3( p )$ and $3( t )$ are not linearly equivalent in $\Gamma$.
\end{lem}

\pf
By symmetry we can assume that the length of the edge $\{u,p\}$ is less than or equal to the length of 
the edge $\{t, w\}$.  
Then there exists a point $t'$ in the 
segment $[t,w]$ so that $3( p )- 3 ( t ) \sim 3 ( u ) - 3( t' )$ ---
see Figure~\ref{fig:g3}(b) ---
and we are led to prove that $D = 3 ( u ) - 3( t' )$ is not linearly equivalent to zero. 
Consider the unique $t'$-reduced divisor $D_{t'}$ linearly equivalent to $D$ in $\Gamma$
(see e.g.~\cite{A09, BakerNorine} for the definition and basic properties of reduced divisors). 
It will be enough to show that $D_{t'} \neq 0$. Three cases can occur, depending on the lengths 
$\ell_z, \ell_w$, and $\ell_{t'}$  in  $\Gamma$ of the edges $\{u,z\}, \{u,w\}$, and the segment $\{u,t'\}$, 
respectively:
\begin{itemize} 
\item If $\ell_z  = \min\Bigl\{ \ell_z, \ell_u, \ell_{t'}\Bigr\}$, 
then there are two points $w'$ and $t''$ on the segments $\{u,w\}$ and $\{u,t'\}$, 
respectively, such that $D_{t'} = ( z ) + (w') + (t'') - 3 ( t ')$. 
\item If $\ell_u  = \min\Bigl\{ \ell_z, \ell_u, \ell_{t'}\Bigr\}$, 
then there are two points $z'$ and $t''$ on the segments $\{u,z\}$ and $\{u,t'\}$, respectively, 
such that $D_{t'} = ( z' ) + (w) + (t'') - 3 ( t ')$. 

\item If $\ell_{t'}  = \min\Bigl\{ \ell_z, \ell_u, \ell_{t'}\Bigr\}$, then there are two
 points $z'$ and $w'$ on the segments $\{u,z\}$ and $\{u,w\}$, respectively, such that $D_{t'} = ( z' ) + (w') - 2 ( t ')$. 
\end{itemize}
In all the three cases, we have $D_{t'} \neq 0$, which shows that $D$ cannot be equivalent to zero in $\Gamma$. 
\qed

\begin{lem}\label{lg6}
Let $\Gamma \in \mathcal M(A_3)$ be a metric graph as
depicted in Figure~\ref{g9}(b).
 For any $a,b\le 2$, the divisor $a ( p )+b ( q )$ has rank
zero on $\Gamma$.
\end{lem}

\pf The arguments are  similar to the ones used in the proof of Lemma~\ref{lg9}. 
Consider the cut-vertex $t$ in $\Gamma$ and denote by $\Gamma_{p}$ and $\Gamma_q$ the corresponding components containing $p$ and $q$, respectively.   
We claim that 
$\eta_{\Gamma_q,2(q)}(1) = 2$. 
This obviously implies the lemma. Indeed, $r_{\Gamma_p} (2 ( p ) - 2 ( t ) )=-1 $ (which can be verified by an analogue of Lemma~\ref{lem:g3} in $\Gamma_p$), and thus~\eqref{eq:ranksum} implies that 
$r_\Gamma (2 ( p )+ 2( q)) \leq 1 +  r_{\Gamma_p} (2 ( p ) - 2 ( t ) ) =0$. 

To show that 
$\eta_{\Gamma_q,2(q)}(1) = 2$,
it will be enough to show that $r_{\Gamma_q}(2( q )+(t)) =0$. 
This can be done in exactly the same way by considering the other cut-vertex $s$ adjacent to $t$ in $\Gamma_q$.
\qed

\begin{lem}\label{lem:non deg 4}
Let $x_1,x_2,$ and
$x_3$ be distinct points in $\P^1(k)$. Then there does not exist a morphism $f:\P^1\to \P^1$ 
of degree four branched over $x_1$, $x_2$ and $x_3$ and having ramification profile $(2,2)$, $(2,2)$, and $(3,1)$ at these three points.
\end{lem}

\pf
Suppose that such a rational map $f:\P^1\to \P^1$
exists. The monodromy group of $f$ is a subgroup of $\frak S_4$, so
its cardinality is of the form $2^a3^b$. In particular, if the
characteristic of $k$ is neither $2$ nor $3$, then $f$ has a tame
monodromy group and the non-existence of $f$ then  comes from the fact
that
$H_{0,0}^4\left((2,2),(2,2),(3,1) \right)=0$
(see~Example~\ref{rem:Hurwitznumbers}). 

\smallskip

Hence it remains to check the lemma for $\mathrm{char}(k)=2,3$. Note
that the same technique we use in this case works in any characteristic,
but the computations are a bit more tedious in characteristic different from $2$ and $3$.

\smallskip

Up to the action of ${\rm GL}(2,k)$ on $\PP^1$ via automorphisms, we may
assume that $x_1=0$, $x_2=\infty$, and $x_3=1$, and that
$$f(X)=a\frac{X^2(X+1)^2}{(X+b)^2} $$
with $a\ne 0$ and $b\ne 0, -1$.
Hence the condition on the ramification profile of $x_3$ translates as
$$aX^2(X+1)^2 - (X+b)^2= c(X-d)^3(X-e) $$
with $c\ne 0$, $d\ne 0, -1, b$, and $e\ne 0, -1 ,b,
d$. Looking at the coefficients of the two polynomials, we obtain the  following
five equations
$$ (E_1): \ a=c ,\quad \quad  (E_2):\ ec=-2a-3cd, \quad  \quad
(E_3):\ a-1=3cd(d +e) 
,$$ 
$$\quad (E_4): \ 2b=cd^2(d+3e) ,\quad  \quad \quad (E_5):\ -b^2=cd^3e
.$$

If $k$ has characteristic $2$, then $(E_2)$ becomes $ec=cd$ which
contradicts the fact that $e\ne d$.

If $k$ has characteristic $3$, then these five equations become
$$ (E_1): \ a=c ,\quad  (E_2):\ ec=a, \quad  (E_3):\ a=1, 
\quad (E_4): \ -b=cd^3 ,\quad  (E_5):\ -b^2=cd^3e.$$
Equations $(E_1),(E_2),(E_3)$ imply $a=c=e=1$.  Then
$(E_4)$ and $(E_5)$ become $-b = d^3 = -b^2$; hence $b=1=e$, which contradicts
our assumptions. 
\qed

\medskip

We can now give the promised proof of Proposition \ref{non4}.

\medskip

\pf (Proof of Proposition \ref{non4})
Suppose that there  exists a metrized complex of
  $k$-curves $\cC_{27}$ of genus $27$
with underlying augmented metric graph $\Gamma_{27}$ in 
  $\mathcal M(G_{27})$, and admitting a  finite harmonic morphism of
metrized complexes of degree four
 $\phi : \cC_{27}\to \mathcal T$, for  
 $ \mathcal T$
 of genus zero with underlying metric tree denoted by $T$.  Without loss of
generality, we may assume that 
$T$ has no infinite vertex $q\in V_\infty(T)$ such that
 any infinite edge $e'$ adjacent to
an infinite vertex $q'\in\phi^{-1}(q)$ has $d_{e'}(\phi)=1$.

We are going to prove below that the local degree at $p$ is $4$.
Assuming that this is the case, we show how the proposition follows. 
Denote by $\Gamma_1,\Gamma_2,$ and $\Gamma_3$ 
the three components of $\Gamma_{27} \setminus \{p\}$ which contain $A_1, A_2,$ and $A_3$, respectively. 
Since the degree of $\phi$ at $p$ is four, we  have $\phi^{-1}(\phi(p)) = \{p\}$. 
Therefore, by the connectivity of 
$\Gamma_i$, the images of $\Gamma_i$ under $\phi$ are pairwise disjoint in $T$. 
This shows that for $x$ sufficiently 
close to $\phi(p)$ in $T$, the support of the divisor $D_x(\phi)$ lives entirely in one of the $\Gamma_i$
for $i \in \{ 1,2,3 \}$.
Choose $x_i$ sufficiently close to $\phi(p)$ such that the support of $D_{x_i}(\phi)$ 
is contained in $\Gamma_i$.
Applying Proposition~\ref{prop:rank-finite}, we see that each divisor $D_{x_i}(\phi)$ has rank one in $\Gamma_i$. Now, 
according to Lemma \ref{lg9}, the degree-four divisor $D_{x_1}(\phi)$ (resp.\ $D_{x_2}(\phi)$)
must be of the form $2( a ) + 2 ( b )$ for two points $a$ and $b$ sufficiently close to $p$ and 
lying on the two different branches of $\Gamma_1$ (resp.\ $\Gamma_2$)
adjacent to $p$.  Similarly, by Lemma~\ref{lg6}, the divisor $D_{x_3}(\phi)$ has to be of the form 
$3 ( a ) + ( b )$ for 
two points $a$ and $b$ sufficiently close to $p$ and 
lying on the two different branches of $\Gamma_3$ adjacent to $p$. This shows that the map $\phi_p$, 
the restriction of $\phi$ to a sufficiently small neighborhood of $p$ in $\Gamma_{27}$, coincides with the map 
depicted in Figure\ref{Map}(a). 
The proposition now follows from Lemma~\ref{lem:non deg 4}.

\medskip

It remains to prove that $d_{p}(\phi)=4$.
We first claim that $\phi$ maps one of the components $\Gamma_i$, for $i=1,2,3$, onto a connected component of
$T\setminus\{\phi(p)\}$. Otherwise, for the sake of contradiction, 
suppose that $\phi^{-1}(\phi(p))$ consists of $p$ and 
one point $p_i$ in each of the components
$\Gamma_i$ for $i=1,2,3$.  Then $\phi$ has local degree one at each of the points $p_i$. 
By Proposition~\ref{prop:rank-finite}, $D_{\phi(p)}(\phi) = (p)+ (p_1)+(p_2)+(p_3)$ has rank one in $\Gamma$. 
By equation~\eqref{eq:ranksum} applied to the cut-vertex $p$ in $\Gamma_{27}$, 
we infer that the divisor $(p)+(p_i)$ has rank one in the metric graph $\overline{\Gamma}_i$, 
the closure of $\Gamma_i$ in $\Gamma_{27}$. In other words, the metric graphs $\overline \Gamma_i$ are 
hyperelliptic, which is clearly not the case.
This gives a contradiction and the claim follows. 

\medskip

Summarizing, there must exist at least one $\Gamma_i$ such that $\phi$ maps $\Gamma_i$ onto one of 
the connected components of
$T\setminus\{\phi(p)\}$. Reasoning again as in the first part of the proof, it follows from 
Proposition~\ref{prop:rank-finite} and Lemmas \ref{lg9} and \ref{lg6}
that the restriction of $\phi$ to $\Gamma_i$ has degree four, which
implies that $d_{p}(\phi)=4$.
\qed

\medskip

\paragraph[Lifting divisors of given rank]\label{par:ranklifting}
First, recall that to  a
smooth proper curve $X$ over $K$ together with a semistable vertex set $V$ and a subset $D_0$ of 
$X(K)$ compatible with $V$, we can naturally associate a metrized complex
of curves $\cC = \Sigma(X, V \cup D_0)$  
with underlying augmented metric graph $\Gamma$.
As in \cite{AB12}, there are natural specialization maps on divisors,
which we denote for simplicity by the same letter $\tau_*$:
$$\tau_{*}:\Div(X)\to\Div(\cC),\quad \text{and}\quad
\tau_{*}:\Div(\cC)\to\Div(\Gamma).$$

Since this discussion is pointless in the case of rational curves, we
may assume that $X$ (equivalently, $\cC$ or the augmented metric graph $\Gamma$) has positive genus.
We will also assume that $\Gamma$ does not have any infinite vertices, i.e., that $D_0$ is empty, which
does not lead to any real loss of generality and which makes various statements easier to write down and understand.
We may also assume without loss of generality that $V$ is a strongly semistable vertex set of $X$.

According to
the Specialization Inequality~\cite{baker:specialization, AC11, AB12}), for
any divisor $D$ in $\Div(X)$ one
has  
\begin{equation}\label{specialization}
r_X(D)  \leq r_\cC(\tau_{*}(D)  ) 
\leq r_\Gamma^{\#}(\tau_{*}(D)) \leq r_\Gamma(\tau_{*}(D)),
\end{equation}
where $r_X$, $r_\cC$ and $r_\Gamma$ denote rank of divisors on
$X$, $\cC$ and (unaugmented) $\Gamma$, respectively, and $r_\Gamma^{\#}$
denotes the weighted rank in the augmented metric graph $(\Gamma, g)$ (see~\parref{par:hyper}).  

\medskip 
We spend the rest of this section discussing the sharpness of the inequalities appearing in (\ref{specialization}).

\begin{defn}
Let $\cC$ be a metrized complex of curves whose underlying metric graph $\Gamma$ has no infinite leaves, 
and let $\mathcal D$ be a $\Lambda$-rational divisor in $\Div_\Lambda(\cC)$.
A \emph{lifting} of the pair $(\cC, \mathcal D)$ consists
of a triple $(X,V; D_X)$ where $X$ is a smooth proper curve over $K$,
$V$ is a strongly semistable vertex set for which $\cC = \Sigma(X, V)$, and
$D_X$ is a divisor in $\Div(X)$ with $\mathcal D \sim \tau_*(D_X)$. We say that the inequality 
$r_X  \leq r_\cC $ is \emph{sharp} if 
for any metrized complex of curves $\cC$ and any divisor
$\mathcal D\in\Div(\cC)$, there exists a lifting $(X,V;D_X)$ of
$(\cC, \mathcal D)$ such that 
$r_X(D_X) = r_\cC(\mathcal D).$

We can define in a similar way what it means to lift a divisor on an
(augmented) metric graph to a divisor on a metrized complex of curves
or to a smooth proper curve over $K$, and what it means for the
corresponding specialization inequalities to be sharp.
\end{defn}

It is easy to see that the inequality $r_\Gamma^\#\leq r_\Gamma$ is not sharp
(see~\cite{AB12} for a precise formula relating the two rank functions). 

The following example is due to Ye Luo (unpublished); we thank him for his permission to
include it here.
Together with Corollary~\ref{cor:morphism.to.harmonic}, it implies that
the inequality $r_X\leq r_\Gamma $ is not sharp.

\begin{eg}[Luo]\label{ex:luo}
Let $\Gamma$ be a metric graph 
in $\mathcal M(G_7)$, 
where $G_7$ is the graph 
of genus seven depicted in Figure \ref{Luoye}(a), 
such that all edge lengths in $\Gamma$ are equal, and let
$D=(p)+(q)+(s)\in\Div(\Gamma)$. Then
$r_{\Gamma}(D)=1$, however there does not exist any finite harmonic morphism of metric
graphs $\phi:\Gamma' \to T$ of degree three to a metric tree for any
$\Gamma' \in \mathcal M(G_7)$. 
In particular, this shows that the stable gonality of an augmented
graph can be greater than its divisorial gonality.

\begin{figure}[h]
\begin{tabular}{ccc}
\scalebox{.32}{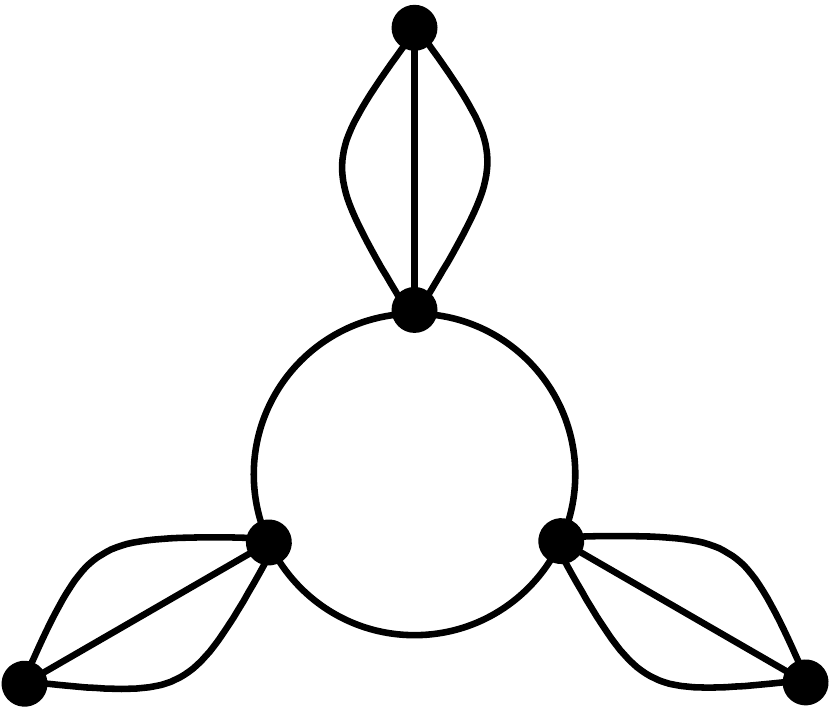}& \hspace{4ex} &
\scalebox{.32}{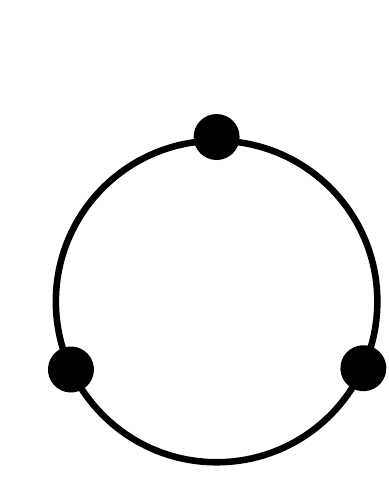}
\\ \\ a) Graph $G_7$ && b) Metric graph $\Gamma'_1 \subset \Gamma'$ 
\end{tabular}
\caption{}\label{Luoye}
\end{figure}

\medskip

We briefly sketch a proof.
Suppose that such a finite harmonic morphism $\phi:\Gamma' \to T$
exists.
Since $\Gamma'$ is not hyperelliptic,  
 one easily verifies that  $D_{\phi(p)}(\phi)=3(p)$, 
$D_{\phi(q)}(\phi)=3(q)$, and $D_{\phi(s )}(\phi)=3(s )$. This shows the existence of a finite
morphism  $\phi':\Gamma'_1\to T'$ of degree $3$ to a metric tree $T'$  where
$\Gamma'_1$ is 
depicted in Figure \ref{Luoye}(b),
so that  $D_{\phi'(p)}(\phi')=3(p)$, $D_{\phi'(q)}(\phi')=3(q)$, and $D_{\phi( s )}(\phi')=3(s )$.
But it is easy to verify by hand that such a morphism $\phi'$ does not
exist.
\end{eg}

\begin{prop}\label{lemma:notsharp}
Neither of the inequalties 
$r_X  \leq r_\cC$ and $r_\cC\leq r_\Gamma^{\#}$
is sharp.
\end{prop}

\pf
To show the non-sharpness of the inequality $r_X \leq r_\cC$,
let $\cC$ 
be a metrized complex of curves  whose underlying metric graph $\Gamma$ 
belongs to the family depicted in Figure \ref{hyper pic}, with first Betti number $\kappa$,
and whose genus function is positive at each vertex.
Consider the divisor $\mathcal D_d=d(p)\oplus d(x)$ in $\cC$ 
for a closed point $x$ in $C_p$ and $d$ a positive integer. 
If $d$ is sufficiently large compared to the genera of the vertices,
then $r_\cC (\mathcal D_d)\geq 1$. If the pair $(\cC, \mathcal D_d)$ lifted 
to a triple $(X,V; D_X)$ with $\tau_*(D_X) \sim \mathcal D_d$, then 
there would exist a finite harmonic morphism
$\phi:\widetilde \cC\to\mathcal T$ from a modification of
$\cC$ to a metrized complex of curves of genus zero.
But this would imply the existence of a degree $d$ morphism $\phi_p: C_p\to
\P^1$ such that the image of $\red_p$ (on edges adjacent to $p$ in $\Gamma$) is contained in the
set of critical values of $\phi_p$. By the Riemann--Hurwitz formula, this is impossible for 
$\kappa$ large enough compared to $d$.

\medskip

To show the non-sharpness of the inequality $r_\cC\leq r^\#$, let again $(\Gamma,g)$ be an augmented metric graph 
 with underlying graph depicted in Figure \ref{hyper pic} with $\kappa\ge 3$ and $2\leq 2g(p) <\kappa-2$,
 and let $D=2(p)$. One easily computes that $r_\Gamma^{\#}(D)=1$. An algebraic
 curve of genus $g(p)\ge 1$ contains at most $2g(p) +2$ distinct points $p$ such that $2(p)$ is in a 
given linear system of degree two, which implies that $(\Gamma,g)$ cannot be lifted to a hyperelliptic metrized
complex of curves.  This shows that the inequality $r_\cC \leq r_\Gamma^{\#}$ is not sharp.
\qed

\bibliographystyle{thesis}
\bibliography{harm_etale}
\bigskip~\bigskip

\immediate\closeout\exportaux

\end{document}